\documentclass[a4paper, 11pt, final]{amsart}

\makeatletter

\numberwithin{equation}{section}
\numberwithin{figure}{section}
\usepackage[a4paper]{geometry}
\usepackage{amsmath,amsfonts,amssymb,amsthm,epsfig}
\usepackage{amsmath}
\usepackage{color,xcolor}
\usepackage[final,linkcolor = blue,citecolor = blue,colorlinks=true]{hyperref}
\usepackage[leqno]{amsmath}
\usepackage[numbers,sort&compress]{natbib}
\usepackage{amscd}
\usepackage{mathrsfs}
\usepackage{}
\usepackage[all]{xy}
\usepackage[OT1]{fontenc}
\usepackage[latin1]{inputenc}
\usepackage{amssymb,latexsym}
\usepackage{type1cm,ae}
\usepackage[notref,notcite]{showkeys}
\usepackage{graphicx}
\usepackage{xspace}
\usepackage{setspace}
\usepackage{enumerate}
\usepackage[english]{babel}
\usepackage{bbm}
\usepackage{esint}
\usepackage{soul}

\usepackage{pifont}

\usepackage{graphicx}

\soulregister{\cite}7
\soulregister{\citep}7
\soulregister{\citet}7
\soulregister{\ref}7
\soulregister{\pageref}7
\sethlcolor{yellow}

\theoremstyle{definition}
\newtheorem{theorem}{Theorem}[section]
\newtheorem{proposition}[theorem]{Proposition}

\newtheorem{lemma}[theorem]{Lemma}

\newtheorem{corollary}[theorem]{Corollary}

\newtheorem{definition}[theorem]{Definition}
\newtheorem{remark}[theorem]{Remark}

\numberwithin{equation}{section}

\newcommand*{\Wert}{\mathord{\mbox{|\kern-1.5pt|\kern-1.5pt|}}}

\newcommand{\N}{\mathbb{N}}	
\newcommand{\Z}{\mathbb{Z}}	
\newcommand{\K}{\mathbb{K}}	
\newcommand{\R}{\mathbb{R}}	
\newcommand{\T}{\mathbb{T}}
\newcommand{\HH}{\mathbb{H}}	
\newcommand{\E}{\mathbb{E}}	
\newcommand{\A}{\mathbb{A}}	
\newcommand{\G}{\mathbb{G}}	
\renewcommand{\P}{\mathbb{P}}	
\newcommand{\V}{\mathscr{V}}	

\newcommand{\J}{\mathcal{J}}

\DeclareMathOperator{\supp}{supp}

\def\rr{{\Bbb R}}
\def\rz{{{\rr}^n}}

\def\nn{{\Bbb N}}
\def\cc{{\Bbb C}}

\def\ss{{\Bbb S}}

\def\B{{\bf B}}

\def\E{\mathbf{E}}

\def\L{\mathcal{L}}

\def\M{\mathcal{M}}

\def\G{\mathcal{G}}
\def\Ez{\mathcal{E}}
\def\Fz{\mathcal{F}}
\def\A{\mathcal{A}}
\def\K{\mathcal{K}}

\def\D{\mathcal{D}}

\def\P{\mathcal{P}}

\def\H{\mathcal{H}}
\def\T{\mathcal{T}}
\def\R{\mathcal{R}}

\def\S{\mathcal{S}}

\def\J{\mathcal{J}}

\def\rdd{\rr_+^{n+1}}

\def\supp{{\rm{\ supp\ }}}

\def\ez{\epsilon}

\def\supp{{\rm supp}}

\def\l{\left}
\def\r{\right}

\def\XXint#1#2#3{{\setbox0=\hbox{$#1{#2#3}{\int}$}
		\vcenter{\hbox{$#2#3$}}\kern-.5\wd0}}

\makeatother

\usepackage{babel}

\begin{document}

\raggedbottom

\allowdisplaybreaks

\title[The Kato problem and extensions for higher-order weighted elliptic operators]{The Kato problem and extensions for degenerate elliptic operators of higher order in weighted spaces}
	
	\author[Guoming Zhang]{Guoming Zhang}
	
	\address[Guoming Zhang]{College of Mathematics and System Science, Shandong University of Science and Technology, Qingdao, 266590, Shandong, People's Republic of China}
	\email{ zhangguoming256@163.com}

	

	\thanks{$^*$Corresponding author: Guoming Zhang}
	\thanks{The author is supported by the National Natural Science Foundation of Shandong Province (No. ~ZR2023QA124).}
	
	\date{}
	
	\begin{abstract} 
	
	We consider the Kato problem and extensions for degenerate elliptic operators of arbitrary order $2m$ ($m\geq 1$), shaped like$$\L_{w}:=(-1)^{m}\sum_{|\alpha|=|\beta|=m}w^{-1}\partial^{\alpha}(a_{\alpha, \beta}\partial^{\beta}),$$ whose coefficients $\{a_{\alpha, \beta}\}_{|\alpha|=|\beta|=m}$ are measurable, complex-valued and satisfy the G$\mathring{a}$rding inequality with respect to a Muckenhoupt $A_{2}$-weight; this generalizes the work of [Cruz-Uribe, Martell and Rios 2018]. 
	
	More precisely, the author identifies intervals that contain the exponents $p$ for which the relations $$\|\L_{w}^{1/2} f\|_{L^{p}(w)}\approx \|\nabla^{m} f\|_{L^{p}(w)}\quad \mbox{ and } \quad \|\L_{w}^{1/2} f\|_{L^{p}(vdw)}\approx \|\nabla^{m} f\|_{L^{p}(vdw)}$$ hold, given some suitable weight $v.$ Moreover, under some extra conditions on $w$ that allow us to take $v=w^{-1},$ the unweighted $L^{p}$-Kato estimate is obtained for $p$ close to $2$. In particular, if $w$ is a power weight $w_{\alpha}:=|x|^{\alpha}$, we prove that there exists $\ez>0,$ depending only on $n, m$ and the ellipticity constants, such that $$\|\L_{w_{\alpha}}^{1/2} f\|_{L^{2}(\rz)}\approx \|\nabla^{m} f\|_{L^{2}(\rz)},\quad \forall\; -\ez<\alpha< \frac{2 mn}{n+2 m}.$$
	
As an application, the unweighted $L^{p}$-Dirichlet, regularity and Neumann boundary value problems associated to $\L_{w}$ are solved when $p$ is sufficiently close to $2.$	
	\end{abstract}

	\maketitle
\section{Introduction}

We study the degenerate operators of order $2m,$ \begin{equation}\label{eq: z00}\L_{w}:=(-1)^{m}\sum_{|\alpha|=|\beta|=m}w^{-1}\partial^{\alpha}(a_{\alpha, \beta}\partial^{\beta}),\end{equation} where $w$ belongs to the Muckenhoupt class $A_{2}=A_{2}(\rz, dx).$ The coefficients $\{a_{\alpha, \beta}\}_{|\alpha|=|\beta|=m}$ are complex-valued and measurable, also satisfying the G$\mathring{a}$rding inequality: \begin{equation}\label{eq: z1} \mbox{Re} \int_{\rz}a_{\alpha, \beta}(x)\partial^{\alpha}f(x)  \overline{\partial^{\beta}f(x)}dx \geq c_{1} \int_{\rz}|\nabla^{m}f(x)|^{2}w(x)dx, \quad \forall f\in H^{m}(w),\end{equation} and \begin{equation}\label{eq: z2}\bigg|\sum_{|\alpha|=|\beta|=m} a_{\alpha, \beta}(x)\xi_{\alpha }\overline{\zeta_{\beta} }\bigg|\leq c_{2}w(x)|\xi| |\zeta|, \quad\forall \; (\xi_{\alpha })_{|\alpha|=m}, \;(\zeta_{\beta})_{|\beta|=m}\in (\cc)^{m},\end{equation} for some positive constants $c_{1}, c_{2}$ and all $x\in \rz.$ In what follows, we use $\Ez(w, c_{1}, c_{2})$ to denote the class of coefficients $\{a_{\alpha, \beta}(x)\}_{|\alpha|=|\beta|=m}$ of complex-valued and measurable functions verifying \eqref{eq: z1}-\eqref{eq: z2}.

The operator defined in \eqref{eq: z00} occurs as a natural higher-order extension of the second-order degenerate operator $L_{w}=-w^{-1}\mbox{div}(A\nabla),$ where $A$ is a real, symmetric and elliptic matrix controlled by a Muckenhoupt $A_{2}$-weight; the operator $L_{w}$ was pioneered in \cite{FJK, FJK1, FJK2}.  For $w\equiv 1,$ the operator $L_{w}$ simplifies to $L,$ a uniformly divergence-form elliptic operator $L.$ Notably, the Kato conjecture for $L-$a long-standing problem asserting that $L^{1/2}f$ is comparable to $\nabla f$ in $L^{2}(\rz)$ for all $f\in H^{1}(\rz)-$was settled in the remarkable paper \cite{AHLT} by Auscher, $et\; al.$ When $w$ is an $A_{2}$-weight, Cruz-Uribe and Rios \cite{CR, CR1, CR2} extended the techniques introduced in \cite{AHLT} to the weighted setting, thereby solving the Kato problem for $L_{w}$ and establishing the comparability between $L_{w}^{1/2}f$ and $\nabla f$ in $L^{2}(w).$ For further results on the Kato estimates for $L$ and $L_{w},$ one can refer to \cite{A, A1, AT, DMR, YZ}. Regarding the higher-order elliptic operators $\L$ (corresponding to $w\equiv 1$ in $\L_{w}$), Auscher, $et\; al.$ \cite{AHMT1} showed that $\L^{1/2}f $ is comparable to $\nabla^{m} f$ in $L^{2}(\rz)$ for all $f\in H^{m}(\rz).$ This result was subsequently generalized by the author \cite{GZ} to the higher-order degenerate operators $\L_{w}$ in the class $\Ez(w, c_{1}, c_{2});$ the autor proved that for all $f\in H^{m}(w),$ \begin{equation}\label{eq: cd1081} \|\L_{w}^{1/2} f\|_{L^{2}(w)}\approx \|\nabla^{m} f\|_{L^{2}(w)}. \end{equation} The estimate \eqref{eq: cd1081} acts as a starting point for our analysis, as the proof strategy outlined below aligns with the approach in \cite{DMR}.  

A central goal of this paper is to identify the conditions on the weight $w$ under which the square root $\L_{w}^{1/2}$ satisfies the unweighted $L^{p}(\rz)$ estimate \begin{equation}\label{eq: cd2000}\|\L_{w}^{1/2} f\|_{L^{p}(\rz)}\approx \|\nabla^{m} f\|_{L^{p}(\rz)}\quad \mbox{for $p$ near 2.}\end{equation} The entire proof can be roughly divided into three parts. First, we determine the ranges of $p$ and the conditions on weights $w, v$ that guarantee the weighted $L^{p}$-boundedness of $\L_{w}^{1/2}$: \begin{equation}\label{eq: cd2001}\|\L_{w}^{1/2} f\|_{L^{p}(w)}\lesssim \|\nabla^{m} f\|_{L^{p}(w)}\quad \mbox{ and } \quad \|\L_{w}^{1/2} f\|_{L^{p}(vdw)}\lesssim \|\nabla^{m} f\|_{L^{p}(vdw)}.\end{equation} Second, we derive the weighted norm estimates for the Riesz transform $\nabla^{m}\L_{w}^{-1/2},$ which corresponds to the reverse direction of the inequalities in \eqref{eq: cd2001}. Importantly, in the course of establishing these results, we also obtain the $L^{p}(w)$ and $L^{p}(vdw)$ estimates for the semigroup $e^{-t\L_{w}}$, its gradient $t^{1/2}\nabla^{m}e^{-t\L_{w}}$ and the functional calculus $\phi(\L_{w}$ with $\phi \in H^{\infty}(\Sigma_{\mu}),$ $\mu\in (\V, \pi)$) associated with the higher-order degenerate operator $\L_{w};$ as a consequence of these results, the weighted estimates for the square functions $g_{\L_{w}}$ and $G_{\L_{w}}$ follow. Third, explicit requirements on weights $w$ are specified for $p$ near 2 (whereas in the second-order case \cite{DMR} they are explicitly stated only for $p=2$), and these conditions permit setting $v=w^{-1}$ to derive \eqref{eq: cd2000}. In particular, when $p=2,$ the following theorem holds for the higher-order degenerate elliptic operators $\L_{w};$, It is a special case of Theorem \ref{theorem: cd1058} and generalizes \cite[Theorem 1.2]{DMR}.


\begin{theorem}\label{theorem: cd2002}\; Let $\L_{w}$ be given by \eqref{eq: z00} with $\{a_{\alpha, \beta}(x)\}_{|\alpha|=|\beta|=m}\in \Ez(w, c_{1}, c_{2}).$ Then, if $w\in A_{1}\cap \mbox{RH}_{1+\frac{n}{2m}},$ we have for every $f\in H^{m}(\rz)$ that \begin{equation}\label{eq: cd2003} \|L_{w}^{1/2} f\|_{L^{2}(\rz)}\approx \| \nabla^{m}f\|_{L^{2}(\rz)},\end{equation} where the implicit constants depend only on $n, m, c_{1}, c_{2}$ and the $A_{1}$ and $ \mbox{RH}_{1+\frac{n}{2m}}$ constants of $w$ (see Section 2.1 for the rigorous definitions of these weight classes). 

In particular, for the power weight $w_{-\alpha}(x):=|x|^{-\alpha},$ there exists a $\ez=\ez(n, m, c_{1}, c_{2})$ with $0<\ez<\frac{1}{2n}$ such that 
$$\|L_{w_{-\alpha}}^{1/2} f\|_{L^{2}(\rz)}\approx \| \nabla^{m}f\|_{L^{2}(\rz)},$$ provided $-\ez<\alpha< \frac{2 mn}{n+2 m}.$
\end{theorem}

This work also provides a solution to the $L^{p}(\rz)$-regularity problem for $p$ close to 2 on $\rdd:=\rz\times [0, \infty):$ \begin{equation}\label{eq: cd2004} \l\{
\begin{aligned}
&\; \partial_{t}^{2}u-\L_{w}u=0\quad\quad\quad\quad\quad\quad\quad\;\mbox{on}\;\rz, \quad\quad\quad\quad\quad\quad\\
&\;\nabla^{l}u(\cdot, t)|_{\partial \rdd}=\nabla^{l} f\;\quad\;\;\;\;\quad\quad\mbox{on}\;\partial \rdd=\rz, \; 0 \leq l\leq m-1,\\
& \sup_{t>0}\l(\|t^{k-1}\partial_{t}^{k}u(\cdot, t)\|_{L^{p}}+\|\nabla^{l}u(\cdot, t)\|_{L^{p}}\r) \lesssim \|f\|_{H^{m, p}}, \; 0 \leq l\leq m, \end{aligned}\r.\end{equation} as a direct application of the unweighted $L^{p}$ estimates for the semigroup, the Riesz transform and the functional calculus of $\L_{w}.$ The corresponding Dirichlet and Neumann problems are also considered, with similar results holding. In fact, building on the main results and proofs of this paper and inspired by \cite{C11, C12, ACMP},  it is natural to expect that the uniform $L^{p}(\rz)$-norm estimates for higher-order derivatives in \eqref{eq: cd2004} can be improved to non-tangential maximal function estimates. This problem will be solved in our next work. 
Another direction for future research is the generalization of the Dirichlet and Neumann problems in \cite{ARR, B2} to higher-order degenerate elliptic operators $\L_{w}.$ 


A significant technical obstruction in our proof, as well as in \cite{DMR}, is that the weight $w$ is only assumed a priori to be in $A_{2}.$ even though this implies the existence of a small $\ez>0$ such that $w\in A_{2-\ez}.$ It prevents us from completely characterizing the interval $\K(\L_{w}),$ which consists of the pairs $(p, q)$ for which $t^{\frac{1}{2m}}\nabla^{m}e^{-t\L_{w}}\in \mho(L^{p}(w)\to L^{q}(w)); $ see Section 6.1. It also prevents a version of Lemma \ref{lemma: cd105} in Section 6.1 with $t^{\frac{k}{2m}}\nabla^{k}e^{-t\L_{w}}$ ($1\leq k\leq m-1$) in place of $t^{\frac{1}{2m}}\nabla^{m}e^{-t\L_{w}},$ marking a key difference from the unweighted case. Fortunately, for our proof,  it is not necessary to handle these intermediate families $\{t^{\frac{k}{2m}}\nabla^{k}e^{-t\L_{w}}\}_{t>0}$ for $1\leq k\leq m-1.$  Compared to the second-order case in \cite{DMR}, the generalized Poincar\'{e}-Sobolev inequalities in Theorem \ref{theorem: cd0} (with $P_{B} f$ replacing $\fint_{B}fdw$) pose a challenge. To overcome this, we introduce, motivated by \cite{Chua}, a refined projection $\pi_{B}^{m}f$ defined by \eqref{eq: cd72}. The projection also avoids the telescoping argument used in \cite{DMR} to treat the integral average $\fint_{B}fdw$, as seen in the proof of \eqref{eq: cd1034}. Throughout our argument, and in contrast to the approaches in \cite{D0, D1, A} for higher-order elliptic operators in the unweighted setting, dividing the proof into two cases $n\geq 2m$ and $n<2m$ is not required. 



The plan of the paper is as follows. In Section 2, we introduce the definitions and properties of Muckenhoupt weights (including power weights $|x|^{\alpha}$), along with the associated higher-order Sobolev spaces and the generalized Poincar\'{e}-Sobolev inequalities defined on them. We also define higher-order off-diagonal estimates (a generalization of the concept in \cite{A, AM, AM2, DMR}), with the section ending by listing key lemmas on off-diagonal estimates and two core theorems (Theorem \ref{theorem: cd18}, Theorem \ref{theorem: cd21}) for our proof. In Section 3, the $H^{\infty}$ functional calculus of $\L_{w}$ in $L^{2}(w)$ is used to build the $L^{2}(w)$-off-diagonal estimates for $e^{-z\L_{w}}$ in the sector $\Sigma_{\frac{\pi}{2}-\V};$ based on these results, we further prove the $L^{p}(w)$- and $L^{p}(vdw)$-off-diagonal estimates for $e^{-t\L_{w}}.$ The $L^{p}(w)$ and $L^{p}(vdw)$ functional calculi for $\L_{w},$ contained in Section 4, form the basis for subsequent analysis. In Section 5 the reverse inequalities \eqref{eq: cd2001} are proved by synthesizing the main results from preceding sections. The proof additionally relies on two higher-order tools: a weighted Calder$\acute{o}$n-Zygmund decomposition and a weighted conservation property, constructed in Sections 5.1 and 5.2, respectively.

In Section 6, we show the existence of the interval $\K(\L_{w})$ and present its key properties, with a focus on showing that $2$ is an interior point of $\K(\L_{w}).$ To carry out the proof, we need a reverse H\"{o}lder inequality (with sharp constants) for solutions of the higher-order degenerate elliptic operator $\L_{w},$ whose proof is given in the Appendix. Section 7 is devoted to the $L^{p}(w)$- and $L^{p}(vdw)$-boundedness of the Riesz transform $\nabla^{m}\L_{w}^{-1/2},$ and Section 8 foucus on proving the $L^{p}(w)$- and $L^{p}(vdw)$-estimates for the vertical square functions $g_{\L_{w}}$ and $G_{\L_{w}},$ associated to the semigroups $e^{-t\L_{w}}$ and $t^{1/2}\nabla^{m}e^{-t\L_{w}}$, respectively. In Section 9, we characterize weight conditions on $w$ that permit the choice  $v\equiv w^{-1}.$ These conditions are then used to establish the unweighted $L^{p}$-boundedness (for $p$ near 2) of several key operators: the semigroup $e^{-t\L_{w}},$ its gradient $t^{1/2}\nabla^{m}e^{-t\L_{w}},$ the functional calculus $\phi(\L_{w}$ ($\phi \in H^{\infty}(\Sigma_{\mu}),$ $\mu\in (\V, \pi)$), the Riesz transform $\nabla^{m}\L_{w}^{-1/2},$ and the vertical square functions $g_{\L_{w}}$ and $G_{\L_{w}}.$ These boundedness results further enable us to solve the corresponding $L^{p}(\rz)$-Dirichlet, regularity and Neumann boundary value problems.

\section{Preliminaries}

We now rigorously define the notations used in the introduction, along with additional symbols needed to present our results.

\subsection{Weights and the associated Sobolev spaces.} Given a set $E\subset \rz (n\geq 2),$ we use the notation $$\fint_{E}h=\frac{1}{|E|}\int_{E}f(x)dx\quad\mbox{and}\quad \fint_{E}hdw=\frac{1}{w(E)}\int_{E}f(x)dw=\frac{1}{w(E)}\int_{E}f(x)w(x)dx$$ for unweighted and weighted averages, respectively. We say that a non-negative locally integrable function $w$ belongs to the Muckenhoupt class $A_{p},$ $1< p<\infty,$ if $$[w]_{p}:=\sup_{B\subset \rz}\l(\fint_{B}w\r)\l(\fint_{B}w^{1-p'}\r)^{p-1}\lesssim 1.$$ Hereafter, for two positive constants $A, B$, the expression $A\lesssim B$ means that there exists a nonessential constant $C,$ depending on $n, m$ and other parameters that will be clear in the text, such that $A\leq C B.$ The notations $A \gtrsim B$ and $A\approx B$ should be understood similarly. When $p=1,$ we say that $w\in A_{1}$ if $$[w]_{1}:=\sup_{B\subset \rz}\l(\fint_{B}w\r)\sup_{x\in B}w(x)^{-1}\lesssim 1.$$ The reverse H\"{o}lder classes are defined in the following way: $w\in \mbox{RH}_{q}, 1<q<\infty,$ if $$[w]_{\mbox{RH}_{q}}:=\sup_{B\subset \rz}\l(\fint_{B}w\r)\l(\fint_{B}w^{q}\r)^{\frac{1}{q}}\lesssim 1,$$ in particular, $w\in \mbox{RH}_{\infty}$ if $$[w]_{\mbox{RH}_{\infty}}:=\sup_{B\subset \rz}\l(\fint_{B}w\r)^{-1}\sup_{x\in B}w(x)\lesssim 1.$$ We also need the new weight class $A_{p}(w)$ and $\mbox{RH}_{q}(w),$ which are defined in \cite{DMR} by substituting Lebesgue measure in the above definitions with $dw=w(x)dx;$ e.g., $v\in A_{p}(w)$ if $$[w]_{A_{p}(w)}:=\sup_{B\subset \rz}\l(\fint_{B}vdw\r)\l(\fint_{B}v^{1-p'}dw\r)^{p-1}\lesssim 1.$$

We summarize some important properties of these classes in the following proposition for easy reference. 
\begin{proposition}\label{proposition: cd3000}\; (\cite{G1, G2})
\begin{equation*}\begin{aligned} &(i)\quad A_{1}\subset A_{p}\subset A_{q},\; \mbox{RH}_{\infty}\subset \mbox{RH}_{q}\subset \mbox{RH}_{p},\;\mbox{for}\; 1<p\leq q<\infty.\\[4pt] 
&(ii)\quad\mbox{If}\; w\in A_{p}, 1<p<\infty,\; \mbox{there exists}\; \ez>0\; \mbox{such that}\; w\in A_{p-\ez}, \;\mbox{similarly if}\; w\in \mbox{RH}_{q}, \\[4pt]\end{aligned} \end{equation*}\end{proposition}
\begin{equation*}\begin{aligned} &\;\quad\quad 1<q<\infty,\; \mbox{there exists}\; \delta>0\; \mbox{such that}\; w\in \mbox{RH}_{q+\delta}.  \\[4pt]
&(iii)\quad\mbox{Given}\;1<p<\infty\; \mbox{and}\; M>0,\; \mbox{there exist }\; C=C(n, p, M)\;\mbox{and}\;\delta=\delta(n, p, M)\\[4pt]
&\;\quad\quad\mbox{such that for all}\; w\in A_{p,}\; [w]_{p}\leq M\; \mbox{implies}\; [w]_{p-\delta}\leq C.\\[4pt]
&(v)\quad A_{\infty}=\cup_{1\leq p<\infty}A_{p}=\cup_{1<q\leq \infty}\mbox{RH}_{q}.\\[4pt]
&(vi)\quad \mbox{If}\; 1<p<\infty, \;\mbox{then}\; w\in A_{p}\; \mbox{if and only if }\;w^{1-p'}\in A_{p'}\\[4pt]
&(vii)\quad\mbox{If}\; w\in A_{p} (v\in A_{p}(u)), 1\leq p<\infty,\; \mbox{then}\; \forall\;\delta>0,\; w\in A_{q} (v\in A_{q}(u))\;\mbox{with}\; q=\delta p+1-\delta.\\[4pt]
&(viii)\quad\mbox{Let}\; w_{1}, w_{2}\in A_{1},\; \mbox{then}\; w_{1}w_{2}^{1-p}\in A_{p}\;\mbox{for any}\; 1<p<\infty.\\[4pt]
&(ix)\quad\mbox{If}\;1\leq q\leq \infty\; \mbox{and}\; 1\leq s<\infty, \; \mbox{then}\;w\in A_{q}\cap\mbox{RH}_{s}\; \mbox{if and only if}\; w^{s}\in A_{s(q-1)+1}. \\[4pt]
&(x)\quad w^{-1}\in A_{p}(w)\; \mbox{if and only if }\;w\in \mbox{RH}_{p'},\;\mbox{and}\; w^{-1}\in  \mbox{RH}_{s}(w)\; \mbox{if and only if }\;w\in A_{s'}.\\[4pt]
\end{aligned} \end{equation*}
It is particularly important to note that, given $w\in A_{p}$ with $1\leq p<\infty,$ there is a constant $D=D(p, n)$ (the doubling order of $w$) such that for any $\lambda \geq 1$ and any ball $B$
\begin{equation}\label{eq: cd66} 
w(\lambda B)\leq \lambda^{D} [w]_{p}  w(B).
\end{equation} As a consequence of \eqref{eq: cd66}, $(\rz, dw, |\cdot|)$ becomes a space of homogeneous type, where $|\cdot|$ denotes the usual Euclidean distance. A canonical example of Muckenhoupt weights is provided by the power weights $w_\alpha(x):= |x|^\alpha$ with $\alpha>-n.$ It is well-known that $w_{\alpha} \in A_1$ for $-n < \alpha \leq 0$, and $w_{\alpha} \in A_p$ $(1 < p < \infty)$ if $-n < \alpha < n(p-1)$; moreover, $w_{\alpha} \in \mbox{RH}_{\infty}$ for $0 \leq \alpha < \infty$, while $w_{\alpha} \in \mbox{RH}_{q}$ holds if $-n/q < \alpha < \infty$. If we define $$r_{w}:=\inf\{p: w\in A_{p}\}\quad \mbox{and}\; s_{w}:=\sup\{q: w\in \mbox{RH}_{q}\},$$ then  \begin{equation}\label{eq: cd1043}r_{w_{\alpha}}=\max\{1, 1+\frac{\alpha}{n}\}, \quad s_{w_{\alpha}}=(\max\{1, (1+\frac{\alpha}{n})^{-1}\} )'.\end{equation}




We will use symbols such as $\alpha, \beta, \gamma$ to denote multi-indices in $(\nn)^{n}.$ (Here, $\nn$ deonotes the non-negative integers.) If $\alpha=(\alpha_{1}, ..., \alpha_{n})$ is a multi-index and $k\in \nn,$ we define $|\alpha|=\alpha_{1}+...+\alpha_{n},$ $\partial^{\alpha}=\partial_{x_{1}}^{\alpha_{1}}\partial_{x_{2}}^{\alpha_{2}}\cdot\cdot\cdot\partial_{x_{n}}^{\alpha_{n}},$ and $\nabla^{k}=(\partial^{\alpha})_{|\alpha|=k}.$ In particular, we introduce the notation $\mbox{div}_{m}:=\sum_{|\alpha|=m}\partial^{\alpha}.$ We also let $(\cc)^{m}:=\{\xi: \xi=(\xi_{\alpha})_{|\alpha|=m}, \;\xi_{\alpha}\in \cc\},$ and for any $\xi, \zeta\in (\cc)^{m},$ $\xi \cdot \overline{\zeta}:=\sum_{|\alpha|=m}\xi_{\alpha} \overline{\zeta_{\alpha}}$ denote the inner product on $(\cc)^{m}.$

Given $w\in A_{2},$ let $\Omega \subset \rz$ be a domain. We denote by $H^{m}(\Omega, w):=W^{m, 2}(\Omega, dw)$ the weighted Sobolev spaces of order $m,$ consisting of distributions for which all $\partial^{\alpha} f$ ($|\alpha|\leq m$) belong to $L^{2}(\Omega, w).$ When $\Omega=\rz,$ we simply write $H^{m}(w)=H^{m}(\rz, w)$ and $L^{2}(w)=L^{2}(\rz, w).$ This space $H^{m}(w)$ is a Hilbert space and coincides with the space defined as the completion of $C_{c}^{\infty}(\rz)$ with respect to the norm $$\|f\|_{H^{m}( w)}:=(\sum_{|\alpha|\leq m}\|\partial^{\alpha} f\|_{L^{2}(w)}^{2})^{1/2};$$ see \cite{NM}. Similarly, we can define $W^{m, p}(\rz, dw)$ ($1\leq p<\infty$) when $w\in A_{p}$  and the unweighted space $W^{m, p}(\rz)$ by taking $w\equiv 1.$ In particular, from the weighted Sobolev interpolation inequality in \cite{GW}: \begin{equation}\label{eq: a2.32}\l(\int_{\rr^{n}}|\partial^{\gamma}v|^{p}w\r)\lesssim \l(\int_{\rr^{n}}|v|^{p}w\r)^{(1-\frac{|\gamma|}{m})}\l(\int_{\rr^{n}}|\nabla^{m}v|^{p}w\r)^{\frac{|\gamma|}{m}} \quad (\forall\; |\gamma|\leq m),\end{equation} it follows that for any $f\in W^{m, p}(\rz, dw),$\begin{equation*}
\|f\|_{W^{m, p}_{w}}\approx (\|f\|_{L^{p}(w)}^{p}+\|\nabla^{m} f\|_{L^{p}(w)}^{p})^{1/p}.\end{equation*}

\subsection{Generalized Poincar\'{e}-Sobolev inequalities.} Repeated application of \cite[Theorem 2.1]{DMR} yields the following weighted Poincar\'{e}-Sobolev inequalities of higher order. 
\begin{theorem}\label{theorem: cd0} Assume $w\in A_{p}$ with $p\geq 1.$  Then, for any $f\in C_{0}^{\infty}(B)$ and any $p\leq q<p^{*, m}_{w},$ \begin{equation}\label{eq: cd35} 
\l(\fint_{B}|f|^{q}dw\r)^{\frac{1}{q}}\lesssim r(B)^{m}\l(\fint_{B}|\nabla^{m}f|^{p}dw\r)^{\frac{1}{p}},
\end{equation} where $\frac{1}{p^{*, m}_{w}}:=\frac{1}{p}-\frac{m}{nr_{w}}$ if $p<\frac{n r_{w}}{m},$ and $p^{*, m}_{w}=\infty$ otherwise. Moreover, if $f\in C^{\infty}(B),$ there exists a polynomial $Q_{B} f$ of degree at most $m-1$ such that \begin{equation}\label{eq: cd1045} \int_{B}D^{\beta}(f-Q_{B} f)dw=0,\quad \forall\; |\beta|\leq m-1,\end{equation} and \begin{equation}\label{eq: cd1044} 
\l(\fint_{B}|f-Q_{B}f|^{q}dw\r)^{\frac{1}{q}}\lesssim r(B)^{m}\l(\fint_{B}|\nabla^{m}f|^{p}dw\r)^{\frac{1}{p}}
\end{equation} for any $p\leq q<p^{*, m}_{w},$ where the implicit constants depend only on $n, m, p$ and the weight constants. \end{theorem} 

 
\begin{remark}\label{remark: cd1082} Defining the projection of a function $u$ onto $\P_{m-1}$ (the collection of polynomials with degree at most $m-1$) by solely requiring \eqref {eq: cd1045} may not always meet our needs. To address this, we introduce a more refined projection denoted by $\pi_{Q}^{m},$ which has an explicit formula given by \eqref{eq: cd72}; this formula plays a crucial role in our proof (see \cite{Chua}). 

Set $\E_{m, w}^{p}:=\{u\in \D'(\rz):  \|\nabla^{m} u\|_{L^{p}_{w}(\rz)}<\infty\},$ and define a projection $\pi_{Q}^{m}: \E_{m, w}^{p}\to \P_{m-1}$ by \begin{equation}\label{eq: cd72}\pi_{Q}^{m}(u)(x)=r^{-n}\sum_{|\beta|\leq m-1}\l(\frac{x-z}{r}\r)^{\beta}\int_{B_{r}(0)}\phi_{\beta}(y/r)u(y+z)dy,\end{equation} where $B_{r}(z)$ is the largest ball \footnote{Indeed, it suffices to require $B_{r}(z) \subset Q$ with $r\approx l(Q)$ and $z$ coinciding with the center of $Q$ (hence $Q$ is starshaped with respect to $B_{r}(z)$); see \cite[Theorem 1.1.10]{Ma}. Of course, the cube $Q$ can be replaced by a ball $B.$} in $Q$  and \begin{equation*}\phi_{\beta}(y)=\sum_{0\leq |\gamma|\leq m-1-|\beta|}\frac{(n+m-1)!}{(n+|\gamma+\beta|)!(m-1-|\gamma+\beta|)!}(-1)^{|\beta|}\frac{1}{\beta!\gamma!}y^{\gamma}D^{\beta+\gamma}v(y)\end{equation*} with $v\in C_{0}^{\infty}(B_{1}(0))$ and $\int v=1.$ It is clear that $\pi_{Q}^{m}u=u$ if $u\in \P_{m-1}.$ Following the argument in \cite[Theorem 4.5; Lemma 4.6]{Chua}, for any $|\gamma|\leq m-1,$ it holds that \begin{equation}\label{eq: cd43}\|D^{\gamma}\pi_{Q}^{m} u\|_{L^{\infty}(Q)}\lesssim r^{-n}\int_{Q}|D^{\gamma} u|\lesssim \l(\fint_{Q}|D^{\gamma} u|^{p}dw\r)^{1/p},\end{equation} while \cite[Theorem 4.7]{Chua} implies \begin{equation}\label{eq: cd44}\|D^{\gamma}(\pi_{Q}^{m} u-u)\|_{L_{w}^{p}(Q)}\lesssim r^{m-|\gamma|}\|\nabla^{m}u\|_{L_{w}^{p}(Q)}.\end{equation}  \end{remark}
 
To ensure the validity of our proof, we further need a Poincar\'{e}-Sobolev inequality featuring a sharp constant estimate. This inequality should be compared to the counterpart in \cite[Remark 2.5]{DMR} for the second-order case. 
 
\begin{theorem}\label{theorem: cd8}\;(\cite[Corollary 2.7]{CMPR}) Assume $ 1\leq p<n$ and $w\in A_{q}$ with $1\leq q\leq p.$ Let $\frac{1}{p^{*}_{w}}:=\frac{1}{p}-\frac{m}{n (q+\log [w]_{A_{q}} )}$ if $p<\frac{n (q+\log [w]_{A_{q}} )}{m},$ and $p^{*}_{w}=\infty$ otherwise, Then for every ball $B$ and $f\in C^{\infty}(B),$ there exists a polynomial $P_{B} f$ of degree at most $m-1$  such that for any $s<p^{*}_{w},$
\begin{equation}\label{eq: cd7}\l(\fint_{B}|f-P_{B}f|^{s}dw\r)^{\frac{1}{s}}\lesssim [w]^{\frac{1}{p}}_{A_{q}}r(B)^{m}\l(\fint_{B}|\nabla^{m}f|^{p}dw\r)^{\frac{1}{p}}.\end{equation}
\end{theorem} 
Set $\tau_{w}:=\inf\{l+\log [w]_{A_{l}}: r_{w}<l\leq q\};$ this value equals $r_{w}$ when $m=1$. As a direct corollary of Theorem \ref{theorem: cd8}, we have 
\begin{corollary}\label{corollary: cd1}\;  Assume $1\leq p<n$ and $w\in A_{q}$ with $1\leq q\leq p.$ Let $\frac{1}{\widetilde{p^{*}_{w}}}:=\frac{1}{p}-\frac{m}{n \tau_{w}}$ if $p<\frac{n \tau_{w}}{m},$ and $\widetilde{p^{*}_{w}}=\infty$ otherwise. Then, for every ball $B$ and any $f\in C^{\infty}(B),$ there exists a polynomial $P_{B} f$ of degree at most $m-1$ such that for any $s<\widetilde{p^{*}_{w}}$, we can find a $q^{\star}$ such that $r_{w}<q^{\star}\leq q$ and \begin{equation}\label{eq: cd6}\l(\fint_{B}|f-P_{B}f|^{s}dw\r)^{\frac{1}{s}}\lesssim [w]^{\frac{1}{p}}_{A_{q^{\star}}}r(B)^{m}\l(\fint_{B}|\nabla^{m}f|^{p}dw\r)^{\frac{1}{p}}.\end{equation}
\end{corollary} 

\subsection{Off-diagonal estimates in higher-order setting.} We now define the higher-order off-diagonal estimates and full off-diagonal estimates on balls, which are the corresponding generalization of those in \cite[Definition 2.1, Definition 3.1]{AM} or \cite[Definition 2.23, Definition 2.33]{DMR}.

For a fixed ball $B,$ we set $C_{j}(B)=2^{j+1}B\setminus 2^{j} B$ for $j\geq 2;$ $C_{1}(B)=4B.$ Since $w(2^{j+1} B)\approx w(C_{j}(B))$ for $w\in A_{2}$ and \eqref{eq: cd66}, we may, by a slight abuse of notation, write $$\fint_{C_{j}(B)}hdw=\frac{1}{w(2^{j+1} B)}\int_{C_{j}(B)}hdw.$$ 

\begin{definition}\label{definition: cd24}\; Given $1\leq p\leq q\leq \infty,$ a family $\{T_{t}\}_{t>0}$ of sublinear operators satisfies $L^{p}(w)\to L^{q}(w)$ off-diagonal estimates on balls, denoted by $$T_{t}\in \mho (L^{p}(w)\to L^{q}(w)),$$ if there exist constants $\theta_{1}, \theta_{2}>0$ and $c>0$ such that for every $t>0$ and for any ball $B,$ setting $r=r(B)$ and $\Upsilon(s):=\max\{s, s^{-1}\}$ for $s>0,$ \begin{equation}\label{eq: cd25}\l(\fint_{B}|T_{t}(f1_{B})|^{q}dw\r)^{\frac{1}{q}}\lesssim \Upsilon\l(\frac{r}{t^{1/2m}}\r)^{\theta_{2}} \l(\fint_{B}|f|^{p}dw\r)^{\frac{1}{p}},\end{equation} and for all $j\geq 2,$
\begin{equation}\label{eq: cd26}\l(\fint_{B}|T_{t}(f1_{C_{j}(B)})|^{q}dw\r)^{\frac{1}{q}}\lesssim 2^{j\theta_{1}}\Upsilon\l(\frac{2^{j}r}{t^{1/2m}}\r)^{\theta_{2}} e^{-c\l(\frac{2^{j}r}{t^{\frac{1}{2m}}}\r)^{\frac{2m}{2m-1}}}\l(\fint_{C_{j}(B)}|f|^{p}dw\r)^{\frac{1}{p}},\end{equation}
\begin{equation}\label{eq: cd27}\l(\fint_{C_{j}(B)}|T_{t}(f1_{B})|^{q}dw\r)^{\frac{1}{q}}\lesssim 2^{j\theta_{1}}\Upsilon\l(\frac{2^{j}r}{t^{1/2m}}\r)^{\theta_{2}} e^{-c\l(\frac{2^{j}r}{t^{\frac{1}{2m}}}\r)^{\frac{2m}{2m-1}}}\l(\fint_{B}|f|^{p}dw\r)^{\frac{1}{p}}.\end{equation}
If the family of sublinear operators $\{T_{z}\}_{z\in \Sigma_{\mu}}$ is defined on a complex sector $\Sigma_{\mu}:=\{z\in \cc: z\neq 0, |\mbox{arg}z|<\mu\}$ ($\mu>0$), we say that it satisfies $L^{p}(w)\to L^{q}(w)$ off-diagonal estimates on balls in $\Sigma_{\mu}$ if \eqref{eq: cd25}-\eqref{eq: cd27} hold for $z\in \Sigma_{\mu}$ with $t$ replaced by $|z|$ in the right-hand side terms. We denote this by $T_{z}\in \mho (L^{p}(w)\to L^{q}(w), \Sigma_{\mu}).$\end{definition}

\begin{definition}\label{definition: cd1083}\; Given $1\leq p\leq q\leq \infty,$ a family of operators $\{T_{t}\}$ satisfies full off-diagonal estimates from $L^{p}(w)$ to $ L^{q}(w),$ denoted by $T_{t}\in \Fz(L^{p}(w)\to L^{q}(w)),$ if there exist constants $c, C>0$ such that for any closed sets $E, F,$ $$\|T_{t}(f1_{E})1_{F}\|_{L^{q}(w)}\leq C t^{-\frac{1}{2m}(\frac{n}{p}-\frac{n}{q})}e^{-c\l(\frac{d(E, F)}{t^{\frac{1}{2m}}}\r)^{\frac{2m}{2m-1}}}\|f1_{E}\|_{L^{p}(w)}.$$
\end{definition}

The results presented below are higher-order generalizations of those in \cite{AM, DMR} and serve as our primary analytical tools. The proofs follow the methodology developed in \cite{AM} and share essential features with the original arguments. As the extension procedure does not pose real difficulties, the detailed proof are omitted.

\begin{lemma}\label{lemma: h3}(\cite[Lemma 2.27]{DMR})\; Given $1\leq p_{i}\leq q_{i}\leq \infty,$ $i=1, 2.$ Assume that $T_{t}\in \mho(L^{p_{1}}(w)\to L^{q_{1}}(w))$ and $T_{t}: L^{p_{2}}(w)\to L^{q_{2}}(w)$ is uniformly bounded. Then $T_{t}\in \mho(L^{p_{\theta}}(w)\to L^{q_{\theta}}(w)),$ $0<\theta<1,$ where $$\frac{1}{p_{\theta}}=\frac{\theta}{p_{1}}+\frac{1-\theta}{p_{2}},\; \frac{1}{q_{\theta}}=\frac{\theta}{q_{1}}+\frac{1-\theta}{q_{2}}.$$
 \end{lemma}

\begin{lemma}\label{lemma: h4}(\cite[Lemma 2.28]{DMR})\; If $1\leq p\leq p_{1}\leq q_{1}\leq q\leq \infty,$ then $$\mho(L^{p_{1}}(w)\to L^{q_{1}}(w))\subset \mho(L^{p}(w)\to L^{q}(w)).$$
 \end{lemma}

\begin{lemma}\label{lemma: h5}(\cite[Lemma 2.29]{DMR})\; Suppose that $\{T_{t}\}_{t>0}$ are a family linear operators and $T_{t}\in \mho(L^{p}(w)\to L^{q}(w))$ with $1\leq p\leq q\leq \infty.$ Then $T_{t}^{*}\in \mho(L^{q'}(w)\to L^{p'}(w)),$ where $T_{t}^{*}$ is the dual operator of $T_{t}$ for the inner product $\int_{\rz}f\overline{g}dw.$
 \end{lemma}

\begin{lemma}\label{lemma: h6}\;(\cite[Theorem 2.3, Theorem 4.3]{AM}) \begin{equation*}\begin{aligned}&(i)\quad \mbox{If}\; T_{z}\in \mho (L^{p}(w)\to L^{p}(w), \Sigma_{\mu}), 0\leq \mu<\pi, 1\leq p\leq \infty, \;\mbox{then}\; T_{z}: (L^{p}(w)\to L^{p}(w) \\[4pt] 
&\quad\quad \mbox{is uniformly bounded on}\; \Sigma_{\mu};\\[4pt] 
&(ii)\quad \mbox{If}\; 1\leq p\leq q\leq r\leq \infty, T_{z}\in \mho (L^{q}(w)\to L^{r}(w), \Sigma_{\mu}) \mbox{and}\; S_{z}\in \mho (L^{p}(w)\to L^{q}(w), \Sigma_{\mu}),\\[4pt] 
&\quad\quad \mbox{then}\; T_{z}\circ S_{z}\in \mho (L^{p}(w)\to L^{r}(w), \Sigma_{\mu}).\end{aligned}\end{equation*} 
\end{lemma}

\begin{lemma}\label{lemma: h7}\;(\cite[Proposition 3.2]{AM})\; Given $1\leq p\leq q\leq \infty.$\begin{equation*}\begin{aligned}&(i)\quad \mbox{If}\; T_{t}\in \Fz (L^{p}(w)\to L^{q}(w)), \; \mbox{then}\; T_{t}: (L^{p}(w)\to L^{q}(w)\;\mbox{is uniformly bounded};\\[4pt] 
&(ii)\quad T_{t}\in \mho (L^{p}(w)\to L^{p}(w)), \mbox{if and only if}\; T_{t}\in \Fz (L^{p}(w)\to L^{p}(w)).\\[4pt] 
\end{aligned}\end{equation*} 
\end{lemma}

\begin{proposition}\label{proposition: h8}\;(\cite[Section 6.5]{AM})\; Let $1\leq p_{0}<q_{0} \leq \infty$ and $T_{t}\in \mho (L^{p}(w)\to L^{q}(w))$ for all $p, q$ with $p_{0}<p\leq q<q_{0}.$ Then, for all $p, q$ with $p_{0}<p\leq q<q_{0}$ and for any $\in A_{\frac{p}{p_{0}}}(w)\cap \mbox{RH}_{(\frac{p_{0}}{q})'}(w),$ we have $T_{t}\in \mho (L^{p}(vdw)\to L^{q}(vdw)).$
\end{proposition}

\begin{lemma}\label{lemma: h1}\;(\cite[Lemma 6.6]{AM})\; If $T_{t}\in \mho(L^{p}(w)\to L^{q}(w))$ with parameters $\theta_{1}, \theta_{2},$ then there exist $\theta_{1}', \theta_{2}'$ such that for any $0<c'<c,$ any ball $B$ with radius $r$ and for every $j\geq 1,$ $$\l(\fint_{B}|T_{t}(f1_{(2^{j} B)^{c}})|^{q}dw\r)^{\frac{1}{q}}\lesssim 2^{j\theta_{1}'}\Upsilon\l(\frac{2^{j}r}{t^{1/2m}}\r)^{\theta_{2}'} e^{-c'\l(\frac{2^{j}r}{t^{\frac{1}{2m}}}\r)^{\frac{2m}{2m-1}}}\l(\fint_{C_{j}(B)}|f|^{p}dw\r)^{\frac{1}{p}}$$ and 
$$\l(\fint_{(2^{j} B)^{c}}|T_{t}(f1_{B})|^{q}dw\r)^{\frac{1}{q}}\lesssim 2^{j\theta_{1}'}\Upsilon\l(\frac{2^{j}r}{t^{1/2m}}\r)^{\theta_{2}'} e^{-c'\l(\frac{2^{j}r}{t^{\frac{1}{2m}}}\r)^{\frac{2m}{2m-1}}}\l(\fint_{B}|f|^{p}dw\r)^{\frac{1}{p}}.$$ \end{lemma}

\begin{theorem}\label{theorem: h9}\;( \cite[Theorem 4.3]{AM})\; Let $1\leq p\leq p_{0}\leq q\leq  \infty$ and $\V_{1}$ with $0\leq \V_{1}<\V_{0}.$ Assume that $\{T_{t}\}_{t>0}\in \mho (L^{p}(w)\to L^{q}(w))$ and that $T_{z}\in \mho (L^{p}(w)\to L^{q}(w), \Sigma_{\V_{0}}).$ Then for any $l\in \nn,$ $z^{l}\frac{d^{l}T_{z}}{dz^{l}} \in \mho (L^{p}(w)\to L^{q}(w), \Sigma_{\V_{1}}).$
\end{theorem}

Indeed, the right hand side of the estimate \eqref{eq: cd25} in Definition \ref{definition: cd24} self-improves, as captured by the following lemma. 

\begin{lemma}\label{lemma: cd103}\;  Given $w\in A_{\infty}$ and a family of sublinear operators $\{T_{t}\}_{t>0}$ such that $T_{t}\in \mho(L^{p}(w)\to L^{q}(w))$ with $1\leq p<q\leq \infty.$ Then, there are $\alpha, \beta>0$ such that for each ball $B$ with radius $r$ and any $t>0,$ \begin{equation}\label{eq: cd104} \l(\fint_{B}|T_{t}(f1_{B})|^{q}dw\r)^{1/q}\lesssim \max\l\{\l(\frac{r}{t^{1/2m}}\r)^{\alpha}, \l(\frac{r}{t^{1/2m}}\r)^{\beta}\r\} \l(\fint_{B}|f|^{p}dw\r)^{1/p}.
\end{equation} \end{lemma}
{\it Proof.}\quad Following \cite[Proposition 2.4]{AM} mutatis mutandis, we note that in Definition \ref{definition: cd24}, the estimates \eqref{eq: cd25}-\eqref{eq: cd27} (for any $t>0$) are equivalent to that for $r\approx t^{1/2m}.$ Furthermore, if these estimates hold for $r\approx t^{1/2m},$ then \eqref{eq: cd25} holds generally with constant $\max\{\l(\frac{r}{t^{1/2m}}\r)^{\alpha}, 1\}$ (some $\alpha>0$) where $1$ applies when $r\leq t^{1/2m}.$ To obtain \eqref{eq: cd104}, it thus suffices to refine this constant: substituting $1$ with $\l(\frac{r}{t^{1/2m}}\r)^{\beta}$ when $r\leq t^{1/2m}.$ For a parallel line of reasoning, consult the argument at the beginning of \cite[Lemma 7.5]{DMR}.


Let $B:=B(x, r)$ with $r\leq t^{1/2m},$ then $B\subset B_{t}:=B(x, t^{1/2m}).$ Since $w\in A_{\infty},$ there exists a $\eta>0$ such that $$\frac{w(B)}{w(B_{t})}\lesssim \l(\frac{|B|}{|B_{t}|}\r)^{\eta}\lesssim \l(\frac{r}{t^{1/2m}}\r)^{\eta n}.$$ From this, together with \eqref{eq: cd25} for $T_{t},$ it follows that \begin{equation*}\begin{aligned} \l(\fint_{B}|T_{t}(f1_{B})|^{q}dw\r)^{1/q}&\lesssim \l(\frac{w(B)}{w(B_{t})}\r)^{1/q}\l(\fint_{B_{t}}|T_{t}(f1_{B})|^{q}dw\r)^{1/q} \\[4pt] 
&\lesssim \l(\frac{w(B)}{w(B_{t})}\r)^{\frac{1}{p}-\frac{1}{q}}\l(\fint_{B}|f1_{B}|^{p}dw\r)^{1/p}\lesssim \l(\frac{r}{t^{1/2m}}\r)^{\beta}\l(\fint_{B}|f1_{B}|^{p}dw\r)^{1/p},\end{aligned}\end{equation*} 
where $\beta:=(\frac{1}{p}-\frac{1}{q})\eta n.$ This yields \eqref{eq: cd104}.

\hfill$\Box$

\subsection{Theorems on weighted boundedness of sublinear operators.} As our proof strategy is consistent with that in \cite{DMR}, the first two theorems below will play a central role.

\begin{theorem}\label{theorem: cd18}\;(\cite[Theorem 2.2]{AM1})\;  Given $w\in A_{2}$ and $1\leq p_{0}<q_{0}\leq \infty,$ let $\T$ be a sublinear operator acting on $L^{p_{0}}(w),$ $\{\A_{r}\}_{r>0}$ a family of operators acting from a subspace $\D$ of $L^{p_{0}}(w)$ into $L^{p_{0}}(w),$ and $S$ an operator from $\D$ into the space of measurable functions on $\rz.$ Suppose that every $f\in \D$ and ball $B$ with radius $r,$ \begin{equation}\label{eq: cd19}\l(\fint_{B}|\T(I-\A_{r})f|^{p_{0}}dw\r)^{1/p_{0}}\leq \sum_{j\geq 1}g(j)\l(\fint_{2^{j+1}B}|Sf|^{p_{0}}dw\r)^{1/p_{0}},\end{equation} 
\begin{equation}\label{eq: cd20}\l(\fint_{B}|\T \A_{r}f|^{q_{0}}dw\r)^{1/q_{0}}\leq \sum_{j\geq 1}g(j)\l(\fint_{2^{j+1}B}|\T f|^{p_{0}}dw\r)^{1/p_{0}},\end{equation} where $\sum_{j\geq 1}g(j)<\infty.$
Then for every $p,$ $p_{0}<p<q_{0},$ and weights $v\in A_{\frac{p}{p_{0}}}(w)\cap \mbox{RH}_{(\frac{q_{0}}{p})'}(w),$ there is a constant $C$ such that for all $f\in \D,$ $$\|\T f\|_{_{L^{p}(vdw)}}\leq C \|\S f\|_{_{L^{p}(vdw)}}.$$

\end{theorem} 

\begin{theorem}\label{theorem: cd21}\;(\cite[Theorem 2.4]{AM1})\;  Given $w\in A_{2}$ with doubling order $D$ and $1\leq p_{0}<q_{0}\leq \infty,$ let $\T: L^{q_{0}}(w)\to L^{q_{0}}(w)$ be a sublinear operator, $\{\A_{r}\}_{r>0}$ a family of operators acting from $L_{c}^{\infty}$ into $L^{q_{0}}(w).$ Suppose that for every ball $B$ with radius $r,$ $f\in L_{c}^{\infty}$ with $\supp f\subset B$ and $j\geq 2,$
\begin{equation}\label{eq: cd22}\l(\fint_{C_{j}(B)}|\T(I-\A_{r})f|^{p_{0}}dw\r)^{1/p_{0}}\leq g(j)\l(\fint_{2^{j+1}B}|f|^{p_{0}}dw\r)^{1/p_{0}}.\end{equation} Suppose further that for every $j\geq 1,$

\begin{equation}\label{eq: cd23}\l(\fint_{C_{j}(B)}|\A_{r}f|^{q_{0}}dw\r)^{1/q_{0}}\leq g(j)\l(\fint_{B}|f|^{p_{0}}dw\r)^{1/p_{0}},\end{equation} where $\sum_{j\geq 1}g(j)2^{jD}<\infty.$
Then for every $p,$ $p_{0}<p<q_{0},$ there is a constant $C$ such that for all $f\in L_{c}^{\infty},$ $$\|\T f\|_{_{L^{p}(w)}}\leq C \| f\|_{_{L^{p}(w)}}.$$
\end{theorem} 

\begin{remark}\label{remark: cd1084}\; In Definition \ref{definition: cd24}-\ref{definition: cd1083} and Theorem \ref{theorem: cd18}-\ref{theorem: cd21}, the case $q=q_{0}=\infty$ is understood as follows: the $L^{q}(w)$ (resp., $L^{q_{0}}(w)$)-average is replaced by the essential supremum. Moreover, if $q_{0}=\infty$ in Theorem \ref{theorem: cd18}, the condition on $v$ becomes $v\in A_{\frac{p}{p_{0}}}(w).$\end{remark} 

We also need the theorem below, a special case of \cite[Theorem 3.1]{AM2} (formulated for spaces of homogeneous type in \cite[Section 5]{AM2}). 

\begin{theorem}\label{theorem: cd1035}\;(\cite[Theorem 9.10]{DMR})\;  Given $1<q<\infty,$ $a\geq 1$ and $u\in \mbox{RH}_{s'}(w),$ $1<s<\infty.$ There exists a $C>1$ with the following property: suppose $F\in L^{1}(w)$ and $G$ are nonnegative measurable functions such that for any ball $B$ there are nonnegative functions $G_{B}$ and $H_{B}$ with \begin{equation}\label{eq: cd5000}F(x)\leq G_{B}(x)+H_{B}(x)\quad \mbox{for a.e. }\; x\in B, \end{equation} and 
 \begin{equation}\label{eq: cd1036} \l(\fint_{B}|H_{B}|^{q}dw\r)^{1/q}\leq a M_{w}(F)(x), \quad \fint_{B}G_{B}dw\leq G(x), \;\mbox{for all }\; x\in B,
\end{equation} where $M_{w}$ is the Hardy-Littlewood function with respect to $dw.$ Then for $1<t<q/s,$
 \begin{equation}\label{eq: cd1037} \|M_{w}(F)\|_{L^{t}(udw)}\leq C\|G\|_{L^{t}(udw)}.
\end{equation}

\end{theorem}

\section{Off-diagonal estimates for the semigroup of $e^{-t\L_{w}}$} 

For $\{a_{\alpha, \beta}(x)\}_{|\alpha|=|\beta|=m}\in \Ez(w, c_{1}, c_{2})$ with $w\in A_{2},$ define $\B(u, v)$ to be the sesquilinear form 
\begin{equation}\label{eq: cd3001}\B(u, v):=\sum_{|\alpha|=|\beta|=m}\int_{\rz}a_{\alpha, \beta}(x)\partial^{\alpha}u(x) \cdot \overline{\partial^{\beta}v(x)}dx.\end{equation} Clearly, $\B$ is a closed, maximally accretive, and continuous sesquilinear form, and there exists an operator $\L_{w}$ (denoted by \eqref{eq: z00}) with domain $\D(\L_{w}):=\{u\in H^{m}(w): \L_{w} u\in L^{2}(w)\}$ such that for all $u\in \D(\L_{w})$ and  $v\in H^{m}(w),$ 
\begin{equation}\label{eq: cd3002}<\L_{w} u, v>=\int_{\rz}\L_{w} u\overline{v}dw=\B(u, v).\end{equation} In particular, $\D(\L_{w})\subset H^{m}(w)$ is dense in $L^{2}(w).$ Similarly, we can define \begin{equation}\label{eq: cd3003}\L_{w}^{*}:=\sum_{|\alpha|=|\beta|=m}(-1)^{|\beta|}w^{-1}(\partial^{\beta}\overline{a_{\alpha, \beta}}\partial^{\alpha})\end{equation} which is the adjoint of $\L_{w}$ with respect to $L^{2}(w)$ via the sesquilinear form $\B^{*}(u, v):=\overline{\B(v, u)}.$ For details on these properties, one may refer to \cite{GZ}. 

Define \begin{equation}\label{eq: cd9001}\V:=\sup \{ |\mbox{arg}\;<\L_{w} f, f>|: f\in \D(\L_{w})\}.\end{equation} From \eqref{eq: z1}-\eqref{eq: z2}, it follows that $0<\V<\frac{\pi}{2}$ and $\L_{w}$ is an operator of type $\V.$ That is, $\L_{w}$ is closed and densely defined, with its spectrum contained in $\Sigma_{\V},$ and its resolvent satisfies $$\|(\xi-\L_{w})^{-1}f\|_{L^{2}(w)}\leq \frac{C_{\mu, \V}}{|\xi|}\|f\|_{L^{2}(w)}\;\quad \mbox{for any }\;\xi\in \cc\;\mbox{with}\; |\mbox{arg}\; \xi|\geq \mu>\V.$$ Then there exists a complex semigroup $e^{-z\L_{w}}$ on $\Sigma_{\frac{\pi}{2}-\V}$ of bounded operators on $L^{2}(w),$ along with an $L^{2}(w)-$functional calculus as in \cite{H, TK, M}.

\subsection{$H^{\infty}$ Functional calculi in $L^{2}(w)$.}  Let $\mu\in (\V, \pi)$ and $\H^{\infty}(\Sigma_{\mu})$ be the collection of bounded holomorphic functions on $\Sigma_{\mu}.$ If $\phi\in \H^{\infty}(\Sigma_{\mu})$ satisfies, for some $s>0,$ $$|\phi(z)|\lesssim \frac{|z|^{s}}{(1+|z|)^{2s}}\quad z\in \Sigma_{\mu},$$ we say that $\phi\in \H_{0}^{\infty}(\Sigma_{\mu}).$ We are able to define $\phi(\L_{w})$ for any $\phi\in \H_{0}^{\infty}(\Sigma_{\mu})$ thanks to the $L^{2}(w)-$functional calculus of $\L_{w}.$ Indeed, $\phi(\L_{w})$ has an intergral representation. Let $\V<\theta<\nu<\min\{\mu, \frac{\pi}{2}\},$ and let $\Gamma_{\pm}, \gamma_{\pm}$ be the half-rays $\rr^{+}e^{\pm i(\frac{\pi}{2}-\theta)}$ and $\rr^{+}e^{\pm i\nu},$ respectively. Then 
\begin{equation}\label{eq: cd3004}\phi(\L_{w}):=\int_{\Gamma_{+}}e^{-z\L_{w}}\eta_{+}(z)dz+\int_{\Gamma_{-}}e^{-z\L_{w}}\eta_{-}(z)dz,\end{equation} where \begin{equation}\label{eq: cd3005}|\eta_{\pm}(z)|=\frac{1}{2\pi i}\int_{\gamma_{\pm}}e^{\xi z}\phi(\xi)d\xi, \quad z\in \Gamma_{\pm}.\end{equation} It is easy to see that the integrals in \eqref{eq: cd3004} converge in $L^{2}(w).$ According to \cite{H, TK, M}, any operator $\L_{w}$ as above admits a bounded holomorphic functional calculus. That is, given $\mu\in (\V, \pi)$: \bigskip

\noindent{\textbf{(a)}\;  for any $\phi\in \H^{\infty}(\Sigma_{\mu}),$  the operator $\phi(\L_{w})$ can be defined and is boounded on $L^{2}(w)$ with \begin{equation}\label{eq: cd3006} \|\phi(\L_{w}) f\|_{L^{2}(w)}\leq C \|\phi\|_{\infty} \|f\|_{L^{2}(w)},\end{equation} where $C$ is independent of $\V$ and $\mu.$   } 

\noindent{\textbf{(b)}\; the product rule $\phi(\L_{w})\psi(\L_{w})=(\phi\psi)(\L_{w})$ holds for any $\phi, \psi\in \H^{\infty}(\Sigma_{\mu}).$  } 

\noindent{\textbf{(c)}\; for any sequence $\{\phi_{k}\}\subset \H^{\infty}(\Sigma_{\mu})$ converging uniformly on compact subsets of $\Sigma_{\mu}$ to $\phi,$ we have that $\phi_{k}(\L_{w})$ converges to $\phi(\L_{w})$ strongly in $L^{2}(w).$} 

\noindent{\textbf{(d)}\; for any operator $\L_{w}$ as above and for any $f\in \H_{0}^{\infty}(\Sigma_{\mu}),$ the following square function estimate holds:\begin{equation}\label{eq: cd3007}\l(\int_{0}^{\infty}\|\phi(t\L_{w})f\|_{L^{2}(w)}^{2}\frac{d t}{t}\r)^{1/2}\leq C \|\phi\|_{\infty} \|f\|_{L^{2}(w)},\end{equation} the same is true for $\L_{w}^{*}.$} 

One can extend the $H^{\infty}$ functional calculus to more general holomorphic functions (such as powers), with $\phi(\L_{w})$ defined as unbounded operators.

\subsection{Off-diagonal estimates in $L^{2}(w)$} Armed with the $L^{2}(w)-$functional calculus for $\L_{w},$ the (full) $L^{2}(w)-$off-diagonal estimates for the complex semigroup $e^{-z\L_{w}}$ and its gradients can be proven. Preceding the proof, the following lemma for the resolvent operators are required.

\begin{lemma}\label{lemma: cd1074} Given $w\in A_{2}$ and $\{a_{\alpha, \beta}(x)\}_{|\alpha|=|\beta|=m}\in \Ez(w, c_{1}, c_{2}).$  Let $E$ and $F$ be two closed sets. Fix $\nu$ such that $0<\nu<\pi-\V$ and $z\in \Sigma_{\nu}.$ Then there exist constants $C$ and $c$ depending only on $n, m, c_{1}, c_{2}, \nu$ such that for all $f\in L^{2}(w)$ and $\overrightarrow{f}=(f_{\beta})_{|\beta|=m}$ with $f_{\beta}\in L^{2}(w),$
\begin{equation*}\begin{aligned} &(i)\quad\quad  \|(1+z^{2m}\L_{w})^{-1}(f1_{E})1_{F}\|_{L^{2}(w)}\leq C e^{-c\frac{d(E, F)}{|z|}}\|f1_{E}\|_{L^{2}(w)},\\[4pt] 
&(ii)\quad\quad\|z^{m}\nabla^{m}(1+z^{2m}\L_{w})^{-1}(f1_{E})1_{F}\|_{L^{2}(w)}\leq C e^{-c\frac{d(E, F)}{|z|}}\|f1_{E}\|_{L^{2}(w)},\\[4pt] 
&(iii)\quad\quad\|z^{m}(1+z^{2m}\L_{w})^{-1}\frac{1}{w}\mbox{div}_{m}(w\overrightarrow{f}1_{E})1_{F}\|_{L^{2}(w)}\leq C e^{-c\frac{d(E, F)}{|z|}}\|f1_{E}\|_{L^{2}(w)},\end{aligned} 
 \end{equation*} where $\mbox{div}_{m}\overrightarrow{f}:=\sum_{|\beta|=m}\partial^{\beta}f_{\beta}.$\end{lemma} 
{\it Proof.}\quad This proof is a variant of the arguments presented in \cite[Lemma 4.1]{CR2} and \cite[Lemma 2.10]{CR1}; for additional reference, see also the proof of \cite[Lemma 4.2]{GZ}.

We first prove $(i)$ and $(ii)$. Assume $0<\nu<\frac{\pi}{2};$ without loss of generality, take $\nu>\frac{\pi}{4}.$ We also assume $\Delta:=\frac{\kappa d(E, F)}{|z|^{1/2m}}\geq 1,$ with $\kappa$ a sufficiently small constant to be determined subsequently. Through the change of variables $z\to z^{2m},$ it suffices to build the following two inequalities: \begin{equation}\label{eq: cd1075} \int_{F}|(1+z\L_{w})^{-1}f|^{2}dw\leq C e^{-c\frac{d(E, F)}{|z|^{1/2m}}}\int_{E}|f|^{2}dw\end{equation} and \begin{equation}\label{eq: cd1079} \int_{F}|z\nabla^{m}(1+z\L_{w})^{-1}f|^{2}dw\leq C e^{-c\frac{d(E, F)}{|z|^{1/2m}}}\int_{E}|f|^{2}dw,\end{equation}
where $f\in L^{2}(w)$ is arbitrary and supported in $E.$

For simplicity, set $u^{z}=(1+z\L_{w})^{-1}f,$ so that $f=u^{z}+z\L_{w} u^{z}.$ By \eqref{eq: cd3002}, we have for all $v\in H^{m}(w)$ that 
$$\int_{\rz}u^{z}(x)\overline{v(x)}wdx+z\sum_{|\alpha|=|\beta|=m}\int_{\rz}a_{\alpha, \beta}(x)\partial^{\alpha}u^{z}(x) \cdot \overline{\partial^{\beta}v(x)}dx=\int_{\rz}f(x)\overline{v(x)}wdx.$$ In the latter equality, we take $v=u^{z}\eta^{2}$ with $\eta=e^{\Delta \tilde{\eta}}-1.$ (Here $\tilde{\eta}\in C_{0}^{\infty}(\rz\setminus E)$ is a non-negative function, satisfying $0\leq \tilde{\eta}\leq 1,$ $\tilde{\eta}\equiv 1$ on $F$ and $|\partial^{\gamma}\tilde{\eta}|\lesssim d(E, F)^{-|\gamma|}$ for any $|\gamma|\leq m.$) Consequently, 
it holds that 
\begin{equation}\label{eq: cd1076}\begin{aligned} \int_{\rz} |u^{z}(\eta+1)|^{2}dw&+z\int_{\rz}a_{\alpha, \beta}(x)\partial^{\beta}(u^{z}(\eta+1))\overline{\partial^{\alpha}(u^{z}(\eta+1))}dx\\[4pt] 
&=z\int_{\rz}a_{\alpha, \beta}(x)\l[\partial^{\beta}(u^{z}(\eta+1))\overline{\partial^{\alpha}(u^{z}(\eta+1))}-\partial^{\beta}u^{z}\overline{\partial^{\alpha}(u^{z}\eta^{2})}\r]dx\\[4pt] 
&\quad \quad+\int_{\rz} |u^{z}|^{2}(2\eta+1)wdx+f(x)\overline{ u^{z}\eta^{2} }dw:=G_{1}+G_{2}+G_{3}.\end{aligned}
 \end{equation} 

To proceed, we split $G_{1}$ into \begin{equation*}\begin{aligned} G_{1}&=z\int_{\rz}a_{\alpha, \beta}(x)\l[\partial^{\beta}(u^{z}(\eta+1))\overline{\partial^{\alpha}(u^{z}(\eta+1))}-\partial^{\beta}u^{z}\overline{\partial^{\alpha}(u^{z}(\eta+1)^{2})}\r]dx\\[4pt] 
&\quad\quad-z\int_{\rz}a_{\alpha, \beta}(x)\partial^{\beta}u^{z}\overline{\partial^{\alpha}(u^{z}(2\eta+1))}:=G_{11}+G_{12},\end{aligned} 
\end{equation*} furthermore, by Leibniz's rule, \begin{equation*}\begin{aligned} G_{11}&=z\sum_{|\tau|+|\gamma|<2m}C_{\alpha}^{\tau}C_{\beta}^{\gamma}\int_{\rz}a_{\alpha, \beta}(x)\partial^{\gamma}u^{z} \overline{\partial^{\tau}u^{z}}\partial^{\beta-\gamma}(\eta+1)\partial^{\alpha-\tau}(\eta+1)dx\\[4pt] 
&-z\sum_{\tau<\alpha}C_{\alpha}^{\tau} \int_{\rz}a_{\alpha, \beta}(x)\partial^{\beta}u^{z}\overline{\partial^{\tau}u^{z}}\partial^{\alpha-\tau}(\eta+1)^{2}dx:=G_{111}+G_{112}.\end{aligned} 
\end{equation*} 

From the definition of $\eta,$ a computation leads to that, for any $|\xi|\leq m,$ \begin{equation}\label{eq: ab2.36}\partial^{\xi}(\eta+1)=(\eta+1)P^{\Delta}_{\xi}(\partial_{1},...,\partial_{n})\tilde{\eta}, \end{equation}where $P^{\Delta}_{\xi}$ denotes a homogeneous polynomial of degree $|\xi|$ ($P^{\Delta}_{0}:=1$) satisfying \begin{equation}\label{eq: ab2.34}|P^{\Delta}_{\xi}(\partial_{1},...,\partial_{n})\tilde{\eta}|\lesssim \l(\frac{\Delta}{d}\r)^{|\xi|},\quad (\Delta\geq 1, d:=d(E, F)),\end{equation} and \begin{equation}\label{eq: ab2.35}\partial^{\xi}u^{z}(\eta+1)=\sum_{\tau\leq \xi}P_{\xi-\tau}^{\Delta}(\partial_{1},...,\partial_{n})\tilde{\eta} \partial^{\tau}(u^{z}(\eta+1)).\end{equation} 
Using \eqref{eq: ab2.36}-\eqref{eq: ab2.35}, an estimate for $G_{111}$ can be obtained by disregarding summations and constant factors, as shown below:
\begin{equation}\label{eq: ab2.38}\begin{aligned}
G_{111}&=z\int_{\rz}w^{-1}a_{\alpha, \beta}\overline{\partial^{\tau}u^{z}}\partial^{\gamma}u^{z}\partial^{\alpha-\tau}(\eta+1)\partial^{\beta-\gamma}(\eta+1)dw\\
&=z\int_{\rz}w^{-1}a_{\alpha, \beta}\overline{\partial^{\tau}u^{z}}\partial^{\gamma}u^{z}(\eta+1)^{2}P^{\Delta}_{\alpha-\tau}(\partial_{1},...,\partial_{n})\tilde{\eta}P^{\Delta}_{\beta-\gamma}(\partial_{1},...,\partial_{n})\tilde{\eta}dw\\
&=z\sum_{S\leq \tau}\sum_{\xi\leq \gamma}\int_{\rz}w^{-1}a_{\alpha, \beta}P_{\gamma-\xi}^{\Delta}(\partial_{1},...,\partial_{n})\tilde{\eta} \partial^{\xi}(u^{z}(\eta+1))\\
&\quad\quad\times P_{\tau-S}^{\Delta}(\partial_{1},...,\partial_{n})\tilde{\eta} \overline{\partial^{S}(u^{z}(\eta+1))}P^{\Delta}_{\beta-\gamma}(\partial_{1},...,\partial_{n})\tilde{\eta}P^{\Delta}_{\alpha-\tau}(\partial_{1},...,\partial_{n})\tilde{\eta}dw\\
&\lesssim\lambda^{2m}\sum_{S\leq \tau}\sum_{\xi\leq \gamma}\l(\frac{\Delta}{d}\r)^{|\alpha-S|+|\beta-\xi|}\|\partial^{\xi}(u^{z}(\eta+1))\|_{L^{2}(w)}\|\partial^{S}(u^{z}(\eta+1))\|_{L^{2}(w)}\quad (\lambda:=|z|^{\frac{1}{2m}})\\
&\lesssim  \kappa \sum_{S\leq \tau}\sum_{\xi\leq \gamma}\l(\lambda^{|S|}\|\partial^{S}(u^{z}(\eta+1))\|_{L^{2}(w)}\r)\l(\lambda^{|\xi|}\|\partial^{\xi}(u^{z}(\eta+1))\|_{L^{2}(w)}\r)\quad (|\xi|+|S|\leq 2m-1)\\\end{aligned}\end{equation}
\begin{equation*}\begin{aligned}
&\lesssim  \kappa\sum_{S\leq \tau}\sum_{\xi\leq \gamma}C(\xi, S, m) \lambda^{|S|}\l(\int_{\rz}|u^{z}(\eta+1)|^{2}w\r)^{\frac{ (1-\frac{|S|}{m}) } {2}}\l(\int_{\rz}|\nabla^{m}(u^{z}(\eta+1)|^{2})w\r)^{\frac{|S|}{2m}}\\
&\quad\quad\times \lambda^{|\xi|}\l(\int_{\rz}|u^{z}(\eta+1)|^{2}w\r)^{\frac{ (1-\frac{|\xi|}{m}) } {2}}\l(\int_{\rz}|\nabla^{m}(u^{z}(\eta+1)|^{2})w\r)^{\frac{|\xi|}{2m}}\\
&\lesssim \kappa(\|u^{z}(\eta+1)\|_{L^{2}(w)}^{2}+|z|\|\nabla^{m}(u^{z}(\eta+1))\|_{L^{2}(w)}^{2}).
\end{aligned}\end{equation*} Similarly, $$G_{112}\lesssim \kappa(\|u^{z}(\eta+1)\|_{L^{2}(w)}^{2}+|z|\|\nabla^{m}(u^{z}(\eta+1))\|_{L^{2}(w)}^{2}).$$ For $G_{12},$ we can apply Young's inequality to derive $$G_{12}\lesssim  |z|\|\nabla^{m}u^{z}\|_{L^{2}(w)}^{2}+\ez |z|\int_{\rz}|\nabla^{m}(u^{z}(\eta+1))|^{2}dw. $$ By the same token, $$G_{3}\lesssim \|f\|_{L^{2}(w)}^{2}+\ez\|u^{z}(\eta+1)\|_{L^{2}(w)}^{2}+\|u^{z}\|_{L^{2}(w)}^{2}.$$ Observing that $\|\eta\|_{\infty}\lesssim e^{\Delta},$ we can bound $G_{2}$ by $$G_{2}\lesssim e^{\Delta}\|u^{z}\|_{L^{2}(w)}^{2}.$$

We now turn to estimating $G_{4},$ written as $G_{4}=z\cdot G_{5},$ where $G_{5}$ is given by: $$G_{5}:=z\sum_{|\alpha|=|\beta|=m}\int_{\rz}a_{\alpha, \beta}(x)\partial^{\beta}(u^{z}(\eta+1))\overline{\partial^{\alpha}(u^{z}(\eta+1))}dx.$$ To the end, we introduce $$\R:=\int_{\rz}|u^{z}(\eta+1)|^{2}dw,\; \S:=\mbox{Re}\; G_{5}, \;\T:=\mbox{Im}\; G_{5} \;\;\mbox{and}\;\; z=s+it.$$ Apparently, by \eqref{eq: z1}-\eqref{eq: z2}, $$\S\geq c_{1} \|\nabla^{m}(u^{z}(\eta+1))\|_{L^{2}(w)}^{2}\quad\mbox{and}\quad |\T|\leq \frac{c_{2}}{c_{1}}\S.$$  Set $\rho:=\frac{c_{1}}{c_{2}\tan \nu};$ then $\rho<1.$ Note also that $|t|\leq s\tan \nu.$ a standard argument yields $$|\R+G_{4}|=|\R+(s+it)(\S+i\T)|\geq \frac{\rho^{1/2} \R}{2}+\frac{|z| \S}{2}.$$ Thus, recalling \eqref{eq: cd1076} and summarizing all estimates we get
 \begin{equation*}\begin{aligned} \int_{\rz}|u^{z}(\eta+1)|^{2}dw &+|z|\int_{\rz}|\nabla^{m}(u^{z}(\eta+1))|^{2}dw\\[4pt] 
&\lesssim \kappa \int_{\rz}|u^{z}(\eta+1)|^{2}dw+(\kappa +\ez)|z|\int_{\rz}|\nabla^{m}(u^{z}(\eta+1))|^{2}dw\\[4pt] 
&\quad\quad +|z|\int_{\rz}|\nabla^{m} u^{z}|^{2}dw+e^{\Delta}\int_{\rz}|u^{z}|^{2}dw,\end{aligned} 
\end{equation*} from which, by letting $\kappa$ and $\ez$ small and also using the property of $\eta, \tilde{\eta},$ it follows that \begin{equation}\label{eq: cd1077}e^{2\Delta}(\|u^{z}\|_{L^{2}(F, w)}^{2}+|z|\|\nabla^{m}u^{z}\|_{L^{2}(w)}^{2})\lesssim |z|\|\nabla^{m} u^{z}\|_{L^{2}(F, w)}^{2}+e^{\Delta}\|u^{z}\|_{L^{2}(w)}^{2}.\end{equation}
Adapting the proof technique from \cite[Lemma 2.8]{CR1} (or \cite[Lemma 4.1]{GZ}), we can prove the uniform bound: \begin{equation}\label{eq: cd1078}\sup_{z\in \Sigma_{\tau}}\l(\|(1+z\L_{w})^{-1}f\|_{L^{2}(w)\to L^{2}(w)}+\|z^{1/2}\nabla^{m}(1+z\L_{w})^{-1}f\|_{L^{2}(w)\to L^{2}(w)}\r)\leq C,\end{equation} where $C$ depends only on $n ,m, c_{1}, c_{2}, \tau.$ Deatils are left to the reader. Inserting \eqref{eq: cd1078} into \eqref{eq: cd1077} we then arrive at \eqref{eq: cd1075} and \eqref{eq: cd1079}. 

It remains to consider the case $\nu\in (\frac{\pi}{2}, \pi-\V).$ Note that there always exist $\nu_{1}<\frac{\pi}{2}$ and $\tau<\frac{\pi}{2}-\V$ such that every $z\in \Sigma_{\nu}$ admits a decomposition $z=z_{1}\xi,$ where $\xi$ is fixed with $|\xi|=1$ and $\mbox{arg} (\xi)\leq \tau,$ and $z_{1}\in \Sigma_{\nu_{1}}.$ Then, substituting $z=z_{1}\xi$ into the left side of \eqref{eq: cd1075} and introducing $\L_{w}^{1}:=\xi \L_{w}$ we have $$\int_{F}|(1+z\L_{w})^{-1}f|^{2}wdx=\int_{F}|(1+z_{1} \L_{w}^{1})^{-1}f.$$ Invoking Lemma \ref{lemma: h2}, we see $\L_{w}^{1} \in \Ez(w, \lambda_{\xi}, \Lambda_{\xi}).$ Repeating the above procedure for $\L_{w}^{1}$ yields the desired estimates.

Eventually, \eqref{eq: cd3003} implies that the estimate $(ii)$ remains valid if $\L_{w}$ is replaced by its adjoint $\L_{w}^{*}.$ From this, a duality argument (as in \cite[Lemma 2.10]{CR1}) leads to conclusion $(iii)$. 

\hfill$\Box$
  
We now elaborate on the proof of the (full) off-diagonal estimates in $L^{2}(w)$ for the complex semigroup $e^{-z\L_{w}}$ and its gradients.
\begin{theorem}\label{theorem: cd1073} Given $w\in A_{2}$ and $\{a_{\alpha, \beta}(x)\}_{|\alpha|=|\beta|=m}\in \Ez(w, c_{1}, c_{2}).$  For all closed sets $E$ and $F,$ $f\in L^{2}(w),$ $0\leq k\leq m$ and $z\in \Sigma_{\nu}$ with $0<\nu<\frac{\pi}{2}-\V,$ we have: 
\begin{equation*}\begin{aligned} &(i)\quad\quad  \|z^{\frac{k}{2m}}\nabla^{k}e^{-z\L_{w}}(f1_{E})1_{F}\|_{L^{2}(w)}\lesssim e^{-c\l(\frac{d(E, F)}{|z|^{\frac{1}{2m}}}\r)^{\frac{2m}{2m-1}}}\|f1_{E}\|_{L^{2}(w)},\\[4pt] 
\end{aligned}
 \end{equation*} 
\begin{equation*}\begin{aligned}  
 &(ii)\quad\quad\|z\L_{w}e^{-z\L_{w}}(f1_{E})1_{F}\|_{L^{2}(w)}\lesssim e^{-c\l(\frac{d(E, F)}{|z|^{\frac{1}{2m}}}\r)^{\frac{2m}{2m-1}}}\|f1_{E}\|_{L^{2}(w)}. \\[4pt] \end{aligned}
 \end{equation*} 
\end{theorem} 
{\it Proof.}\quad Set $d:=d(E, F).$ Proving $(i)$ and $(ii)$ for $d^{2m}\geq |z|$ is sufficient. Fix $\theta$ with $\frac{\pi}{2}+|\mbox{arg} z|< \theta<\pi-\V$ and a parameter $\rho>0$ (to be determined later), and define $$\Gamma_{\theta}^{\pm}:=\{re^{\pm i  \theta}: r\geq \rho\}\quad \mbox{and}\quad \Gamma_{\theta}:= \{re^{i\phi}: |\phi|\leq \theta\}.$$ Using the $L^{2}(w)-$functional calculus of $\L_{w}$ again, we may express $e^{-z\L_{w}}$ through the intergral $$e^{-z\L_{w}}f=\frac{1}{2\pi}\int_{\Gamma_{\theta}^{\pm}\cup\Gamma_{\theta}}e^{z\xi}(\xi+\L_{w})^{-1}fd\xi.$$ From this formular, in conjunction with Lemma \ref{lemma: cd1074}, it follows that 
\begin{equation*}\begin{aligned} \l(\int_{F} \bigg| \int_{\Gamma_{\theta}^{\pm}} e^{z\xi} (\xi+\L_{w})^{-1}(f1_{E})d\xi\bigg|^{2}w(x)dx\r)^{\frac{1}{2}}&\lesssim\int_{\Gamma_{\theta}^{\pm}} |e^{z\xi}|\l(\int_{F} |(\xi+\L_{w})^{-1}(f1_{E}) |^{2}dw \r)^{\frac{1}{2}}|d\xi| \\[4pt] 
&\lesssim \int_{\Gamma_{\theta}^{\pm}} |e^{z\xi}| |\xi|^{-1} e^{-cd |\xi|^{\frac{1}{2m}}}\|f1_{E}\|_{L^{2}(w)}|d\xi|\\[4pt] 
&\lesssim e^{-cd \rho^{\frac{1}{2m}}} (|z|\rho)^{-1}e^{-c'\rho |z|},\end{aligned}
 \end{equation*} moreover, 
\begin{equation*}\begin{aligned} \l(\int_{F} \bigg| \int_{\Gamma_{\theta}} e^{z\xi} (\xi+\L_{w})^{-1}(f1_{E})d\xi\bigg|^{2}w(x)dx\r)^{\frac{1}{2}}&\lesssim \int_{-\theta}^{\theta} \rho^{-1}e^{|z|\rho}e^{-cd \rho^{\frac{1}{2m}}} \rho d\phi\|f1_{E}\|_{L^{2}(w)} \\[4pt] 
&\lesssim e^{-c''d \rho^{\frac{1}{2m}}} e^{ |z| \rho}.\end{aligned}\end{equation*} Collecting the above two estimates we get 
\begin{equation}\label{eq: cd1080}\l(\int_{F} |e^{-z\L_{w}}(f1_{E})|^{2}w(x)dx\r)^{\frac{1}{2}} \lesssim e^{-cd \rho^{\frac{1}{2m}}} (|z|\rho)^{-1}e^{-c'\rho |z|}+e^{-c''d \rho^{\frac{1}{2m}}} e^{ |z| \rho}. \end{equation} By \eqref{eq: cd1080}, if we let $\rho=\ez \frac{d^{\frac{2m}{2m-1} }}{|z|^{\frac{2m}{2m-1}}}$ with $\ez$ small enough, then $$e^{-cd \rho^{\frac{1}{2m}}} (|z|\rho)^{-1}e^{-c'\rho |z|}+e^{-c''d \rho^{\frac{1}{2m}}} e^{ |z| \rho}\lesssim e^{-c\l(\frac{d(E, F)}{|z|^{\frac{1}{2m}}}\r)^{\frac{2m}{2m-1}}}.$$ Hence, we conclude with \begin{equation}\label{eq: cd3008}\|e^{-z\L_{w}}(f1_{E})1_{F}\|_{L^{2}(w)}\lesssim e^{-c\l(\frac{d(E, F)}{|z|^{\frac{1}{2m}}}\r)^{\frac{2m}{2m-1}}}\|f1_{E}\|_{L^{2}(w)}.\end{equation} 

A similar argument leads to \begin{equation}\label{eq: cd3009}\|z^{\frac{1}{2}}\nabla^{m}e^{-z\L_{w}}(f1_{E})1_{F}\|_{L^{2}(w)}\lesssim e^{-c\l(\frac{d(E, F)}{|z|^{\frac{1}{2m}}}\r)^{\frac{2m}{2m-1}}}.\end{equation} Conclusion $(i)$ in Theorem \ref{theorem: cd1073} therefore follows by \eqref{eq: a2.32} and \eqref{eq: cd3008}-\eqref{eq: cd3009}. 

Observe that $$(z^{2m}\L_{w}(1+z^{2m}\L_{w})^{-1}(f1_{E}))1_{F}=-(1+z^{2m}\L_{w})^{-1}(f1_{E})1_{F}$$ since the two sets $E$ and $F$ are disjoint. Then, by Lemma \ref{lemma: cd1074}, $$\|z^{2m}\L_{w}(1+z^{2m}\L_{w})^{-1}(f1_{E})1_{F}\|_{L^{2}(w)}\lesssim e^{-c\frac{d(E, F)}{|z|}}\|f1_{E}\|_{L^{2}(w)}.$$
The above argument, applied similarly, yields conclusion $(ii)$ in Theorem \ref{theorem: cd1073}.
\hfill$\Box$

\subsection{Off-diagonal estimates in $L^{p}(w)$} Owing to Theorem \ref{theorem: cd1073}, Definition \ref{definition: cd24}-\ref{definition: cd1083} and Lemma \ref{lemma: h7}, we see that $2\in \widetilde{\J}(\L_{w})$ ($\widetilde{\J}(\L_{w}):=\{p\in [1, \infty]: \sup_{t>0}\|e^{-t\L_{w}}\|_{L^{p}(w)\to L^{p}(w)}\lesssim 1\}$) and $e^{-t\L_{w}}\in \mho(L^{2}(w)\to L^{2}(w)).$ Then, if $\widetilde{\J}(\L_{w})$ has more than one point, it is an interval by interpolation; the next proposition further shows that it actually contains a right triangle (see \cite[Figure 1]{DMR}).

\begin{proposition}\label{proposition: cd12}\;  There exists an interval $\J(\L_{w})\subset [1, \infty]$ such that $p, q\in \J(\L_{w})$ if and only if $e^{-t\L_{w}} \in \mho( L^{p}(w)\to L^{q}(w)).$ Furthermore, $\J(\L_{w})$ has the following properties: \begin{equation*}\begin{aligned} (i)&\; \J(\L_{w})\subset \widetilde{\J}(\L_{w}); 
(ii)&\;\mbox{Int}\;\J(\L_{w})=\mbox{Int}\;\widetilde{\J}(\L_{w}); (iii)\; p_{-}(\L_{w})\leq (2^{*, m}_{w})'\; \mbox{and}\; p_{+}(\L_{w})\geq 2^{*, m}_{w},\end{aligned}
\end{equation*}
where $p_{-}(\L_{w})$ and  $p_{+}(\L_{w})$ denote the left and right endpoints of $\J(\L_{w}),$ respectively.  
\end{proposition} 
\begin{remark}\label{remark: cd1085}\; If $w\in A_{1}$ (i.e. $r_{w}=1$), we have $p_{-}(\L_{w})\leq \frac{2n}{n+2m}$ and $p_{+}(\L_{w})\geq \frac{2n}{n-2m}.$ We refer the reader to \cite[Section 8.2]{A} for more precise control over the endpoints $p_{-}(\L_{w})$ and $p_{+}(\L_{w})$ in the case $w\equiv 1$. \end{remark} 

{\it Proof.}\quad We first prove that $e^{-t\L_{w}}\in \mho(L^{2}(w)\to L^{q}(w))$ for any $q$ with $2<q<2^{*, m}_{w}.$ To the end, we need to show (by Definition \ref{definition: cd24}) that \begin{equation}\label{eq: cd9}\l(\fint_{B}e^{-t\L_{w}}(f1_{B})|^{q}dw\r)^{\frac{1}{q}}\lesssim \Upsilon\l(\frac{r}{t^{\frac{1}{2m}}}\r)^{\theta}\l(\fint_{B}|f|^{2}dw\r)^{\frac{1}{2}},\end{equation} \begin{equation}\label{eq: cd1086}\l(\fint_{B}e^{-t\L_{w}}(f1_{C_{j}(B)})|^{q}dw\r)^{\frac{1}{q}}\lesssim 2^{j\theta_{1}}\Upsilon\l(\frac{2^{j}r}{t^{\frac{1}{2m}}}\r)^{m+\theta_{2}}e^{-c\l(\frac{2^{j}r}{t^{\frac{1}{2m}}}\r)^{\frac{2m}{2m-1}}}\l(\fint_{C_{j}(B)}|f|^{2}dw\r)^{\frac{1}{2}},\end{equation} and \begin{equation}\label{eq: cd11}\l(\fint_{C_{j}(B)}e^{-t\L_{w}}(f1_{B})|^{q}dw\r)^{\frac{1}{q}}\lesssim 2^{j\theta_{1}}\Upsilon\l(\frac{2^{j}r}{t^{\frac{1}{2m}}}\r)^{m+\theta_{2}}e^{-c\l(\frac{2^{j}r}{t^{\frac{1}{2m}}}\r)^{\frac{2m}{2m-1}}}\l(\fint_{B}|f|^{2}dw\r)^{\frac{1}{2}}.\end{equation}

We start by proving \eqref{eq: cd9}. Let $g:=e^{-t\L_{w}}(f1_{B}).$ Then, the left-hand side of \eqref{eq: cd9} is controlled by \begin{equation}\label{eq: cd10}\begin{aligned}\l(\fint_{B}e^{-t\L_{w}}(f1_{B})|^{q}dw\r)^{\frac{1}{q}}&\lesssim \l(\fint_{B}|g-Q_{B}g|^{q}dw\r)^{\frac{1}{q}}+\l(\fint_{B}|Q_{B}g-\pi^{m}_{B}g|^{q}dw\r)^{\frac{1}{q}}\\[4pt] 
&+\l(\fint_{B}|\pi^{m}_{B}g|^{q}dw\r)^{\frac{1}{q}} :=J_{1}+J_{2}+J_{3},
 \end{aligned}
 \end{equation} where $Q_{B}g, \pi^{m}_{B}g$ are two polynomials of degree at most $m-1,$ defined in Theorem \ref{theorem: cd0} and Remark \ref{remark: cd1082}. Form \eqref{eq: cd43} and $w\in A_{2}$, 
it follows that 
\begin{equation*}
J_{3}\lesssim \|\pi^{m}_{B}g\|_{L^{\infty}(B)}\lesssim \fint_{B}|g|dx\lesssim  \l(\fint_{B}|g|^{2}dw\r)^{1/2};\end{equation*} furthermore, by Theorem \ref{theorem: cd1073} and Lemma \ref{lemma: h7}, we know $e^{-t\L_{w}} \in \mho( L^{2}(w)\to L^{2}(w)),$ which implies $$\l(\fint_{B}|g|^{2}dw\r)^{1/2} \lesssim \Upsilon\l(\frac{r}{t^{\frac{1}{2m}}}\r)^{\theta_{2}}\l(\fint_{B}|f|^{2}dw\r)^{\frac{1}{2}}.$$ Connecting the two inequalities we reach $$J_{3}\lesssim \Upsilon\l(\frac{r}{t^{\frac{1}{2m}}}\r)^{\theta_{2}}\l(\fint_{B}|f|^{2}dw\r)^{\frac{1}{2}}.$$ To bound $J_{1},$ we apply \eqref{eq: cd1044} along with the property $t^{\frac{1}{2}} \nabla^{m}e^{-t\L_{w}} \in \mho( L^{2}(w)\to L^{2}(w))$ (a result from Theorem \ref{theorem: cd1073} and Lemma \ref{lemma: h7}, as before) to deduce 
$$J_{1} \lesssim r(B)^{m}\l(\fint_{B}|\nabla^{m}e^{-t\L_{w}}(f1_{B}) |^{2}dw\r)^{1/2}\lesssim \frac{r(B)^{m} }{t^{\frac{1}{2}}}\Upsilon\l(\frac{r}{t^{\frac{1}{2m}}}\r)^{\theta_{2}}\l(\fint_{B}|f|^{2}dw\r)^{\frac{1}{2}}.$$ Here we make the simplifying assumption that both operators $e^{-t\L_{w}}$ and $ t^{\frac{1}{2}} \nabla^{m}e^{-t\L_{w}}$ share the exponents $\theta_{1}, \theta_{2}$ from Definition \ref{definition: cd24}. By recalling the definition of $\pi_{B}^{m}$ and the inclusion $A_{2}\subset A_{q},$ we can reduce the estimate of $J_{2}$ to that of $J_{1}:$ $$J_{2}\approx \l(\fint_{B}|\pi^{m}_{B}(Q_{B}g-g)|^{q}dw\r)^{\frac{1}{q}}\lesssim\|\pi^{m}_{B}(Q_{B}g-g)\|_{L^{\infty}(B)} \lesssim J_{1}.$$
Gathering all the above estimates we find \begin{equation*}\begin{aligned}\l(\fint_{B}e^{-t\L_{w}}(f1_{B})|^{q}dw\r)^{\frac{1}{q}}&\lesssim \l(1+\l(\frac{r(B) }{t^{\frac{1}{2m}}}\r)^{m}\r)\Upsilon\l(\frac{r}{t^{\frac{1}{2m}}}\r)^{\theta_{2}}\l(\fint_{B}|f|^{2}dw\r)^{\frac{1}{2}}\\[4pt]
&\lesssim \Upsilon\l(\frac{r}{t^{\frac{1}{2m}}}\r)^{m+\theta_{2}}\l(\fint_{B}|f|^{2}dw\r)^{\frac{1}{2}}.\end{aligned}
\end{equation*} This proves \eqref{eq: cd9}.

An analogous argument results in \eqref{eq: cd1086} and we leave the details to the reader. 

Consider \eqref{eq: cd11} next. For any $j\geq 2,$ the annulus $C_{j}(B)$ can always be covered by a family of balls $\{B_{k}\}_{k=1}^{N}.$ Each ball satisfies $r(B_{k})=2^{j-2} r$ and has its center $x_{k}\in C_{j}(B),$ where the constant $N$ depends solely on $n.$ Repeating the above arguments again and using \eqref{eq: cd66} we can deduce \begin{equation*}\begin{aligned}\l(\fint_{B_{k}}e^{-t\L_{w}}(f1_{B})|^{q}dw\r)^{\frac{1}{q}}&\lesssim \l(\fint_{B_{k}}|e^{-t\L_{w}}(f1_{B}) |^{2}dw\r)^{\frac{1}{2}} +r(B_{k})^{m}\l(\fint_{B_{k}}|\nabla^{m}e^{-t\L_{w}}(f1_{B}) |^{2}dw\r)^{\frac{1}{2}}\\[4pt]
&\lesssim  \l(\fint_{2^{j+1}B\setminus2^{j-1}B }|e^{-t\L_{w}}(f1_{B}) |^{2}dw\r)^{\frac{1}{2}} \\[4pt]
&\quad\quad +(2^{j}r)^{m}\l(\fint_{2^{j+1}B\setminus2^{j-1}B}|\nabla^{m}e^{-t\L_{w}}(f1_{B}) |^{2}dw\r)^{\frac{1}{2}}:=I+II.\end{aligned}
\end{equation*} 
Fix $j\geq 3,$ then $2^{j+1}B\setminus2^{j-1}B=C_{j+1}(B)\cup C_{j}(B)\cup C_{j-1}(B).$ Recall that both $e^{-t\L_{w}}$ and $t^{\frac{1}{2}} \nabla^{m}e^{-t\L_{w}}$ satisfy \eqref{eq: cd27} with $p=q=2$ on each $C_{i}(B)$ for all $i$ satisfying $j-1\leq i\leq j+1.$ Then, we have $$I+II\lesssim 2^{j\theta_{1}}\Upsilon\l(\frac{2^{j}r}{t^{\frac{1}{2m}}}\r)^{m+\theta_{2}}e^{-c\l(\frac{2^{j}r}{t^{\frac{1}{2m}}}\r)^{\frac{2m}{2m-1}}}\l(\fint_{B}|f|^{2}dw\r)^{\frac{1}{2}}.$$ When $j=2,$ we split $2^{4}B\setminus B=C_{3}(B)\cup C_{2}(B)\cup (4B\setminus 2B).$ The preceding arguments extend to $C_{3}(B)$ and $C_{2}(B);$ on $4B\setminus 2B$ we can follow the proof of \cite[Lemma 6.5]{AM}. In summary, it is not difficult to derive  \begin{equation*}\begin{aligned} \l(\fint_{4B\setminus 2B }|e^{-t\L_{w}}(f1_{B}) |^{2}dw\r)^{\frac{1}{2}}&+(2^{2}r)^{m}\l(\fint_{4B\setminus 2B}|\nabla^{m}e^{-t\L_{w}}(f1_{B}) |^{2}dw\r)^{\frac{1}{2}}\\[4pt]
&\lesssim \Upsilon\l(\frac{2r}{t^{\frac{1}{2m}}}\r)^{m+\theta_{2}}e^{-c\l(\frac{2r}{t^{\frac{1}{2m}}}\r)^{\frac{2m}{2m-1}}}\l(\fint_{B}|f|^{2}dw\r)^{\frac{1}{2}}.
\end{aligned}
\end{equation*} 
Summing up these estimates, we arrive at \begin{equation*}\begin{aligned}\l(\fint_{C_{j}(B)}e^{-t\L_{w}}(f1_{B})|^{q}dw\r)^{\frac{1}{q}}&\lesssim  \sum_{k=1}^{N}\l(\fint_{B_{k}}e^{-t\L_{w}}(f1_{B})|^{q}dw\r)^{\frac{1}{q}}\\[4pt] &\lesssim 2^{j\theta_{1}}\Upsilon\l(\frac{2^{j}r}{t^{\frac{1}{2m}}}\r)^{m+\theta_{2}}e^{-c\l(\frac{2^{j}r}{t^{\frac{1}{2m}}}\r)^{\frac{2m}{2m-1}}}\l(\fint_{B}|f|^{2}dw\r)^{\frac{1}{2}}.\end{aligned} 
\end{equation*} This is exactly \eqref{eq: cd11}.

Note that all the estimates just established hold for $\L_{w}^{*}$ due to \eqref{eq: cd3003}. Consequently, $e^{-t\L_{w}^{*}}\in \mho(L^{2}(w)\to L^{q}(w))$ for any $q$ with $2<q<2^{*, m}_{w}.$ Then, by Lemma \ref{lemma: h5}, $e^{-t\L_{w}}\in \mho(L^{q'}(w)\to L^{2}(w)).$ Using this result, along with Lemma \ref{lemma: h6} and the identity $e^{-t\L_{w}}=e^{-t/2\L_{w}}\circ e^{-t/2\L_{w}},$ it holds that $e^{-t\L_{w}}\in \mho(L^{q'}(w)\to L^{q}(w)).$ From this, an argument completely analogous to that in \cite[Proposition 4.1]{AM} yields that there exists an interval $\J(\L_{w})\subset [1, \infty]$ such that $p, q\in \J(\L_{w})$ if and only if $e^{-t\L_{w}} \in \mho( L^{p}(w)\to L^{q}(w)),$ with properties $ (i)$ and $ (ii)$ satisfied. In particular, $[q', q] \subset \J(\L_{w})$ for all $q$ with $2<q<2^{*, m}_{w},$ thereby proving property $ (iii).$


\hfill$\Box$

\begin{corollary}\label{corollary: cd13}\;  Assume $p_{-}(\L_{w})<p\leq q< p_{+}(\L_{w}).$ If $v\in A_{\frac{p}{p_{-}(\L_{w})}}(w)\cap \mbox{RH}_{(\frac{p_{+}(\L_{w})}{q})'}(w), $ then $e^{-t\L_{w}} \in \mho( L^{p}(vdw)\to L^{q}(vdw)).$\end{corollary}
{\it Proof.}\quad Clearly, $e^{-t\L_{w}} \in \mho( L^{p}(w)\to L^{q}(w))$ by Proposition \ref{proposition: cd12}, then Corollary \ref{corollary: cd13} follows instantly from Lemma \ref{proposition: h8}.  

\hfill$\Box$

\begin{corollary}\label{corollary: cd14}\;  For any $\nu$ with $0<\nu<\frac{\pi}{2}-\V$ and any $p\leq q$ such that $e^{-t\L_{w}} \in \mho( L^{p}(w)\to L^{q}(w)),$ we have for all $k\in \N\cup \{0\},$ $(z\L_{w})^{k}e^{-z\L_{w}} \in \mho( L^{p}(w)\to L^{q}(w), \Sigma_{\nu}).$
\end{corollary}
{\it Proof.}\quad Recall that $e^{-z\L_{w}} \in \Fz( L^{2}(w)\to L^{2}(w), \Sigma_{\frac{\pi}{2}-\V})$ by Theorem \ref{theorem: cd1073}. This corollary is a consequence of the characterization of $\J(\L_{w})$ in Proposition \ref{proposition: cd12} and Theorem \ref{theorem: h9}.

 \hfill$\Box$

\section{The weighted $L^{p}$ functional calculus for $\L_{w}$ } 
In Section 3.1, we showed that $\phi(\L_{w})$ is well-defined in $L^{2}(w)$ for any $\phi\in \H^{\infty}(\Sigma_{\mu})$ with $\mu\in (\V, \pi),$ and that it has an $H^{\infty}$ functional calculus as specified in \eqref{eq: cd3006}. However, this result is insufficient for our purpose; we must further define $\phi(\L_{w})$ on $L^{p}(w)$ (and even on $L^{p}(vdw)$) and prove that it satisfies a $L^{p}(w)$-version (and $L^{p}(vdw)$-version) of \eqref{eq: cd3006} to complete the analysis in the subsequent sections. 

\begin{proposition}\label{proposition: cd15}\;  Let $p_{-}(\L_{w})<p<p_{+}(\L_{w})$ and  $\mu\in (\V, \pi).$ There exists a constant $C,$ independent of $\phi$ and $f,$ such that \begin{equation}\label{eq: cd16}\|\phi(\L_{w})f\|_{L^{p}(w)}\leq C\|\phi\|_{\infty}\|f\|_{L^{p}(w)}\end{equation} for any $\phi \in \H_{0}^{\infty}(\Sigma_{\mu});$ that is, $\L_{w}$ has a bounded holomorphic functional calculus on $L^{p}(w).$ If $v\in A_{\frac{p}{p_{-}(\L_{w})}}(w)\cap \mbox{RH}_{(\frac{p_{+}(\L_{w})}{p})'}(w),$ we also have \begin{equation}\label{eq: cd17}\|\phi(\L_{w})f\|_{L^{p}(vdw)}\leq C\|\phi\|_{\infty}\|f\|_{L^{p}(vdw)},\end{equation} with $C$ independent of $\phi$ and $f.$ 
 \end{proposition}
 
\begin{remark}\label{remark: cd91} Although \eqref{eq: cd16} is stated for $\phi\in \H_{0}^{\infty}(\Sigma_{\mu}),$ it in fact holds for all $\phi \in \H^{\infty}(\Sigma_{\mu});$ see \cite{H, M}.
 \end{remark}  
 
{\it Proof.}\quad The proof is quite similar to that in \cite[Proposition 4.3]{DMR}; however we provide the details for the sake of readability. Hereafter, we simplify the notation by setting $p_{-}:=p_{-}(\L_{w})$ and $p_{+}:=p_{+}(\L_{w}).$ 

We first show \eqref{eq: cd16} for any $f\in L_{c}^{\infty}$ when $p\in (p_{-}, 2),$ then prove \eqref{eq: cd17} for $p\in (p_{-}, p_{+});$ notably, \eqref{eq: cd16} will be recovered by taking $v\equiv 1.$ Without loss of generality, we assume $\|\phi\|_{L^{\infty}}=1$ throughout the entire proof. 

We will use Theorem \ref{theorem: cd21} to prove \eqref{eq: cd16} when $p\in (p_{-}, 2).$ To the end, fix $p_{0}$ with $p_{-}<p_{0}<p<2,$ and let $q_{0}=2,$ $\T=\phi(\L_{w}),$ along with the operator $$\A_{r}f(x)=(I-(I-e^{-r^{2m} \L_{w}})^{N})f(x),$$ where $N$ is a sufficiently large integer to be chosen later. Note that $$\A_{r}=\sum_{k=1}^{N}C_{N}^{k}(-1)^{k+1}e^{-kr^{2m} \L_{w}},$$ and that for any $1\leq k\leq N$ and $t, s>0,$ $$\Upsilon\l(\frac{s}{k^{\frac{1}{2m}}}\r)\leq N^{\frac{1}{2m}}\Upsilon(s)\quad \mbox{and}\quad e^{-c\l(\frac{2^{j}r}{(kt)^{\frac{1}{2m}}}\r)^{\frac{2m}{2m-1}}}\leq e^{-\frac{c}{N^{\frac{1}{2m-1}}}\l(\frac{2^{j}r}{t^{\frac{1}{2m}}}\r)^{\frac{2m}{2m-1}}}.$$ As a consequence of Proposition \ref{proposition: cd12}, $$\A_{r} \in \mho (L^{p}(w) \to L^{q}(w)),\quad \forall\; p_{-}<p\leq q<p_{+}.$$ 

We now verify that condition \eqref{eq: cd23} is satisfied for the operators $T=\T, \A_{r}$ and exponents $p_{0}, q_{0}.$ For every ball $B$ with radius $r,$ any $f\in L_{c}^{\infty}$ with $\supp f\subset B$ and $j\geq 1,$ it is easy to see that \begin{equation}\label{eq: cd28}\l(\fint_{B}|\A_{r} (f1_{B})|^{q}dw\r)^{\frac{1}{q}}\lesssim \l(\fint_{B}|f|^{p}dw\r)^{\frac{1}{p}},\end{equation} and for all $j\geq 2,$

\begin{equation}\label{eq: cd29}\l(\fint_{B}|\A_{r}  (f1_{C_{j}(B)})|^{q}dw\r)^{\frac{1}{q}}\lesssim 2^{j\theta_{1}}\Upsilon(2^{j})^{\theta_{2}} e^{-c(2^{j\frac{2m}{2m-1}}})\l(\fint_{C_{j}(B)}|f|^{p}dw\r)^{\frac{1}{p}},\end{equation} and 

\begin{equation}\label{eq: cd30}\l(\fint_{C_{j}(B)}|\A_{r}  (f1_{B})|^{q}dw\r)^{\frac{1}{q}}\lesssim 2^{j\theta_{1}}\Upsilon(2^{j})^{\theta_{2}} e^{-c(2^{j\frac{2m}{2m-1}}})\l(\fint_{B}|f|^{p}dw\r)^{\frac{1}{p}}\end{equation} hold for any $p_{-}<p<q<p_{+}$ and any $1\leq k\leq N.$ Apparently, \eqref{eq: cd30} with $q=q_{0}$ and $p=p_{0}$ implies \eqref{eq: cd23}, where \eqref{eq: cd23} involves the function $g(j):=C2^{j(\theta_{1}+\theta_{2})}e^{-c(2^{j\frac{2m}{2m-1}}})$ satisfying \begin{equation}\label{eq: cd1088}\sum_{j\geq 1} g(j) 2^{jD}<\infty\quad (D\;\mbox{is the doubling constant in \eqref{eq: cd66} }).\end{equation}


Next, we seek to build condition \eqref{eq: cd22}. As $(I-e^{-r^{2m} z})^{N}$ is bounded on $\Sigma_{\frac{\pi}{2}},$ then $\varphi(z):=\phi(z)(I-e^{-r^{2m} z})^{N} \in \H_{0}^{\infty}(\Sigma_{\min\{\mu, \frac{\pi}{2}\}}).$ By \eqref{eq: cd3004}-\eqref{eq: cd3005}, we can write \begin{equation}\label{eq: cd56} \T(I-\A_{r})f=\int_{\Gamma_{+}}e^{-z\L_{w}}f\eta_{+}(z)dz+\int_{\Gamma_{-}}e^{-z\L_{w}}f\eta_{-}(z)dz,\end{equation} where $\Gamma_{\pm}=\rr^{+} e^{\pm i(\frac{\pi}{2}-\theta)}$, $ \eta_{\pm}(z):=\frac{1}{2\pi i}\int_{\gamma_{\pm}}e^{\xi z}\varphi (\xi) d\xi,$ $\gamma_{\pm}:=\rr^{+}e^{\pm i\nu}$ and $0<\V<\theta<\nu<\min\{\mu, \frac{\pi}{2}\}.$ Utilizing the mean value inequality, a straightforward calculation gives \begin{equation}\label{eq: cd31} |\eta_{\pm}(z)|\lesssim \frac{r^{2mN}}{|z|^{N+1}}.\end{equation} By Corollary \ref{corollary: cd14} and the definition of $\Gamma_{\pm},$ $e^{-z\L_{w}} \in \mho (L^{p_{0}}(w) \to L^{p_{0}}(w))$ for any $z\in \Gamma_{\pm}.$ Therefore, for every ball $B$ of radius $r,$ any $f\in L_{c}^{\infty}$ with $\supp f\subset B$ and $j\geq 2,$ if we choose $N$ large enough such that $2mN>\theta_{2}+1,$ then \begin{equation}\label{eq: cd32} \begin{aligned}&\l(\fint_{C_{j}(B)}|\T(I-\A_{r})f|^{p_{0}}dw\r)^{\frac{1}{p_{0}}}\\[4pt]&\quad\quad\quad  \lesssim \l(\fint_{C_{j}(B)}\bigg|\int_{\Gamma_{\pm}}e^{-z\L_{w}}f\eta_{\pm}(z)dz\bigg|^{p_{0}}dw\r)^{\frac{1}{p_{0}}}\\[4pt]
&\quad\quad\quad \lesssim\int_{\Gamma_{\pm}} \l(\fint_{C_{j}(B)} |e^{-z\L_{w}}f|^{p_{0}}dw\r)^{1/p_{0}} \frac{r^{2mN}}{|z|^{N+1}} |dz| \\[4pt]
&\quad\quad\quad\lesssim \l(\fint_{B}|f|^{p_{0}}dw\r)^{1/p_{0}}\int_{\Gamma_{\pm}}2^{j\theta_{1}}\frac{r^{2mN}}{|z|^{N+1}}\Upsilon\l(\frac{2^{j} r}{|z|^{1/2m}}\r)^{\theta_{2}}e^{-c\l(\frac{2^{j}r}{|z|^{\frac{1}{2m}}}\r)^{\frac{2m}{2m-1}}}|dz|\\[4pt]
&\quad\quad\quad\approx \l(\fint_{B}|f|^{p_{0}}dw\r)^{1/p_{0}} 2^{j(\theta_{1}-2mN)}\int_{0}^{\infty}\Upsilon(\tau)^{\theta_{2}} \tau^{2mN} e^{-c\tau^{\frac{2m}{2m-1}} }\frac{d\tau}{\tau}\\[4pt]
&\quad\quad\quad\lesssim2^{j(\theta_{1}-2mN)} \l(\fint_{B}|f|^{p_{0}}dw\r)^{1/p_{0}}.\end{aligned}
\end{equation} Further imposing $2mN>1+\theta_{1}+\theta_{2}+D,$ we have \eqref{eq: cd1088} satisfied with $g(j)=C 2^{j(\theta_{1}-2mN)}.$ Invoking Theorem \ref{theorem: cd21}, it follows that \eqref{eq: cd16} holds for all $p_{-}<p\leq 2.$ 

We now establish \eqref{eq: cd17} for $p\in (p_{-}, p_{+})$ by applying Theorem \ref{theorem: cd18}. 
Since $v\in A_{\frac{p}{p_{-}}}(w)\cap \mbox{RH}_{(\frac{p_{+}}{p})'}(w),$ by Proposition \ref{proposition: cd3000}, there are $p_{0}, q_{0}$ (by letting $p_{0}\to p_{-}$ and $q_{0}\to p_{+}$) such that \begin{equation}\label{eq: cd5005}p_{-}<p_{0}<\min\{p, 2\}\leq p<q_{0}<p_{+}\quad \mbox{and}\quad A_{\frac{p}{p_{0}}}(w)\cap \mbox{RH}_{(\frac{q_{0}}{p})'}(w).\end{equation} In the sequel, we let the operator $\S$ in Theorem \ref{theorem: cd18} be the identity operator $I.$ Recall that $\T$ is bounded on $L^{p_{0}}(w),$ as established in the preceding argument. To apply Theorem \ref{theorem: cd18}, it remains to verify consitions \eqref{eq: cd19}- \eqref{eq: cd20} for the operators $\T$ and $\S.$ 

Given a ball $B$ of radius $r,$ decompose $f$ as $f=\sum_{j=1}^{\infty}f1_{C_{j}(B)}:=\sum_{j=1}^{\infty}f_{j}.$ A similar argument as in \eqref{eq: cd32} contributes to, for all $f\in L_{c}^{\infty},$ 
\begin{equation}\label{eq: cd33}\begin{aligned}
\l(\fint_{B}|\phi(\L_{w})(I-\A_{r})f|^{p_{0}}dw\r)^{\frac{1}{p_{0}}}&\lesssim \sum_{j\geq 1}\l(\fint_{B}|\phi(\L_{w})(I-\A_{r})f_{j}|^{p_{0}}dw\r)^{\frac{1}{p_{0}}} \\[4pt]
&\lesssim \sum_{j\geq 1}2^{j(\theta_{1}-2mN)}\l(\fint_{C_{j}(B)}|\S f|^{p_{0}}dw\r)^{1/p_{0}}\\[4pt]\end{aligned} 
\end{equation} with the restriction $2mN>\theta_{2}+1.$ This thus leads to \eqref{eq: cd19} with $g(j):=C2^{j(\theta_{1}-2mN)}.$ Here, the series $\sum_{j}g(j)<\infty$ converges, provided we choose $N$ such that $2mN>\theta_{1}+\theta_{2}+2.$ 

Exploiting the commutativity of $\T$ and $\A_{r},$ together with \eqref{eq: cd28}-\eqref{eq: cd29}, we can deduce 
\begin{equation}\label{eq: cd34}\begin{aligned}
\l(\fint_{B}|\T\A_{r}f|^{q_{0}}dw\r)^{\frac{1}{q_{0}}}&\lesssim \l(\fint_{B}|\A_{r}\T f|^{q_{0}}dw\r)^{\frac{1}{q_{0}}} \\[4pt]&\lesssim \sum_{j\geq 1}\l(\fint_{B}|\A_{r}[(\T f)_{j}]|^{p_{0}}dw\r)^{\frac{1}{p_{0}}} \\[4pt]&\lesssim \sum_{j\geq 1} 2^{j\theta_{1}}\Upsilon(2^{j})^{\theta_{2}} e^{-c(2^{j\frac{2m}{2m-1}})}\l(\fint_{C_{j}(B)}|\T f|^{q_{0}}dw\r)^{1/q_{0}}\\[4pt]
&\lesssim \sum_{j\geq 1}2^{j(\theta_{1}+\theta_{2})}e^{-c(2^{j\frac{2m}{2m-1}})}\l(\fint_{2^{j+1}B}|\T f|^{p_{0}}dw\r)^{1/p_{0}},
 \end{aligned} 
\end{equation} which implies \eqref{eq: cd20} due to $\sum_{j\geq 1}2^{j(\theta_{1}+\theta_{2})}e^{-c(2^{j\frac{2m}{2m-1}})}<\infty.$ Thus, Theorem \ref{theorem: cd18} applies, and the proof of \eqref{eq: cd17} is complete.

Lastly, as $L_{c}^{\infty}$ is dense both in $L^{p}(w)$ and $L^{p}(vdw)$, \eqref{eq: cd16} and \eqref{eq: cd17} extends to $L^{p}(w)$ and $L^{p}(vdw),$ respectively, via a limiting argument. 

\hfill$\Box$

\section{Reverse inequalities for square roots in weighted spaces} Building on the preparations in the previous sections, we are now in a position to identify the intervals for which the reverse square root inequalities (cf. \eqref{eq: cd2001}) are satisfied. The endpoints of these intervals will depend on the exponents $p_{-}, p_{+}$ and $ r_{w},$ owing to the reliance of our proof on the generalized Poincar\'{e}-Sobolev inequalities (Theorem \ref{theorem: cd0} and Remark \ref{remark: cd1082}), the off-diagonal estimates for the semigroup $e^{-z\L_{w}}$ (Corollary \ref{corollary: cd14}) and the $H^{\infty}$ functional calculus (Proposition \ref{proposition: cd15}).

Prior to proving \eqref{eq: cd2001}, two technical lemmas are needed. The first one is a higher-order generalization of the weighted Calder$\acute{o}$n-Zygmund decomposition from \cite[lemma 6.6]{AM1}, and additionally constitutes a weighted extension of \cite[Lemma 16]{A1}.  

\subsection{The higher-order weighted Calder$\acute{o}$n-Zygmund decomposition } 


\begin{lemma}\label{lemma: cd36}\;  Given $w\in A_{p}$ with $1\leq p<\infty.$ Assume that $f\in \ss(\rz)$ 
satisfies $\|\nabla^{m} f\|_{L^{p}_{w}(\rz)}<\infty.$ Fix $\alpha>0.$ Then there exist a collection of cubes $\{Q_{i}\}$ (or balls $\{B_{i}\}$), functions $g\in L^{1}_{loc}(w)$ and $b_{i}$ such that 

\begin{equation}\label{eq: cd37}f=g+\sum_{i}b_{i},\end{equation} and the following properties hold: 

\begin{equation}\label{eq: cd38}\|\nabla^{m}g\|_{\infty}\leq C \alpha,\end{equation} \begin{equation}\label{eq: cd39}b_{i}\in W^{m, p}_{0}(Q_{i}) \quad \mbox{and}\quad \int_{Q_{i}}|\nabla^{m} b_{i}|^{p}dw\leq C \alpha^{p} w(Q_{i}),\end{equation} 

\begin{equation}\label{eq: cd40}\sum_{i}w(Q_{i})\leq \frac{C}{\alpha^{p}} \int_{\rz}|\nabla^{m} f|^{p}dw,\end{equation} 

\begin{equation}\label{eq: cd41}\sum_{i}1_{Q_{i}}\leq \M,\end{equation} and for all $1\leq q <p^{*, m}_{w},$ \begin{equation}\label{eq: cd42}\l(\fint_{Q_{i}}|b_{i}|^{q}dw\r)^{1/q}\leq C \alpha l(Q_{i})^{m},\end{equation}
where $C$ and $\M$ depends only on $p, q, m,$ the doubling constant of $w$ and dimension. 

\end{lemma} 
{\it Proof.}\quad We define the uncentered maximal operator $M_{w}$ with respect to the wight $w$ as follows: $$M_{w}f(x):=\sup_{x\ni Q}\fint_{Q}|f(x)|dw.$$ Let $\Omega:=\{x\in \rz: M_{w}(|\nabla^{m}f|^{p})(x) >\alpha^{p}\}.$ If $\Omega$ is empty, we may directly define $g$ to be equal to $f.$ Since $w\in A_{p},$ then $dw$ is doubling (see \eqref{eq: cd66} ). By the maximal theorem, this implies $$w(\Omega)\leq \frac{C}{\alpha^{p}} \int_{\rz}|\nabla^{m} f|^{p}dw. $$ In the sequel, we denote the complement of $\Omega$ by $F.$ By the Lebesgue differentiation theorem, we readily obtain that $$|\nabla^{m}f(x)|\leq C \alpha,\quad \mbox{for}\; dw-a.e. \;x \in F.$$

To continue, we decompose $\Omega$ into a collection of dyadic Whitney boxes $\{Q_{i}\}.$ This decomposition satisfies three key properties: $\Omega$ is the disjoint union of the $Q_{i};$ each $Q_{i}$ satisfies $2Q_{i}\subset \Omega;$ the family $\{Q_{i}\}$ has bounded overlap and every cube $4Q_{i}$ intersects $F.$ Furthermore,  \begin{equation}\label{eq: cd47}\mbox{if}\; Q_{i}\cap Q_{j} \neq \varnothing,\;\mbox{then}\; l(Q_{i})\approx l(Q_{j})\;\mbox{ and}\; |z-y|\leq C l(Q_{j})\;\mbox{ for any} \; z\in Q_{i}, \; y\in Q_{j}.\end{equation} Using this decomposition and the aforementioned two inequalities, \eqref{eq: cd40}-\eqref{eq: cd41} for the cubes $2Q_{i}$ follow directly. 

For the proof of \eqref{eq: cd39}, we consider a sequence of smooth functions with compact supports $\{\eta_{i}\}$, induced by the partition of unity on $\Omega$ for the covering $\{Q_{i}\}.$ Clearly, $\supp\; \eta_{i} \subset 2Q_{i}$ with the estimate $$l(Q_{i})^{|\gamma|}\|D^{\gamma} \eta_{i}\|_{\infty}\leq C$$ holds for all $ |\gamma|\leq m.$ If we define $$b_{i}=(f-\pi_{2Q_{i}}^{m}f)\eta_{i},$$ then $\supp \;b_{i} \subset 2Q_{i}.$ Moreover, by the Leibniz rule and \eqref{eq: cd44}, for all $|\gamma|\leq m$ we derive \begin{equation}\label{eq: cd45}\begin{aligned}
\|D^{\gamma} b_{i}\|_{L^{p}_{w}(\rz)}&\lesssim \sum_{\beta\leq \gamma}C_{\beta}^{\gamma}l(Q_{i})^{-(|\gamma|-|\beta|)}\|D^{\beta}(f-\pi_{2Q_{i}}^{m}f)\|_{L^{p}_{w}(2Q_{i})} \\[4pt]
&\lesssim l(Q_{i})^{(m-|\gamma|)}\|\nabla^{m}f \|_{L^{p}_{w}(2Q_{i})}.\\[4pt] \end{aligned}\end{equation} This yields \eqref{eq: cd39} becasue $4Q_{i}\cap F$ is nonempty.

It remains to prove \eqref{eq: cd37}-\eqref{eq: cd38}. First, we show that $\sum_{i}b_{i}$ converges in  $L^{p}_{loc}(\rz, dw).$ Indeed, fix a compact set $E\subset \rz,$ then the cubes $Q_{i}$ that intersect $E$ have uniformly bounded sidelengths. 
From \eqref{eq: cd45}, together with the bounded overlap property of the $Q_{i}$'s, it follows that $$\sum_{i}\|b_{i}\|_{L^{p}_{w}(E)}\lesssim \sum_{i,\;Q_{i}\cap E\neq \varnothing }l(Q_{i})^{m} \|\nabla^{m}f \|_{L^{p}_{w}(\rz)}<\infty.$$ This ensures that $$g:=f-\sum_{i}b_{i} \quad dw-a.e. \;x$$ is well-defined. Second, as \eqref{eq: cd39}-\eqref{eq: cd41} imply the convergence of $\sum_{i}|\nabla^{m} b_{i}|$ in $L^{p}_{w}(\rz),$ we have \begin{equation}\label{eq: cd46}\nabla^{m}g=\nabla^{m}f-\sum_{i}\nabla^{m}b_{i} \quad dw-a.e.\;x.\end{equation}  From the equality, our goal is to compute $\nabla^{m} g$ in order to deduce \eqref{eq: cd38}.

Given that $\sum_{i} \eta_{i}=1$ on $\Omega,$ and $\sum_{i} \eta_{i}=$ on $F,$ and the sum is locally finite, we have $$\sum_{i} D^{\gamma}\eta_{i}=0\quad \mbox{on}\; \Omega \quad \mbox{for any}\; 1\leq |\gamma|\leq m.$$ Then, for any $|\gamma|=m,$ applying Leibniz's rule and the aforementioned estimate, we arrive at $$\sum_{i}D^{\gamma} b_{i}=D^{\gamma}f \sum_{i}\eta_{i}+\sum_{i}\sum_{\beta<\gamma}C_{\beta, \gamma}D^{\beta}\pi_{2Q_{i}}^{m}f D^{\gamma-\beta}\eta_{i}.$$ To bound the abve two sums, we introduce the notation $h:=h_{\beta, \gamma}$ with $$h_{\beta, \gamma}:=\sum_{i}D^{\beta}\pi_{2Q_{i}}^{m}f D^{\gamma-\beta}\eta_{i}.$$ Observing that, if \begin{equation}\label{eq: cd1089}\|h\|_{\infty}\leq C \alpha \;\mbox{for any}\; \beta, \gamma,\end{equation} then, by \eqref{eq: cd46}, we see that $$D^{\gamma} g=(D^{\gamma} f)1_{F}-\sum_{\beta<\gamma}C_{\beta, \gamma}h_{\beta, \gamma}$$ holds almost everywhere.\footnote{Bear in mind that if $w\in A_{p},$ then for any measurable set $E\subset \rz,$ $w(E)=0$ if and only if $|E|=0$.} From this, \eqref{eq: cd38} is deduced.  

Let us turn to the proof of \eqref{eq: cd1089}. Note that the sum defining $h$ is locally finite in $\Omega,$ with $h(x)=0$ whenever $x\in F.$ If $Q_{j}$ is the Whitney cube containing $x\in \Omega$ and $I_{x}$ denotes the set of indices $i$ such that $x\in 2Q_{i},$ then $\sharp I_{x}\leq \M.$ Choose $x_{j}\in 4Q_{j}\cap F,$ and let $\widetilde{Q}_{j}$ be a dilation of $Q_{j}$ that contains all cubes $2Q_{i}$ for $i\in I_{x}$ (as guaranteed by \eqref{eq: cd47}) and the point $x_{j}$. As $\gamma-\beta\neq 0,$ we may write $$h(x)=\sum_{i\in I_{x}}D^{\beta}(\pi_{2Q_{i}}^{m}f-\pi_{\widetilde{Q}_{j}}^{m}f)(x) D^{\gamma-\beta}\eta_{i}(x).$$ Then, there exists a constant $C,$ independent of $x$ and $f,$ such that \begin{equation}\label{eq: cd48}
J:=|D^{\beta}(\pi_{2Q_{i}}^{m}f-\pi_{\widetilde{Q}_{j}}^{m}f)(x)|\leq Cl(Q_{j})^{m-|\beta|}\l(\fint_{\widetilde{Q}_{j}}|\nabla^{m} f|^{p}dw\r)^{1/p}.
\end{equation} Admit \eqref{eq: cd48} for the moment. We then have the following estimate:  
\begin{equation}\label{eq: cd49}\begin{aligned}
|h(x)|&\leq C\sum_{i\in I_{x}} l(Q_{j})^{m-|\beta|}\l(\fint_{\widetilde{Q}_{j}}|\nabla^{m} f|^{p}dw\r)^{1/p}
l(Q_{i})^{-(|\gamma|-|\beta|)} \\[4pt]
&\leq Cl(Q_{j})^{(m-|\gamma|)}\l(\fint_{\widetilde{Q}_{j}}|\nabla^{m} f|^{p}dw\r)^{1/p}\leq C M_{w}(|\nabla^{m} f|^{p})^{1/p} \leq C \alpha, 
\end{aligned}\end{equation} which contributes to \eqref{eq: cd1089} due to the arbitrariness of $x$ in $\Omega.$ Therefore, the proof will be complete once we establish \eqref{eq: cd48}. Employing \eqref{eq: cd72}, it follows that  $\pi_{2Q_{i}}^{m}(\pi_{\widetilde{Q}_{j}}^{m}f)=\pi_{\widetilde{Q}_{j}}^{m}f.$ Utilizing this and \eqref{eq: cd43}-\eqref{eq: cd44}, we can derive
\begin{equation*}\begin{aligned}
J\leq C \|D^{\beta}\pi_{2Q_{i}}^{m}(f-\pi_{\widetilde{Q}_{j}}^{m}f)\|_{L^{\infty}(2Q_{i})}&\leq Cl(Q_{i})^{-n}\int_{2Q_{i}}|D^{\beta} (f-\pi_{\widetilde{Q}_{j}}^{m}f)|dx\\[4pt] 
&\leq C  l(Q_{j})^{-n}\int_{\widetilde{Q}_{j}}|D^{\beta} (f-\pi_{\widetilde{Q}_{j}}^{m}f)|dx\\[4pt] 
&\leq C  \l(\fint_{\widetilde{Q}_{j}}|D^{\beta} (f-\pi_{\widetilde{Q}_{j}}^{m}f)|^{p}dw\r)^{1/p}\\[4pt] 
&\leq C l(Q_{j})^{m-|\beta|} \l(\fint_{\widetilde{Q}_{j}}|\nabla^{m} f|^{p}dw\r)^{1/p}.\end{aligned}\end{equation*} This suffices.

\hfill$\Box$

\subsection{The weighted conservation property in higher-order case} The second technical lemma concerns a conservation property for higher-order weighted elliptic operators. Its proof generalizes the arguments found in \cite[Lemma 3.1]{AHMT1} and \cite[Section 3.5]{A}. 

\begin{lemma}\label{lemma: cd67}\;  Let $w\in A_{2}.$  Then for every polynomial $P$ with degree $d$ not exceeding $m-1,$ the equality $$e^{-t\L_{w}}P=P$$ holds in the sense of $L^{2}_{loc}(w).$ 
\end{lemma}
{\it Proof.}\quad Let $\eta\in C_{0}^{\infty}(B_{2}(0))$ such that $\eta\equiv 1$ on $B_{1}(0).$ For $R>0,$ define $\eta_{R}(x):=\eta(x/R)$  For any $\phi\in C_{0}^{\infty}(\rz),$ and for all $t>0$ and sufficiently large $R$, we decompose the integral as \begin{equation}\label{eq: cd68} \int_{\rz}P(x) \overline{e^{-t\L_{w}^{*}} \phi} dw(x)=\int_{\rz}P\eta_{R}\overline{e^{-t\L_{w}^{*}} \phi}dw(x) +\int_{\rz}P(1-\eta_{R})\overline{e^{-t\L_{w}^{*} }\phi}dw(x):=I+II.\end{equation} 

The integral $I$ is well-defined thanks to $P\eta_{R}\in L^{2}(w)$ and $e^{-t\L_{w}^{*}}\in \mho(L^{2}(w)\to L^{2}(w))$ by Proposition \ref{proposition: cd12} and Lemma \ref{lemma: h5}. On the other hand, an application of Lemma \ref{lemma: h7} shows that $e^{-t\L_{w}^{*}}\in \Fz (L^{2}(w)\to L^{2}(w)).$ Choosing $R$ large enough such that $\supp\;\phi \subset B_{R}(0)$ and applying \eqref{eq: cd66},  we can bound the integral $II$ as follows:  
\begin{equation}\label{eq: cd69}\begin{aligned}
II &\lesssim \sum_{j\geq 0}\int_{C_{j}(B_{R})}|P(x)| |e^{-t\L_{w}^{*}} \phi|dw(x)\\[4pt] 
&\lesssim \sum_{j\geq 0} (2^{j} R)^{d}\l(\int_{C_{j}(B_{R})} |e^{-t\L_{w}^{*}} \phi|^{2}dw(x)\r)^{1/2}w(B_{2^{j} R})^{1/2}\\[4pt] 
&\lesssim \sum_{j\geq 0} (2^{j} R)^{d} (2^{j} R)^{D/2} w(B_{1})^{1/2}e^{-c\l(\frac{2^{j}R}{t^{\frac{1}{2m}}}\r)^{\frac{2m}{2m-1}}}\|\phi1_{\supp\; \phi}\|_{L^{2}(w)} \\[4pt] 
& \lesssim \sum_{j\geq 0} (2^{j} R)^{d} (2^{j} R)^{D/2} w(B_{1})^{1/2} e^{-c'(2^{j\frac{2m}{2m-1} })}  e^{-c''\l(\frac{R}{t^{\frac{1}{2m}}}\r)^{\frac{2m}{2m-1}}}\|\phi 1_{\supp\; \phi}\|_{L^{2}(w)} \\[4pt] 
&\leq C (t) \|\phi 1_{\supp\; \phi}\|_{L^{2}(w)}w(B_{1})^{1/2}<\infty.\end{aligned}\end{equation} Thus the equality \eqref{eq: cd68} makes sense. 

Note also that $t\L^{*}_{w} e^{-t\L_{w}^{*}} \in \Fz (L^{2}(w)\to L^{2}(w))$ by Corollary \ref{corollary: cd14}, Lemma \ref{lemma: h5} and Lemma \ref{lemma: h7}. From a similar argument to \eqref{eq: cd69}, it holds that $$\int_{\rz}P(1-\eta_{R})\overline{\frac{d }{dt} e^{-t\L_{w}^{*} }\phi}dw(x)= \frac{d }{dt} \int_{\rz}P(1-\eta_{R})\overline{ e^{-t\L_{w}^{*} }\phi}dw(x).$$ This expression is well-defined and tends to zero as $R\to \infty.$ Furthermore, using the definition of $\L_{w}^{*}$ (\eqref{eq: cd3003}) and the Leibniz rule, we get 
\begin{equation*}\begin{aligned}\frac{d }{dt} \int_{\rz}P\eta_{R}\overline{ e^{-t\L_{w}^{*} }\phi}dw(x)&=\sum_{|\alpha|=|\beta|=m}\int_{\rz}w^{-1}a_{\alpha, \beta}\partial^{\beta}(P\eta_{R})\overline{\partial^{\alpha} e^{-t\L_{w}^{*} }\phi}dw(x)\\[4pt] 
&=\sum_{|\alpha|=|\beta|=m}\sum_{\gamma<\beta}C_{\beta}^{\gamma}\int_{\rz}w^{-1}a_{\alpha, \beta}\partial^{\gamma}P \partial^{\beta-\gamma}\eta_{R}\overline{ \partial^{\alpha}e^{-t\L_{w}^{*} }\phi}dw(x):=\M, \end{aligned}\end{equation*} where in the last step we also used the fact that the degree of the polynomial $P$ is less than $m.$ Since $|\beta-\gamma|\geq 1$ and $\supp\; (\partial^{\beta-\gamma}\eta_{R}) \subset B_{2R}\setminus B_{R},$ we obtain 
\begin{equation}\label{eq: cd70}\begin{aligned}
\M&\lesssim \sum_{|\alpha|=|\beta|=m}\sum_{\gamma<\beta}C_{\beta}^{\gamma}R^{d-|\gamma|}R^{-(m-|\gamma|)}\int_{B_{2R}\setminus B_{R}}|\partial^{\alpha}e^{-t\L_{w}^{*} }\phi| dw(x)\\[4pt] 
&\lesssim R^{d-m} w(B_{2R})^{\frac{1}{2}} \l(\int_{B_{2R}\setminus B_{R}} |\nabla^{m}e^{-t\L_{w}^{*}} \phi|^{2}dw(x)\r)^{\frac{1}{2}}\;(t^{1/2}\nabla^{m}e^{-t\L_{w}^{*}}\in \Fz(L^{2}(w)\to L^{2}(w)))\\[4pt] 
&\lesssim R^{d-m} w(B_{2})^{1/2} t^{-1/2}R^{D/2} e^{-c\l(\frac{R}{t^{\frac{1}{2m}}}\r)^{\frac{2m}{2m-1} } } \|\phi 1_{\supp\; \phi}\|_{L^{2}(w)}w(B_{1})^{1/2},
\end{aligned}\end{equation} which tends to zero as  $R\to \infty.$ Putting all these estimates together, we conclude that the left hand side of \eqref{eq: cd68} is independent of $t>0.$

To conclude the proof of Lemma \ref{lemma: cd67}, it suffices to show that \begin{equation}\label{eq: cd71}\int_{\rz}P(x) \overline{e^{-t\L_{w}^{*}} \phi} dw(x)=\int_{\rz} P\overline{\phi}dw(x)\end{equation} for all compactly supported $\phi\in L^{2}(w).$ Choose $R$ large enough so that the supports of $\phi$ and $(1-\eta_{R})$ are far apart. Exploiting a similar argument to \eqref{eq: cd69} we obtain
\begin{equation*}\begin{aligned}II &\lesssim \sum_{j\geq 0} (2^{j} R)^{d} (2^{j} R)^{D/2} w(B_{1})^{1/2}e^{-c\l(\frac{2^{j}R}{t^{\frac{1}{2m}}}\r)^{\frac{2m}{2m-1}}}\|\phi1_{\supp\; \phi}\|_{L^{2}(w)} \\[4pt] 
& \lesssim \sum_{j\geq 0} (2^{j} R)^{d} (2^{j} R)^{D/2} w(B_{1})^{1/2} e^{-c'(2^{j}R )^{\frac{2m}{2m-1} } } e^{-c''\l(\frac{1}{t^{\frac{1}{2m}}}\r)^{\frac{2m}{2m-1}}}\|\phi 1_{\supp\; \phi}\|_{L^{2}(w)} \\[4pt] 
& \lesssim e^{-c''\l(\frac{1}{t^{\frac{1}{2m}}}\r)^{\frac{2m}{2m-1}}}\|\phi 1_{\supp\; \phi}\|_{L^{2}(w)} w(B_{1})^{1/2},\end{aligned}\end{equation*} where the right hand side tends to 0 as $t \to 0$. 
In addition, because $e^{-t\L_{w}^{*}}$ forms a continuous semigroup on $L^{2}(w)$ at $t=0,$ it follows that $$I\to \int_{\rz} P\eta_{R} \overline{\phi}dw(x)=\int_{\rz} P\overline{\phi}dw(x)\quad \mbox{as}\; t\to 0.$$ Combining this with \eqref{eq: cd68}, we arrive at \eqref{eq: cd71}.

\hfill$\Box$

\subsection{Proof of the reverse inequalities in $L^{p}(w)$ and $L^{p}(vdw)$} We are now ready to present the proof of \eqref{eq: cd2001}. 
Define $$(p_{-})_{w, m, *}:=\frac{nr_{w}p_{-}}{nr_{w}+mp_{-}}.$$ Clearly, $(p_{-})_{w, m, *}<p_{-}<2.$

\begin{proposition}\label{proposition: cd81}\;  Let $\max\{r_{w}, (p_{-})_{w, m, *}\}<p<p_{+}$. Then for all $f\in \ss(\rz),$ \begin{equation}\label{eq: cd82} \|\L_{w}^{1/2}f\|_{L^{p}(w)} \leq C \|\nabla^{m}f\|_{L^{p}(w)},\end{equation}  furthermore, if $$\max\{r_{w}, p_{-}\}<p<p_{+}\quad \mbox{and}\quad v\in A_{\frac{p}{\max\{r_{w}, p_{-}\}   }}(w)\cap \mbox{RH}_{(\frac{p_{+}}{p})'}(w),$$ then \begin{equation}\label{eq: cd83} \|\L_{w}^{1/2}f\|_{L^{p}(vdw)} \leq C \|\nabla^{m}f\|_{L^{p}(vdw)},\end{equation} where the constant $C$ is independent of $f.$ \end{proposition}
{\it Proof.}\quad Our argument proceeds along the same lines as the proof of \cite[Proposition 6.1]{DMR}. Given $p$ with $\max\{r_{w}, (p_{-})_{w, m, *}\}<p<2,$ and $f\in \ss(\rz).$ Our first objective is to establish \begin{equation}\label{eq: cd84} \|\L_{w}^{1/2}f\|_{L^{p, \infty}(w)} \lesssim \|\nabla^{m}f\|_{L^{p}(w)}.\end{equation} Note that $w\in A_{p}\subset A_{2}$ if $2>p>r_{w},$ by the definition of $r_{w}.$ Then, from \eqref{eq: cd37} in Lemma \ref{lemma: cd36}, it suffices to build the corresponding weak-type estimates in \eqref{eq: cd84} with $f$ replaced by $g$ and $b_{i}.$


For $g,$ using successively the $L^{2}(w)$-Kato estimate \eqref{eq: cd1081}, \eqref{eq: cd38}, \eqref{eq: cd41} and \eqref{eq: cd39}-\eqref{eq: cd40}, we can derive
\begin{equation*}\begin{aligned}w(\{|\L_{w}^{1/2}g |>\alpha/3\}) &\lesssim \frac{1}{\alpha^{2}}\int_{\rz} |\nabla^{m} g|^{2}dw \lesssim  \frac{1}{\alpha^{p}}\int_{\rz} |\nabla^{m} g|^{p}dw\\[4pt] 
& \lesssim  \frac{1}{\alpha^{p}}\int_{\rz} |\nabla^{m} f|^{p}dw + \frac{1}{\alpha^{p}}\int_{\rz} |\sum_{i}\nabla^{m} b_{i}|^{p}dw \lesssim  \frac{1}{\alpha^{p}}\int_{\rz} |\nabla^{m} f|^{p}dw.
\end{aligned}\end{equation*} For $b_{i}$ with $\supp\; b_{i}\subset B_{i}$, we first observe that there is a $k\in \Z$ such that $2^{k}\leq r(B_{i})<2^{k+1}.$ Then, for all $i,$ $r_{i}\approx r(B_{i})$ if we let $r_{i}=2^{k}.$ By virtue of \cite{H, TK, M}, the square root $\L_{w}^{1/2}$ has the integral representation: \begin{equation}\label{eq: cd1026}\L_{w}^{1/2}=\frac{1}{\sqrt{\pi}}\int_{0}^{\infty} t^{1/2} \L_{w} e^{-t\L_{w}}\frac{d t}{t}.\end{equation} Thus, we can write $$\L_{w}^{1/2}=\frac{1}{\pi^{1/2}}\int_{0}^{r_{i}^{2m}}\L_{w}e^{-t\L_{w}}\frac{dt}{t^{1/2}}+\frac{1}{\pi^{1/2}}\int_{r_{i}^{2m}}^{\infty}\L_{w}e^{-t\L_{w}}\frac{dt}{t^{1/2}}:=T_{i}+S_{i}.$$
Then, by \eqref{eq: cd40},
\begin{equation*}\begin{aligned}w(\{|\sum_{i}\L_{w}^{1/2}b_{i} |>2\alpha/3\}) &\lesssim w(\cup_{i} 4B_{i})+ w(\{|\sum_{i} S_{i} b_{i}|>\alpha/3\})\\[4pt] 
& \quad +  w((\cup_{i} 4B_{i} )^{c}\cap \{|\sum_{i} T_{i} b_{i}|>\alpha/3\})\\[4pt] 
& \lesssim  \frac{1}{\alpha^{p}}\int_{\rz} |\nabla^{m} f|^{p}dw + J_{1}+J_{2}.
\end{aligned}\end{equation*} where $$J_{1}:= w(\{|\sum_{i} S_{i} b_{i}|>\alpha/3\})\quad\mbox{and}\quad J_{2}:=w((\cup_{i} 4B_{i} )^{c}\cap \{|\sum_{i} T_{i} b_{i}|>\alpha/3\}).$$ First, we bound $J_{2}.$ Because $p>(p_{-})_{w, m, *},$ this implies $p^{*, m}_{w}>((p_{-})_{w, m, *}) ^{*, m}_{w}=p_{-}.$ As a consequence, there exists a $q\in \J(\L_{w})$ for which \eqref{eq: cd42} holds. Moreover, applying Corollary \ref{corollary: cd14} to this exponent $q,$ we further have $t\L_{w}e^{-t\L_{w}}\in \mho(L^{q}(w) \to L^{q}(w))$. By this property, \eqref{eq: cd66}, \eqref{eq: cd42} and \eqref{eq: cd40}, we derive

\begin{equation*}\begin{aligned} J_{2} &\lesssim \frac{1}{\alpha}\sum_{i}\sum_{j\geq 2}\int_{C_{j}(B_{i})}|T_{i} b_{i}|d w \\[4pt] 
&\lesssim \frac{1}{\alpha}\sum_{i}\sum_{j\geq 2}2^{jD}w( B_{i})\int_{0}^{r_{i}^{2m}}\l(\fint_{C_{j}(B_{i})}|t\L_{w}e^{-t\L_{w}} b_{i}|^{q}d w\r)^{1/q}\frac{dt}{t^{3/2}}\\[4pt]& \lesssim \frac{1}{\alpha}\sum_{i}\sum_{j\geq 2}2^{jD}w( B_{i})\int_{0}^{r_{i}^{2m}} 2^{j\theta_{1}}\Upsilon\l(\frac{2^{j}r_{i}}{t^{1/2m}}\r)^{\theta_{2}} e^{-c\l(\frac{2^{j} r_{i}}{t^{\frac{1}{2m}}}\r)^{\frac{2m}{2m-1}}}\frac{dt}{t^{3/2}} \l(\fint_{B_{i}}|b_{i}|^{q}d w\r)^{1/q}\\[4pt]
& \lesssim  \frac{1}{\alpha}\sum_{i}\sum_{j\geq 2}2^{j\theta_{1}}2^{jD}w( B_{i}) e^{-c(2^{j\frac{2m}{2m-1}})}  r_{i}^{-m} \l(\fint_{B_{i}}|b_{i}|^{q}dw\r)^{1/q}\lesssim \frac{1}{\alpha^{p}}\int_{\rz} |\nabla^{m} f|^{p}dw.\end{aligned}\end{equation*} 

Second, we handle $J_{1}.$ To the end, set $$\psi(z):=\frac{1}{\pi^{1/2}}\int_{ 1}^{\infty}ze^{-tz}\frac{dt}{t^{1/2}}\quad \mbox{and}\quad \beta_{k}:=\sum_{i: r_{i}=2^{k}} \frac{b_{i}}{r_{i}^{m}}.$$ Hence, $S_{i}=r_{i}^{-m} \psi(r_{i}^{2m} \L_{w})$ and $$\sum_{i} S_{i} b_{i}=\sum_{k\in \Z} \psi(2^{2mk} \L_{w})\beta_{k}.$$ An application of \eqref{eq: cd6000} (a higher-order extension of \cite[Proposition 5.14]{DMR} or weighted analogue of \cite[Lemma 21]{A1}; see Section 8 for the proof) yields 
$$\|\sum_{k\in \Z} \psi(2^{2mk} \L_{w})\beta_{k}\|_{L^{q}(w)}\lesssim \|\l(\sum_{k\in \Z} |\beta_{k}|^{2}\r)^{1/2}\|_{L^{q}(w)}.$$ From this,   together with \eqref{eq: cd40}-\eqref{eq: cd42}, it holds that 
\begin{equation*}\begin{aligned} J_{1}\lesssim \frac{1}{\alpha^{q}}\|\sum_{i} S_{i} b_{i}\|_{L^{q}(w)}^{q}&\lesssim \frac{1}{\alpha^{q}}\|\l(\sum_{k\in \Z} |\beta_{k}|^{2}\r)^{1/2}\|_{L^{q}(w)}^{q}\\[4pt] 
&\lesssim\frac{1}{\alpha^{q}} \int_{\rz} \sum_{i}|\frac{b_{i}}{r_{i}^{m}}|^{q}d w  \lesssim \frac{1}{\alpha^{p}}\int_{\rz} |\nabla^{m} f|^{p}dw. \\[4pt] 
\end{aligned}\end{equation*} By integrating the foregoing estimates, we thus reach \eqref{eq: cd84}.

Next, we prove \eqref{eq: cd82} via \eqref{eq: cd84}. To accomplish this, we need to generalize the interpolation technique developed in \cite{AM1} to accommodate higher-order scenarios. For any $p$ and $r$ such that $\max\{r_{w}, (p_{-})_{w, m, *}\}<r<p<2,$ the $L^{2}(w)$-Kato estimate \eqref{eq: cd1081}, along with \eqref{eq: cd84} implies that for all $ f\in \ss(\rz),$ \begin{equation}\label{eq: cd90} \|\L_{w}^{1/2}f\|_{L^{r, \infty}(w)} \lesssim \|\nabla^{m}f\|_{L^{r}(w)},\quad \|\L_{w}^{1/2}f\|_{L^{2}(w)} \lesssim \|\nabla^{m}f\|_{L^{2}(w)}.\end{equation} Additionally, for every $q>r_{w},$ we can adapt the proof from \cite[Lemma 6.7]{AM1} to demenstrate that  $$\Ez=\{(-\Delta)^{m/2}: f\in \ss(\rz), \;\supp\; \hat{f}\subset \rz\setminus \{0\}\}$$ is dense in $L^{q}(w),$ where $\hat{f}$ denotes the Fourier transform of $f.$ Furthermore, since $r>r_{w},$ we have $w\in A_{r}.$ Then, employing the properties of Riesz transforms, it follows that \begin{equation}\label{eq: cd1090}\|g\|_{L^{r}(w)}\approx \|\nabla^{m}(-\Delta)^{m/2}g\|_{L^{r}(w)}.\end{equation} Thus, for $g\in \Ez,$ \eqref{eq: cd1090} and $f:=(-\Delta)^{m/2}g$ imply $\L_{w}^{1/2}(-\Delta)^{m/2}g=\L_{w}^{1/2}f$ with $$\|\nabla^{m} f\|_{L^{r}(w)}\approx \|g\|_{L^{r}(w)}, \; \;\forall\; r>r_{w}.$$ Defining $T:=\L_{w}^{1/2}(-\Delta)^{m/2},$ we can rewrite \eqref{eq: cd90} as $$\|Tf\|_{L^{r, \infty}(w)} \lesssim \|f\|_{L^{r}(w)},\quad \|Tf\|_{L^{2}(w)} \lesssim \|f\|_{L^{2}(w)}, \quad \forall \; f\in \Ez.$$ Of course, by density arguments, we can extend these last two estimates to each $L^{q}(w),$ noting that their restrictions to the space of simple functions coincide. This allows us to apply Marcinkiewicz interpolation. For any $r<p<2,$ we conclude
$$\|Tf\|_{L^{p}} \lesssim \|f\|_{L^{p}(w)}, \quad \forall \; f\in \ss(\rz),$$ which is equivalent to $$\|\L_{w}^{1/2}f\|_{L^{p}(w)} \lesssim \|\nabla^{m}f\|_{L^{p}(w)},\quad \forall \; f\in \ss(\rz).$$ Using density once more, this gives \eqref{eq: cd82} on $L^{p}(w)$ for all $r<p<2.$ By the arbitrariness of $r$, \eqref{eq: cd82} holds for $p\in (\max\{r_{w}, (p_{-})_{w, m, *}\}, 2).$

To prove \eqref{eq: cd83} for $p$ satisfying $\max\{r_{w}, p_{-}\}<p<p_{+},$ we proceed as before by applying Theorem \ref{theorem: cd18}. Granting \eqref{eq: cd83}, then \eqref{eq: cd82} holds for $2\leq p<p_{+}$ by letting $v\equiv 1.$ 

Set $\tilde{p}_{-}=\max\{r_{w}, p_{-}\}<2$ and choose $p$ such that $\tilde{p}_{-}<p<p_{+}.$ Recall from Proposition \ref{proposition: cd3000} that there are $p_{0}, q_{0}$ such that $$\tilde{p}_{-}<p_{0}<\min\{p, 2\}\leq p<q_{0}<p_{+}\quad \mbox{and}\quad v\in A_{\frac{p}{p_{0}}}(w)\cap \mbox{RH}_{(\frac{q_{0}}{p})'}(w).$$ In order to use Theorem \ref{theorem: cd18}, we need to construct \eqref{eq: cd19}-\eqref{eq: cd20} for the operators $$\T=\L_{w}^{1/2}, \quad \S=\nabla^{m}, \quad \mbox{and}\quad \A_{r}=I-(I-e^{-r^{2m}\L_{w}})^{N}.$$ Since $p_{0}, q_{0}\in \J(\L_{w}),$ it holds that $\A_{r}\in \mho(L^{p_{0}}(w) \to L^{q_{0}}(w) )$ with estimates \eqref{eq: cd28}-\eqref{eq: cd30}. Combining this with the fact that $\A_{r}$ and $T$ commute, and using similar arguments as in \eqref{eq: cd34}, we can derive 
\eqref{eq: cd20} with $g(j):=C2^{j(\theta_{1} +\theta_{2})} e^{-c(2^{j\frac{2m}{2m-1}})}.$ 

At this stage, we are left to show \eqref{eq: cd19}. Given $f\in \ss(\rz),$ let $\phi(z)=z^{1/2}(1-e^{-r^{2m} z})^{N}.$ Clearly, $\phi(\L_{w})=\T(I-\A_{r}).$ Then, by Lemma \ref{lemma: cd67}, $$\phi(\L_{w})f=\phi(\L_{w})(f-\pi_{4B}^{m}(f))=\sum_{j\geq 1}\phi(\L_{w})h_{j},$$ with $\pi_{4B}^{m}(f) $ from \eqref{eq: cd72}, $h_{j}=(f-\pi_{4B}^{m}(f))\psi_{j},$ $\psi_{j}=1_{C_{j}(B)}$ for $j\geq 3,$ $\psi_{1}\in C_{0}^{\infty}(4B)$ ($1$ on $2B,$ $0\leq \psi \leq 1$ and $\|D^{\gamma} \psi_{1}\|_{\infty}\lesssim \frac{1}{r^{|\gamma|}}$ for any $|\gamma|\leq m$), and $\psi_{2}\in C_{0}^{\infty}(8B\setminus 2B)$ satisfying $\sum_{j\geq 1}\psi_{j}=1.$ To establish \eqref{eq: cd19}, we are required to handle each of these terms $$\l(\fint_{B}|\phi(\L_{w}) \psi_{j} |^{p_{0}}dw\r)^{1/p_{0}}\quad\mbox{for}\; j=1,2,....$$ 

Observe that $\phi(\L_{w}) \psi_{1}=(1-e^{-r^{2m} \L_{w}})^{N} \L_{w}^{1/2}\psi_{1}.$ Since $(1-e^{-r^{2m} z})^{N}\in \H^{\infty}(\Sigma_{\mu})$ with $\mu<\pi/2,$ utilizing Proposition \ref{proposition: cd15} (Remark \ref{remark: cd91}), we then get $$\l(\int_{\rz}|\phi(\L_{w}) \psi_{1} |^{p_{0}}dw\r)^{1/p_{0}}\lesssim \l(\int_{\rz}|\L_{w}^{1/2}\psi_{1} |^{p_{0}}dw\r)^{1/p_{0}}.$$ For $\tilde{p}_{-}<p_{0}<2,$ substituting $p=p_{0}$ into \eqref{eq: cd82} and applying the Leibniz rule along with \eqref{eq: cd44} allows us to deduce   
\begin{equation*}\begin{aligned} \l(\int_{\rz}|\L_{w}^{1/2}\psi_{1} |^{p_{0}}dw\r)^{1/p_{0}}&\lesssim \|\nabla^{m}\psi_{1}\|_{L^{p_{0}}(w)} \\[4pt] 
&\lesssim \sum_{\gamma\leq m}C_{m}^{\gamma} r^{-(m-|\gamma|)}\|D^{\gamma}(f-\pi_{4B}^{m}(f))\|_{L_{w}^{p_{0}}(4B)}\lesssim \|\nabla^{m}f\|_{L_{w}^{p_{0}}(4B)},\end{aligned}\end{equation*} which in turn implies $$\l(\fint_{B}|\phi(\L_{w}) \psi_{1} |^{p_{0}}dw\r)^{1/p_{0}}\lesssim \l(\fint_{4B}|\nabla^{m}f |^{p_{0}}dw\r)^{1/p_{0}}.$$
When $j\geq 3,$ we rewrite $\phi(\L_{w})\psi_{j}$ by the intergral representation from \eqref{eq: cd3004}-\eqref{eq: cd3005}
with $$|\eta_{\pm}(z, t)|\lesssim \frac{ r^{2mN}}{|z|^{N+3/2}},\quad z\in \Gamma_{\pm},$$ where $0<\V<\theta<\nu<\mu.$ Bear in mind that Corollary \ref{corollary: cd14} guarantees $e^{-z\L_{w}}\in \mho(L^{p_{0}}(w)\to L^{p_{0}}(w), \Sigma_{\frac{\pi}{2} -\theta})$ for $z\in \Gamma_{\pm}.$ Therefore, 
\begin{equation*}\begin{aligned}  &\l(\fint_{B}|\phi(\L_{w})\psi_{j} |^{p_{0}}dw\r)^{1/p_{0}}\\[4pt] 
&\quad\quad\lesssim \int_{\Gamma_{\pm}} \l(\fint_{B}|e^{-z\L_{w}}\psi_{j} |^{p_{0}}dw\r)^{1/p_{0}}|\eta_{\pm}(z)||dz|\\[4pt] 
&\quad\quad \lesssim 2^{j\theta_{1}}\l(\fint_{C_{j}(B)}|\psi_{j} |^{p_{0}}dw\r)^{1/p_{0}}\int_{\Gamma_{\pm}} \Upsilon\l( \frac{2^{j} r}{|z|^{\frac{1}{2m}}} \r)^{\theta_{2}} e^{-c( \frac{2^{j} r}{|z|^{\frac{1}{2m}}} )^{\frac{2m}{2m-1}}} \frac{ r^{2mN}}{|z|^{N+3/2}} |dz| \\[4pt] 
&\quad\quad \lesssim 2^{j(\theta_{1} -2mN-m)}r^{-m}\l(\fint_{2^{j+1}B} |f-\pi_{4B}^{m}(f)|^{p_{0}}dw\r)^{1/p_{0}}\quad (2mN>\theta_{2}+1)\\[4pt] 
&\quad\quad \lesssim 2^{j(\theta_{1} -2mN-m)}r^{-m}\l(\fint_{2^{j+1}B} |f-\pi_{2^{j+1}B}^{m}(f)|^{p_{0}}dw\r)^{1/p_{0}}\\[4pt]
\end{aligned}\end{equation*} 
\begin{equation*}\begin{aligned}&\quad+2^{j(\theta_{1} -2mN-m)}r^{-m}\l(\fint_{2^{j+1}B} |\pi_{2^{j+1}B}^{m}(f)-\pi_{4B}^{m}(f)|^{p_{0}}dw\r)^{1/p_{0}}.
\end{aligned}\end{equation*} 

Using key properties of the polynomial $\pi_{4B}^{m}(f)-$specifically \eqref{eq: cd72} and the argument in the proof of \cite[ Lemma 4.6]{Chua}-we can derive \begin{equation}\label{eq: cd1034}\begin{aligned} \l(\fint_{2^{j+1}B} |\pi_{2^{j+1}B}^{m}(f)-\pi_{4B}^{m}(f)|^{p_{0}}dw\r)^{1/p_{0}}&=\l(\fint_{2^{j+1}B} |\pi_{4B}^{m}(\pi_{2^{j+1}B}^{m}(f)-f)|^{p_{0}}dw\r)^{1/p_{0}}\\[4pt] 
&\lesssim 2^{jm} \fint_{4B}|\pi_{2^{j+1}B}^{m}(f)-f|dx\\[4pt] 
&\lesssim2^{j(m+n)} \fint_{2^{j+1}B}|\pi_{2^{j+1}B}^{m}(f)-f|dx\\[4pt] 
\end{aligned}\end{equation}
\begin{equation*}\begin{aligned}
&\lesssim2^{j(m+n)} \l(\fint_{2^{j+1}B}|\pi_{2^{j+1}B}^{m}(f)-f|^{p_{0}}dw\r)^{1/p_{0}}\\[4pt]&\lesssim2^{j(m+n)}r^{m}\l(\fint_{2^{j+1}B}|\nabla^{m}f|^{p_{0}}dw\r)^{1/p_{0}},
\end{aligned}\end{equation*} with the last step employing \eqref{eq: cd44}. Connecting the above two inequalities and using \eqref{eq: cd44} once more, we have $$\l(\fint_{B}|\phi(\L_{w})\psi_{j} |^{p_{0}}dw\r)^{1/p_{0}}\lesssim 2^{j(\theta_{1}+n-2mN)} \l(\fint_{2^{j+1}B}|\nabla^{m}f|^{p_{0}}dw\r)^{1/p_{0}}.$$ The case $j=2$ can be managed in the same manner, and we leave the details to the interested reader. 

Summarizing the previous estimates, we actually arrive at $$\l(\fint_{B}|\phi(\L_{w}) \psi |^{p_{0}}dw\r)^{1/p_{0}}\leq C \sum_{j\geq 1}2^{j(\theta_{1}+n-2mN)} \l(\fint_{2^{j+1}B}|\nabla^{m}f|^{p_{0}}dw\r)^{1/p_{0}}.$$ This inequality leads directly to \eqref{eq: cd19} under the condition that $2mN>\theta_{1}+\theta_{2}+n+1,$ thereby completing the entire proof for Proposition \ref{proposition: cd81}.

\hfill$\Box$

\section{Off-diagonal estimates for $t^{1/2}\nabla^{m} e^{-t\L_{w}}$ and key properties of $\K(\L_{w})$ } 
In this section, we provide the necessary preliminaries for proving the weighted $L^{p}$-boundedness of the Riesz transform $\nabla^{m}\L_{w}^{-1/2},$ which will be the main topic of the next section. We first connect the off-diagonal estimates for $e^{-t\L_{w}}$ and $t^{1/2}\nabla^{m} e^{-t\L_{w}}.$ Using this connection, we show the existence of an interval $\K(\L_{w})$ consisting of pairs $(p, q)$ for which $t^{1/2}\nabla^{m} e^{-t\L_{w}}\in \mho(L^{p}(w)\to L^{2}(w)),$ and establish its basic properties. Finally, we focus on showing that $2$ is an interior point of $\K(\L_{w}),$ which serves as a prerequisite for the arguments in the subsequent section. 

\subsection{The connection between off-diagonal estimates for $e^{-t\L_{w}}$ and $t^{1/2}\nabla^{m} e^{-t\L_{w}}$} For $p<2,$ the following lemma relates the off-diagonal estimates for $e^{-t\L_{w}}$ and $t^{1/2}\nabla^{m} e^{-t\L_{w}}.$ 

\begin{lemma}\label{lemma: cd105}\;  Given $1\leq p<2.$ The following are equivalent: 
 \begin{equation*}\begin{aligned} &(i)\quad e^{-t\L_{w}}\in \mho(L^{p}(w)\to L^{2}(w)).\\[4pt] 
&(ii)\quad t^{1/2}\nabla^{m} e^{-t\L_{w}}\in \mho(L^{p}(w)\to L^{2}(w)).\\[4pt] 
&(iii)\quad t\L_{w} e^{-t\L_{w}}\in \mho(L^{p}(w)\to L^{2}(w)).\end{aligned}\end{equation*}  \end{lemma}
{\it Proof.}\quad The proof proceeds similarly to that in \cite[Lemma 7.7]{DMR}, which originates from \cite[Lemma 5.3]{AM}. First, we show that $(i)$ implies $(ii)$. Indeed, Theorem \ref{theorem: cd1073} and Lemma \ref{lemma: h7} yield that $t^{1/2}\nabla^{m} e^{-t\L_{w}}\in \mho(L^{2}(w)\to L^{2}(w)).$ Consequently, $(ii)$ follows by applying Lemma \ref{lemma: h6} and the composition $(t^{1/2}\nabla^{m} e^{-t/2\L_{w}})\circ e^{-t/2\L_{w}}$. 

Second, we show that $(ii)$ implies $(iii)$. For any $\overrightarrow{f}:=(f_{\beta})_{|\beta|=m},$ define $$S_{t}\overrightarrow{f}:=t^{1/2} e^{-t\L_{w}}((-1)^{m}\sum_{|\alpha|=|\beta|=m}w^{-1}\partial^{\alpha}(a_{\alpha, \beta} f_{\beta}).$$ By duality in $L^{2}(w)$, the following holds:
 \begin{equation*}\begin{aligned} <S_{\lambda}\overrightarrow{f}, g>_{L^{2}(w)}&=<(-1)^{m}\sum_{|\alpha|=|\beta|=m}w^{-1}\partial^{\alpha}(a_{\alpha, \beta} f_{\beta}), t^{1/2} e^{-t\L_{w}^{*}}g>_{L^{2}(w)}\\[4pt] 
\end{aligned}\end{equation*}  \begin{equation*}\begin{aligned}&= \sum_{|\beta|=|\alpha|=m}< f_{\beta}, t^{1/2} w^{-1}\overline{a_{\alpha, \beta}}\partial^{\alpha}(e^{-t\L_{w}^{*}}g)>_{L^{2}(w)}.\end{aligned}\end{equation*} 
From this, together with Lemma \ref{lemma: h5}, $w^{-1}\overline{a_{\alpha, \beta}}\in L^{\infty}$ and $t^{1/2}\nabla^{m} e^{-t\L_{w}^{*}}\in \mho(L^{2}(w)\to L^{2}(w)),$ it follows that $S_{t}\in \mho(L^{2}(w)\to L^{2}(w)).$ Clearly, $S_{t}\circ (t^{1/2}\nabla^{m}e^{-t\L_{w}})=t\L_{w}e^{-2t\L_{w}},$ so the implication $(ii)\Longrightarrow (iii)$ follows from Lemma \ref{lemma: h6} and the semigroup property.

We now prove that $(iii)$ implies $(i),$ and so Lemma \ref{lemma: cd105} is concluded. In light of Definition \ref{definition: cd24}, we are required to construct \eqref{eq: cd25}-\eqref{eq: cd27} with $T=e^{-t\L_{w}}$. We first show \eqref{eq: cd25}. Fix a ball $B,$ and choose two functions $f, g$ in $ L^{2}(B, dw)$ such that $$\l(\fint_{B}|f|^{p}dw\r)^{1/p}=\l(\fint_{B}|g|^{2}dw\r)^{1/2}=1.$$ Then, by duality once more, it suffices to prove \begin{equation}\label{eq: cd106}|h(t)|\lesssim \Upsilon\l(\frac{r}{t^{1/2m}}\r)^{\theta}\end{equation} for some $\theta>0,$
where $$h(t):=\fint_{B}e^{-t\L_{w}}(1_{B} f)(x)\overline{g(x)} dw(x).$$ 
Since $e^{-tz}$ converges to $0$ on compact subsets of $\mbox{Re}\; z>0,$ we have $\lim_{t\to \infty}h(t)=0$ by the bounded holomorphic functional calculus of $\L_{w}$ on $L^{2}(w).$ Thus, \begin{equation}\label{eq: cd4000}h(t)=-\int_{t}^{\infty}sh'(s)\frac{ds}{s}.\end{equation} Let $ \widetilde{\Upsilon}(s)=\max\{s^{\alpha}, s^{\beta}\}.$ Applying Lemma \ref{lemma: cd103} with $t\L_{w} e^{-t\L_{w}}\in \mho(L^{p}(w)\to L^{2}(w))$, we see $$t|h'(t)|\lesssim \widetilde{\Upsilon}\l(\frac{r}{t^{1/2m}}\r).$$ From this and \eqref{eq: cd4000}, it follows that $$|h(t)|\lesssim \int_{t}^{\infty} \widetilde{\Upsilon}\l(\frac{r}{s^{1/2m}}\r)\frac{ds}{s}\lesssim \widetilde{\Upsilon}\l(\frac{r}{t^{1/2m}}\r)\lesssim \Upsilon\l(\frac{r}{t^{1/2m}}\r)^{\alpha+\beta}.$$ This gives \eqref{eq: cd106}, hence \eqref{eq: cd25}.  

The proof of \eqref{eq: cd26} is analogous to that described above. Fix $f\in L^{2}(C_{j}(B), dw)$ and $g\in L^{2}(B, dw)$ with $$\l(\fint_{C_{j}(B)}|f|^{p}dw\r)^{1/p}=\l(\fint_{B}|g|^{2}dw\r)^{1/2}=1,$$ and let $$h(t):=\fint_{B}e^{-t\L_{w}}(1_{C_{j}(B)} f)(x)\overline{g(x)} dw(x).$$ Since $e^{-t\L_{w}}\in \mho(L^{2}(w)\to L^{2}(w))$, we have $\lim_{t\to 0}h(t)=0$, and thus \eqref{eq: cd4000}. Using assumption $(iii)$ and Lemma \ref{lemma: z1}, we obtain \eqref{eq: cd26} through the following derivation: 
 \begin{equation*}\begin{aligned} |h(t)|&\lesssim 2^{j\theta_{1}}\int_{0}^{t}\Upsilon\l(\frac{2^{j}r}{s^{1/2m}}\r)^{\theta_{2}} e^{-c\l(\frac{2^{j}r}{s^{\frac{1}{2m}}}\r)^{\frac{2m}{2m-1}}}\frac{ds}{s}\\[4pt] 
&\lesssim 2^{j\theta_{1}}\int_{\frac{2^{j}r}{t^{1/2m}}}^{\infty} \Upsilon^{\theta_{2}}(s)e^{-cs^{\frac{2m}{2m-1}}}\frac{ds}{s}\lesssim 2^{j\theta_{1}} \Upsilon\l(\frac{2^{j}r}{t^{1/2m}}\r)^{\theta_{2}} e^{-c\l(\frac{2^{j}r}{t^{\frac{1}{2m}}}\r)^{\frac{2m}{2m-1}}}.\end{aligned}\end{equation*} Exchanging the roles of $B$ and $C_{j}(B),$ an analogous argument leads to \eqref{eq: cd27}. The details are omitted. 

\hfill$\Box$

\subsection{Basic properties of the interval $\K(\L_{w})$} We introduce the set $\widetilde{\K}(\L_{w}):=\{p\in [1, \infty]: \sup_{t>0}\|t^{1/2}\nabla^{m} e^{-t\L_{w}}\|_{L^{p}(w)\to L^{p}(w)}<\infty\}.$ Like the set $\widetilde{\J}(\L_{w}),$ $2\in \widetilde{\K}(\L_{w})$ thanks to Theorem \ref{theorem: cd1073} and Lemma \ref{lemma: h7}, and $\widetilde{\K}(\L_{w})$ will be an interval if it contains more than one point. As mentioned earlier, we denote by $\K(\L_{w})$ the set of all pairs $(p, q)$ such that $t^{1/2}\nabla^{m} e^{-t\L_{w}}\in \mho(L^{p}(w)\to L^{q}(w)).$ According to \cite[Remark 7.2]{DMR}, a complete characterization of $\K(\L_{w})$ is not possible, in contrast to the unweighted setting \cite{A}, due to the absence of a proof that $p<q<2$ and $t^{1/2}\nabla^{m} e^{-t\L_{w}}\in \mho(L^{p}(w)\to L^{q}(w))$ imply $p, q\in \K(\L_{w})$ (with $w\in A_{2}$ and $p, q$ possibly close to $1$).

\begin{proposition}\label{proposition: cd100}\;There exists an interval $\K(\L_{w})$ such that if $p, q\in \K(\L_{w}),$ $p\leq q,$ then $t^{1/2}\nabla^{m} e^{-t\L_{w}}\in \mho(L^{p}(w)\to L^{q}(w)).$ Moreover, $\K(\L_{w})$ has the following properties: 
\begin{equation*}\begin{aligned} 
&(i)\quad \K(\L_{w})\subset \widetilde{\K}(\L_{w}). \\[4pt] 
&(ii)\quad\mbox{If}\; q_{-}(\L_{w})\;\mbox{and}\; q_{+}(\L_{w})\;\mbox{denote the left and right endpoints of}\;\K(\L_{w}),\;\mbox{then}\; q_{-}(\L_{w})=p_{-},\\[4pt] 
&\quad\quad  2\leq q_{+}(\L_{w})\leq (q_{+}(\L_{w}))^{*, m}_{w}\leq p_{+},\;2\in \K(\L_{w})\;\mbox{and}\; \K(\L_{w})\subset \J(\L_{w}). \\[4pt]
&(iii)\quad \mbox{If }\; q\geq 2,\; p\leq q,\;\mbox{and}\; t^{1/2}\nabla^{m} e^{-t\L_{w}}\in \mho(L^{p}(w)\to L^{q}(w)),\;\mbox{then}\; p, q\in \K(\L_{w}).\\[4pt]
&(v)\quad \sup  \widetilde{\K}(\L_{w})=q_{+}(\L_{w}).\end{aligned}\end{equation*} 
\end{proposition}
{\it Proof.}\quad Let $\K(\L_{w}):=\K_{-}(\L_{w})\cup \K_{+}(\L_{w}),$ where $$\K_{-}(\L_{w}):=\{p\in [1, 2]: t^{1/2}\nabla^{m} e^{-t\L_{w}}\in \mho(L^{p}(w)\to L^{2}(w))\}$$ and $$\K_{+}(\L_{w}):=\{p\in [2, \infty]: t^{1/2}\nabla^{m} e^{-t\L_{w}}\in \mho(L^{2}(w)\to L^{p}(w))\}.$$ Clearly, by Proposition \ref{proposition: cd12}, Lemma \ref{lemma: cd105} and Lemma \ref{lemma: h4}, $\K_{-}(\L_{w})$ is an interval, and so is $\K(\L_{w}).$ For any $p, q\in \K(\L_{w})$ with $p<q,$ by applying Lemma \ref{lemma: h4}, Lemma \ref{lemma: cd105} and Lemma \ref{lemma: h6} in place of \cite[Lemma 2.28]{DMR}, \cite[Lemma 2.30]{DMR} and \cite[Lemma 7.7]{DMR} respectively, and following a similar argument to that in \cite[Proposition 7.1]{DMR}, we can show that $t^{1/2}\nabla^{m} e^{-t\L_{w}}\in \mho(L^{p}(w)\to L^{q}(w)).$

Property $(i)$ follows immediately from Lemma \ref{lemma: h4}.

We now prove property $(ii).$ When $p<2,$ it follows from Lemma \ref{lemma: cd105} and Proposition \ref{proposition: cd12} that $p\in \J(\L_{w})$ if and only if $p\in \K_{-}(\L_{w}).$ Hence $\J(\L_{w})\cap [1, 2]=\K_{-}(\L_{w}),$ which implies $q_{-}(\L_{w})=p_{-}(\L_{w}).$ 

Note that, if $q_{+}(\L_{w})=2,$ $(q_{+}(\L_{w}))^{*, m}_{w}=2^{*, m}_{w}\leq p_{+}(\L_{w})$ by Proposition \ref{proposition: cd12}. For $q_{+}(\L_{w})>2,$ choose $p, q$ such that $2<p<q_{+}(\L_{w})$ and $p<q<p^{*, m}_{w}.$ As $2, p\in \K_{+}(\L_{w}),$ we have $ e^{-t\L_{w}}\in \mho(L^{2}(w)\to L^{2}(w))$ and $t^{1/2}\nabla^{m} e^{-t\L_{w}}\in \mho(L^{2}(w)\to L^{p}(w)).$ Since $A_{2}\subset A_{p},$ by adapting the approach used to handle $J_{1}, J_{2}, J_{3}$ in \eqref{eq: cd10}, we obtain \begin{equation*}\begin{aligned} \l(\fint_{B}e^{-t\L_{w}}(f1_{B})|^{q}dw\r)^{\frac{1}{q}}& \lesssim \l(\fint_{B}|e^{-t\L_{w}}(f1_{B})|^{2}dw\r)^{\frac{1}{2}}+ r^{m}\l(\fint_{B}\nabla^{m}e^{-t\L_{w}}(f1_{B})|^{p}dw\r)^{\frac{1}{p}}\\[4pt] &\lesssim \Upsilon\l(\frac{r}{t^{1/2m}}\r)^{m+\theta_{2}} \l(\fint_{B}|f|^{2}dw\r)^{\frac{1}{2}}.\end{aligned}\end{equation*} This is \eqref{eq: cd25} in Definition \ref{definition: cd24}. Similarly, \eqref{eq: cd26}-\eqref{eq: cd27} can be proved. Thus, $ e^{-t\L_{w}}\in \mho(L^{2}(w)\to L^{q}(w)),$ so $(q_{+}(\L_{w}))^{*, m}_{w}\leq p_{+}(\L_{w})$ by letting $p\nearrow q_{+}(\L_{w})$ and $q\nearrow p^{*, m}_{w}.$ 

Of course, $q_{+}(\L_{w})\leq p_{+}(\L_{w}).$ If $q_{+}(\L_{w})<\infty,$ then $q_{+}(\L_{w})<(q_{+}(\L_{w}))^{*, m}_{w} \leq p_{+}(\L_{w})$ and so $\K_{+}(\L_{w})\subset \J(\L_{w}).$ Otherwise, $p_{+}=\infty,$ which yields $\K_{+}(\L_{w})\subset \J(\L_{w})$ trivially. This completes the proof for property $(ii).$

Property $(iii)$ and $(v)$ follow similarly to page 642 of \cite{DMR} (or \cite[Proposition 5.6]{AM}), using Lemma \ref{lemma: h3} and Lemma \ref{lemma: h4} instead of \cite[Lemma 2.28]{DMR} and \cite[Lemma 2.27]{DMR}. Details are left to the reader.

\hfill$\Box$

\begin{corollary}\label{corollary: cd1010}\;  Let $q_{-}(\L_{w})<p\leq q <q_{+}(\L_{w}).$ If  $v\in A_{\frac{p}{q_{-}(\L_{w})}}(w)\cap \mbox{RH}_{(\frac{q_{+}(\L_{w})}{q})'}(w), $ then $t^{1/2}\nabla^{m}e^{-t\L_{w}} \in \mho( L^{p}(vdw)\to L^{q}(vdw))$ and $z^{1/2}\nabla^{m} e^{-z\L_{w}} \in \mho( L^{p}(dw)\to L^{q}(dw))$ for all $z\in \Sigma_{\mu}$ with $0<\mu<\frac{\pi}{2}-\V.$
\end{corollary}{\it Proof.}\quad See the proof of Corollary \ref{corollary: cd13}-\ref{corollary: cd14}.\hfill$\Box$

\subsection{A key interior point of $\K(\L_{w})$ } We now prove that $2\in \mbox{Int}\;\K(\L_{w}).$ For this purpose, we first recall a reverse H\"{o}lder inequality with sharp constants for solutions of $\L_{w},$ 
established in \eqref{eq: cd1012} and \eqref{eq: cd1016}. More precisely, for a fixed ball $B_{0},$ if $u\in H^{m}(4B_{0}, w)$ is a solution of $\L_{w} u=0$ in $4B_{0},$ then 
\begin{equation}\label{eq: cd1018}\l(\fint_{B}|\nabla^{m} u|^{2}dw\r)^{1/2} \leq  \frac{C_{1}}{r^{m}}\l(\fint_{2B } |u-P_{2B}(u)|^{2}dw\r)^{1/2}\end{equation} holds for any ball $B$ such that $3B\subset 4B_{0},$ where $C_{1}:=C[w]_{A_{2}}^{m/2}$ and $P_{2B}$ is as defined in Corollary \ref{corollary: cd1}. Since $r_{w}<2,$  we can always find a $q$ such that $$\max\{r_{w}, \frac{2n\tau_{w}}{n\tau_{w}+2m}\}<q<2\leq n.$$ With this choice of $q,$ we have $2<q^{*, m}_{w}.$ Consequently, by Corollary \ref{corollary: cd1}, there exists $q^{\star}\in (r_{w}, 2)$ such that \begin{equation}\label{eq: cd1020}\frac{1}{r^{m}}\l(\fint_{2B } |u-P_{2B}(u)|^{2}dw\r)^{1/2}\leq C_{2} \l(\fint_{2 B}|\nabla^{m}u|^{q}dw\r)^{\frac{1}{q}},\end{equation} where $C_{2}:=C [w]^{\frac{1}{q}}_{A_{q^{\star}}}.$
Combining \eqref{eq: cd1018} and \eqref{eq: cd1020}, we get $$\l(\fint_{B}|\nabla^{m} u|^{2}dw\r)^{1/2} \leq C_{1}C_{2} \l(\fint_{2 B}|\nabla^{m}u|^{q}dw\r)^{\frac{1}{q}}.$$ From this, \cite[Theorem 3.22]{BB} applies, so there exists a $p_{0}>2$ such that for every admissible ball $B,$ 
\begin{equation}\label{eq: cd1021}\l(\fint_{B}|\nabla^{m} u|^{p_{0}}dw\r)^{1/p_{0}} \leq C_{3} \l(\fint_{2 B}|\nabla^{m}u|^{2}dw\r)^{\frac{1}{2}},\end{equation} where $C_{3}:=8^{1/q}C_{1}C_{2} (2^{D}[w]_{A_{2}})^{31/q}$ and 
\begin{equation}\label{eq: cd1022}p_{0}:=2+\frac{2-q}{2^{\frac{4}{q} +1} C_{1}^{2} C_{2}^{2}(2^{D}[w]_{A_{2}})^{\frac{6}{q} +17} }.\end{equation} 

To proceed, as shown in \cite[Section 8]{DMR}, we need to introduce the Riesz transform $\nabla^{m} \L_{w}^{1/2}$ associated with the higher-order weighted elliptic operator $\L_{w}.$ In fact, $\nabla^{m} \L_{w}^{1/2}$ can be defined by \begin{equation}\label{eq: cd1027}\nabla^{m}\L_{w}^{-1/2}=\frac{1}{\sqrt{\pi}}\int_{0}^{\infty} t^{1/2}\nabla^{m} e^{-t\L_{w}}\frac{d t}{t}.\end{equation} To verify this definition, we must show that this integral is well-posed, meaning it converges at both $0$ and $\infty.$ For this aim, for any $\ez>0,$ we introduce \begin{equation}\label{eq: cd1032}S_{\ez}:=S_{\ez}(\L_{w}):=\frac{1}{\sqrt{\pi}}\int_{\ez}^{1/\ez}t^{1/2}e^{-t\L_{w}}\frac{d t}{t}.\end{equation} It is easy to see that, for each $0<\ez<1,$ the function $S_{\ez}(z):=\frac{1}{\sqrt{\pi}}\int_{\ez}^{1/\ez}t^{1/2}e^{-tz}\frac{d t}{t}$ is holomorphic and uniformly bounded on the right half-plane. From the results in Section 3.1, it follows that $$\|S_{\ez}(\L_{w}) f\|_{L^{2}(w)}\leq C \|S_{\ez}(z)\|_{\infty} \|f\|_{L^{2}(w)} \leq C \|f\|_{L^{2}(w)}$$ with the constant $C$ independent of both $\ez$ and $f.$ Note that, given $f\in C_{0}^{\infty}(\rz),$ $S_{\ez} f\in \D(\L_{w}) \subset \D(\L_{w}^{1/2}),$ so \begin{equation}\label{eq: cd1033}\|\nabla^{m}S_{\ez} f \|_{L^{2}(w)}\lesssim \|\L_{w}^{1/2}S_{\ez} f \|_{L^{2}(w)} =\|\phi_{\ez}(\L_{w}) f\|_{L^{2}(w)}, \end{equation} where $$\phi_{\ez}(z):=\frac{1}{\sqrt{\pi}}\int_{\ez}^{1/\ez}t^{1/2}z^{1/2}e^{-tz}\frac{d t}{t}.$$ Then, we can deduce that $\L_{w}^{1/2}S_{\ez}f\to f$ strongly in $L^{2}(w),$ as $\{\phi_{\ez}\}_{0<\ez<1}$ is uniformly bounded and converges uniformly to 1 on compact subsets of the sector $\Sigma_{\mu}$ with $0<\mu<\frac{\pi}{2}.$ Combining this and \eqref{eq: cd1033}, we see that $\{\nabla^{m} S_{\ez} f\}$ is a Cauchy sequence in $L^{2}(w)$. Hence, we can define $$\nabla^{m} \L_{w}^{-1/2}f=\lim_{\ez\to 0}\nabla^{m} S_{\ez} f$$ with the limit interpreted in $L^{2}(w),$ thereby proving \eqref{eq: cd1027}. In what follows, when considering $L^{2}(w)$ 
estimates for $\nabla^{m} \L_{w}^{1/2}$, we actually establish estimates for $\nabla^{m}S_{\ez}$ with constants independent of $\ez.$ These arguments are implicit unless details need to be emphasized.

Having established the above, we are able to define the Hodge projection operator by $$\H:= \nabla^{m} \L_{w}^{-1/2}(\nabla^{m} (\L_{w}^{*})^{-1/2})^{*},$$ where adjoints are taken with respect to the $L^{2}(w)$ inner product. As the Riesz transform is bounded on $L^{2}(w)$ by the $L^{2}(w)-$Kato estimate \eqref{eq: cd1081}, the Hodge projection $\H$ is also bounded. Moreover, \begin{equation}\label{eq: cd4001}\H =(-1)^{m} \nabla^{m}\L_{w}^{-1}(w^{-1} \mbox{div}_{m} (w \cdot))\end{equation} since we have $(\nabla^{m} (\L_{w}^{*})^{-1/2})^{*}=(-1)^{m} \L_{w}^{-1/2}(w^{-1} \mbox{div}_{m} (w \cdot))$ by duality. 

Fix a function $\overrightarrow{f}=(f_{\beta})_{|\beta|=m}\in L^{2}(w)\cap L^{p_{0}}(w)$ with $\supp\; \overrightarrow{f} \subset \rz\setminus 4B_{0},$ and let $u:=L_{w}^{-1}(w^{-1} \mbox{div}_{m} (w \overrightarrow{f})).$ Then, we can prove $\nabla^{m} u\in L^{2}(w)$ using duality arguments and the $L^{2}(w)$ boundedness of the Riesz transform. As a consequence of \eqref{eq: cd4001}, $$\H \overrightarrow{f}=(-1)^{m}  \nabla^{m} u$$ holds in the sense of distributions. Clearly, $L_{w}u=0$ on $4B_{0}$ since $\supp\; \overrightarrow{f} \subset \rz\setminus 4B_{0}.$ Indeed, exploiting a standard Lax-Milgram argument (guaranteed by \eqref{eq: z1}-\eqref{eq: z2}) along with the generalized Poincar\'{e}-Sobolev inequality in Theorem \ref{theorem: cd0}, we can derive $u\in H^{m}(4B_{0}, w).$ 
This allows us to use \eqref{eq: cd1021} to deduce that for any ball $B$ such that $3B\subset 4B_{0},$ $$\l(\fint_{B}|\H \overrightarrow{f}|^{p_{0}}dw\r)^{1/p_{0}} =\l(\fint_{B}|\nabla^{m} u|^{p_{0}}dw\r)^{1/p_{0}} \leq C_{3} \l(\fint_{2 B}|\H \overrightarrow{f}|^{2}dw\r)^{\frac{1}{2}}.$$ Thus, invoking \cite[Theorem 3.14]{AM2} in the context of homogeneous spaces (see \cite[Section 5]{AM2}), we immediately get that $\H : L^{q}(w)\to L^{q}(w)$ for all $q,$ $2\leq q<p_{0}.$ 
Equivalently, $\H^{*} : L^{q'}(w)\to L^{q'}(w)$ for all $q',$ $p_{0}'<q'\leq 2.$




Utilizing the $L^{q'}(w)$ boundedness of $\H^{*}$ just established, we can show that the Reisz transform $\nabla^{m} \L_{w}^{1/2}$ is bounded on $L^{q}(w)$ for all $q$ satisfying $$2<q<\min\{p_{+}(\L_{w}), r_{w}', p_{0}\}:=q_{w}.$$ Equivalently, 
$(\nabla^{m} \L_{w}^{1/2})^{*}$ is bounded on $L^{q'}(w)$ for $q'$ such that $$(p_{-}(\L_{w}^{*}))_{w, m, *}\leq \max\{p_{-}(\L_{w}^{*}), r_{w}, p_{0}'\}<q'\leq 2,$$ as $p_{-}(\L_{w}^{*})'=p_{+}(\L_{w})$ by Lemma \ref{lemma: h5}. The proof is straightforward. Note that \eqref{eq: cd4001} implies $\H^{*}\overrightarrow{f}=(-1)^{m} (\L^{*}_{w})^{-1}(w^{-1} \mbox{div}_{m} (w \cdot))$ by duality. Therefore, \begin{equation*}\begin{aligned}\|(\nabla^{m} \L_{w}^{-1/2})^{*}\overrightarrow{f}\|_{L^{q'}(w)} & =\|(\L^{*}_{w})^{-1/2}(w^{-1} \mbox{div}_{m} (w \overrightarrow{f}))\|_{L^{q'}(w)} \\[4pt] 
&\lesssim \|\nabla^{m}(\L^{*}_{w})^{-1}(w^{-1} \mbox{div}_{m} (w \overrightarrow{f}))\|_{L^{q'}(w)} \approx \| \H^{*}\overrightarrow{f}\|_{L^{q'}(w)}\lesssim \| \overrightarrow{f}\|_{L^{q'}(w)},\end{aligned}\end{equation*} where we used Proposition \ref{proposition: cd81} in the derivation.



We claim that $t^{1/2}\nabla^{m} e^{-t\L_{w}}: L^{q}(w)\to L^{q}(w)$ for all $q\in (2, q_{w}).$ If this claim holds, then by Proposition \ref{proposition: cd100}, $$q_{+}(\L_{w})=\sup \widetilde{\K}(\L_{w})\geq q_{w}>2,$$ which implies that $2$ is a interior point of $\K(\L_{w})$ as desired. To prove the claim, let $\phi_{t}(z):=(tz)^{1/2}e^{-tz}$ for any $t>0.$ It is easy to see that $\phi_{t}(z)$ is holomorphic and uniformly bounded on compact subsets of the right half-plane, with $\|\phi_{t}\|_{\infty, \Sigma_{\mu}}\leq C_{\mu}$ for any $\V<\mu <\pi/2.$ Hence, for any $q\in (2, q_{w}),$ applying Proposition \ref{proposition: cd15} and the previous estimates for the Reisz transform yields \begin{equation*}\begin{aligned}\|t^{1/2}\nabla^{m} e^{-t\L_{w}}f\|_{L^{q}(w)} & =\|\nabla^{m}\L_{w}^{-1/2} t^{1/2}\L_{w}^{1/2}e^{-t\L_{w}}f\|_{L^{q}(w)} \\[4pt] 
&\lesssim \|t^{1/2}\L_{w}^{1/2}e^{-t\L_{w}}f\|_{L^{q}(w)} \approx \|\phi_{t}(\L_{w})e^{-t\L_{w}}f\|_{L^{q}(w)}\lesssim \|f\|_{L^{q}(w)},\end{aligned}\end{equation*} where the implicit constants are independent of $t.$ This proves the claim.

\section{Estimates for higher-order Reisz transform in weighted spaces}
This section is devoted to proving the weighted $L^{p}$-estimates for the Reisz transform $\nabla^{m} \L_{w}^{1/2},$ which represents the reverse direction of the inequalities \eqref{eq: cd2001}. We will follow the approach in \cite[Proposition 9.1]{DMR} (stemming from \cite{AM3}), whose novelty lies in avoiding the use of (generalized) Poincar\'{e} inequalities, thereby accommodating the case where $p$ is close to 1 and the weight $w$ is in $A_2$ only.

\begin{proposition}\label{proposition: cd1023}\; Let $q_{-}(\L_{w})<p<q_{+}(\L_{w}).$ Then 
\begin{equation}\label{eq: cd1024}\|\nabla^{m} \L_{w}^{-1/2}f\|_{L^{p}(w)}\lesssim \|f\|_{L^{p}(w)},\end{equation} moreover, if $v\in A_{\frac{p}{q_{-}(\L_{w})}}(w)\cap \mbox{RH}_{(\frac{q_{+}(\L_{w})}{p})'}(w),$ \begin{equation}\label{eq: cd1025}\|\nabla^{m} \L_{w}^{-1/2}f\|_{L^{p}(vdw)}\lesssim \|f\|_{L^{p}(vdw)},\end{equation} with the implicit constants independent of $f.$
\end{proposition}
{\it Proof.}\quad For brevity, we set $q_{-}:=q_{-}(\L_{w})$ and $q_{+}:=q_{+}(\L_{w})$ throughout the following. We begin by proving \eqref{eq: cd1024} in the interval $(2,\; q_{+})$ (ensured by Section 6.3). To this end, we proceed by invoking Theorem \ref{theorem: cd18}, as in Proposition \ref{proposition: cd15}.

 Fix $p$ with $2<p<q_{+},$ and let $p_{0}=2$ and $q_{0}$ such that $2<p<q_{0}<q_{+}.$ We will show that the two conditions \eqref{eq: cd19}- \eqref{eq: cd20} of Theorem \ref{theorem: cd18} are satisfied with the choices: $(p_{0}, q_{0}),$ $\T:=\nabla^{m} \L_{w}^{-1/2},$ $\S=I$ and $\D:=L_{c}^{\infty}.$ 
As previously, we define $\A_{r}:=I-(I-e^{-r^{2m} \L_{w}})^{N}$, with $N$ to be specified later. For $f \in \D$, consider its decomposition $$f=\sum_{j\geq 1}f1_{C_{j}(B)}:=\sum_{j\geq 1}f_{j}.$$ 
Then, it is easy to see $$\l(\fint_{B}|\T (I-\A_{r})f|^{p_{0}} dw\r)^{1/p_{0}}\leq \sum_{j\geq 1}\l(\fint_{B}|\nabla^{m} \L_{w}^{-1/2}(I-e^{-r^{2m} \L_{w}})^{N}f_{j}|^{p_{0}} dw\r)^{1/p_{0}}.$$ 
First, it follows from \eqref{eq: cd1081}, Theorem \ref{theorem: cd1073} and Lemma \ref{lemma: h7} that $$\l(\fint_{B}|\nabla^{m} \L_{w}^{-1/2}(I-e^{-r^{2m} \L_{w}})^{N}f_{1}|^{p_{0}} dw\r)^{1/p_{0}}\lesssim \l(\fint_{4B}|f|^{p_{0}} dw\r)^{1/^{p_{0}} }.$$
Second, for any $j\geq 2$ and $h\in L^{2}(w),$ an application of \eqref{eq: cd1027} yields \begin{equation}\label{eq: cd1028}\nabla^{m} \L_{w}^{-1/2}(I-e^{-r^{2m} \L_{w}})^{N}h=c\int_{0}^{\infty} t^{1/2} \nabla^{m} \phi(\L_{w}, t)h \frac{dt}{t},\end{equation} where $\phi(z, t):=e^{-tz} (I-e^{-r^{2m} z})^{N}\in \H_{0}^{\infty}(\Sigma_{\mu}) (\V<\mu<\frac{\pi}{2}).$ Moreover, $\phi(\L_{w}, t)$ admits a representation given by \eqref{eq: cd3004}-\eqref{eq: cd3005}, with $\eta_{\pm}(z, t)$ satisfying
$$|\eta_{\pm}(z, t)|\lesssim \frac{r^{2mN}}{(|z|+t)^{N+1}}, \quad \mbox{for all}\; z\in \Gamma_{\pm}, \;t>0.$$ Then, by this and Theorem \ref{theorem: cd1073} (or Corollary \ref{corollary: cd1010}), we can deduce that 
\begin{equation}\label{eq: cd1029}\begin{aligned} &\l(\fint_{B}\bigg |\int_{\Gamma_{\pm}}t^{1/2}\nabla^{m}e^{-z\L_{w}}f_{j}\eta_{\pm}(z, t)dz\bigg |^{p_{0}}dw\r)^{\frac{1}{p_{0}}}\\[4pt]
&\quad\quad\quad\quad\lesssim \int_{\Gamma_{\pm}}\l(\fint_{B}|  z^{1/2} \nabla^{m} e^{-z\L_{w}}f_{j} |^{p_{0}}dw\r)^{\frac{1}{p_{0}}}\frac{t^{1/2}}{|z|^{1/2}}|\eta_{\pm}(z, t)| |dz|\\[4pt]&\quad\quad\quad\quad\lesssim 2^{j\theta_{1}}\l(\fint_{C_{j}(B)}|f|^{p_{0}}dw\r)^{1/p_{0}}\int_{\Gamma_{\pm}}\Upsilon\l(\frac{2^{j} r}{|z|^{1/2m}}\r)^{\theta_{2}}e^{-c\l(\frac{2^{j}r}{|z|^{\frac{1}{2m}}}\r)^{\frac{2m}{2m-1}}}\frac{t^{1/2}}{|z|^{1/2}}|\eta_{\pm}(z, t)| |dz|\\[4pt]
&\quad\quad\quad\quad\lesssim \l(\fint_{C_{j}(B)}|f|^{p_{0}}dw\r)^{1/p_{0}} 2^{j\theta_{1}}\int_{0}^{\infty}\frac{r^{2mN}}{(s+t)^{N+1}}\Upsilon\l(\frac{2^{j} r}{s^{1/2m}}\r)^{\theta_{2}}e^{-c\l(\frac{2^{j}r}{s^{\frac{1}{2m}}}\r)^{\frac{2m}{2m-1}}}\frac{t^{1/2}}{s^{1/2}} ds.\end{aligned}
\end{equation} Combining \eqref{eq: cd1029} and \eqref{eq: cd1028}, we achieve
\begin{equation}\label{eq: cd1030}\begin{aligned} &\l(\fint_{B}|\nabla^{m} \L_{w}^{-1/2}(I-e^{-r^{2m} \L_{w}})^{N}f_{j}|^{p_{0}} dw\r)^{1/p_{0}}\\[4pt]
&\lesssim \l(\fint_{C_{j}(B)}|f|^{p_{0}}dw\r)^{1/p_{0}} 2^{j\theta_{1}} \int_{0}^{\infty}\int_{0}^{\infty}\frac{r^{2mN}}{(s+t)^{N+1}}\Upsilon\l(\frac{2^{j} r}{s^{1/2m}}\r)^{\theta_{2}}e^{-c\l(\frac{2^{j}r}{s^{\frac{1}{2m}}}\r)^{\frac{2m}{2m-1}}}\frac{t^{1/2}}{s^{1/2}} \frac{ds dt}{t}\\[4pt]
&\lesssim 2^{-2mNj} 2^{j\theta_{1}}\l(\fint_{C_{j}(B)}|f|^{p_{0}}dw\r)^{1/p_{0}},\\[4pt]
\end{aligned}
\end{equation} provided $2mN>\theta_{2}+1.$ Summing over all $j\geq 1$ and using \eqref{eq: cd1030} we get \eqref{eq: cd19} with $g(j):=C2^{j(\theta_{1}-2mN)}$ if we further impose the condition $2mN>\theta_{1}+\theta_{2}+1.$

The proof of \eqref{eq: cd20} relies on the following key estimate: for every $f\in H^{m}(w)$ and $1\leq k\leq N,$ 
\begin{equation}\label{eq: cd1031}\l(\fint_{B}|\nabla^{m} e^{-kr^{2m} \L_{w}}f|^{q_{0}} dw\r)^{1/q_{0}}\lesssim \sum_{j\geq 1}g(j)\l(\fint_{2^{j+1}B}|\nabla^{m}f|^{p_{0}} dw\r)^{1/^{p_{0}} },\end{equation} where $g(j):=2^{j(\theta_{1}+\theta_{2}+m+n)} e^{-c(2^{j})^{\frac{2m}{2m-1}}}.$ To verify this estimate, fix $1\leq k\leq N$ and $f\in H^{m}(w),$ and define $h:=f-\pi_{4B}^{m}(f).$
Then, from Lemma \ref{lemma: cd67}, it follows that $$\nabla^{m} e^{-kr^{2m} \L_{w}}f=\nabla^{m} e^{-kr^{2m} \L_{w}}h:=\sum_{j\geq 1}\nabla^{m} e^{-kr^{2m} \L_{w}}h_{j},$$ where $h_{j}:=h1_{C_{j}(B)}.$ Therefore, $$\l(\fint_{B}|\nabla^{m} e^{-kr^{2m} \L_{w}}f|^{q_{0}} dw\r)^{1/q_{0}}\lesssim \sum_{j\geq 1}\l(\fint_{B}|\nabla^{m} e^{-kr^{2m} \L_{w}}h_{j}|^{q_{0}} dw\r)^{1/^{q_{0}} }.$$ As $2<q_{0}<q_{+},$ Proposition \ref{proposition: cd100} implies that $t^{1/2}\nabla^{m} e^{-t\L_{w}}\in \mho(L^{p_{0}}(w)\to L^{q}(w)).$ By this, together with \eqref{eq: cd44}, it holds that for each $j\geq 1,$
\begin{equation*}\begin{aligned}\l(\fint_{B}|\nabla^{m} e^{-kr^{2m} \L_{w}}  (h1_{C_{j}(B)})|^{q}dw\r)^{\frac{1}{q}}&\lesssim \frac{2^{j(\theta_{1}+\theta_{2})} e^{-c(2^{j\frac{2m}{2m-1}})} }{r^{m}}\l(\fint_{2^{j+1}B}|h|^{p_{0}}dw\r)^{\frac{1}{p_{0}}}  \\[4pt] 
&\lesssim \frac{2^{j(\theta_{1}+\theta_{2}+m+n) } e^{-c(2^{j\frac{2m}{2m-1}})} }{r^{m}} \l(\fint_{2^{j+1}B}|\pi_{2^{j+1}B}^{m}(f)-f|^{p_{0}}dw\r)^{\frac{1}{p_{0}}} \\[4pt]
&\lesssim 2^{j(\theta_{1}+\theta_{2}+m+n)}  e^{-c(2^{j\frac{2m}{2m-1}}) } \l(\fint_{2^{j+1}B}|\nabla^{m}f|^{p_{0}}dw\r)^{\frac{1}{p_{0}}},\end{aligned}\end{equation*} 
where in the last second step we have employed the same reasoning as in \eqref{eq: cd1034}. This gives us \eqref{eq: cd1031}. 

Note that for any fixed $\ez>0,$ the function $S_{\ez}f$ defined in \eqref{eq: cd1032} belongs to $H^{m}(w)$ because of the $L^{2}(w)$-boundedness of $e^{-t\L_{w}}$ and $t^{1/2}\nabla^{m}e^{-t\L_{w}}.$ By this, the commutativity of $\A_{r}$ and $S_{\ez}$, along with \eqref{eq: cd1031}, implies that $$\l(\fint_{B}|\nabla^{m} S_{\ez }\A_{r}f|^{q_{0}} dw\r)^{1/q_{0}}\lesssim \sum_{j\geq 1}g(j)\l(\fint_{2^{j+1}B}|\nabla^{m}S_{\ez}f|^{p_{0}} dw\r)^{1/^{p_{0}} },$$ where the implicit constant is independent of $\ez.$ Letting $\ez \to 0$ in the above inequality and using an analogous argument to the one after \eqref{eq: cd1032}, we can derive \eqref{eq: cd20}. This is justified because the series $\sum_{j \geq 1} g(j)$ is finite. Consequently, applying Theorem \ref{theorem: cd18} with $v\equiv 1$ leads to \eqref{eq: cd1024} for every $f\in \D$ and any $p\in (2, \;q_{+}).$

At this stage, we are left only to show that \eqref{eq: cd1025} holds for all $q,$ $q_{-}<p<q_{+},$ and $v\in A_{\frac{p}{q_{-}(\L_{w})}}(w)\cap \mbox{RH}_{(\frac{q_{+}(\L_{w})}{p})'}(w),$ as \eqref{eq: cd1024} in the interval $(2,\; q_{+})$ will follow directly from \eqref{eq: cd1025} by taking $v\equiv 1.$ 

Recall that, by Proposition \ref{proposition: cd3000}, there exist $p_{0}, q_{0}$ such that $$\tilde{q}_{-}<p_{0}<\min\{p, 2\}\leq \max\{p, 2\} <q_{0}<q_{+}\quad \mbox{and}\quad v\in A_{\frac{p}{p_{0}}}(w)\cap \mbox{RH}_{(\frac{q_{0}}{p})'}(w).$$
From this, along with \cite[Lemma 4.4]{AM2}, it follows that $$u:=v^{1-p'}\in A_{\frac{p'}{ q_{0}' }}(w)\cap \mbox{RH}_{(\frac{p_{0}'}{p'})'}(w).$$ Then, \eqref{eq: cd1025} holds by duality, once we prove that \begin{equation}\label{eq: cd10380}\|\T^{*}\overrightarrow{f}\|_{L^{p'}(\rz; udw)}\lesssim \|\overrightarrow{f}\|_{L^{p'}(\rz; (\cc)^{m}, udw)}.\end{equation} 
In contrast to the strategy used in Propositions \ref{proposition: cd15} and \ref{proposition: cd81}, the proof of \eqref{eq: cd10380} requires us to appeal to Theorem \ref{theorem: cd1035} rather than Theorem \ref{theorem: cd18}.






We now verify that the conditions of Theorem \ref{theorem: cd1035} are satisfied. Let $\overrightarrow{f}\in L^{\infty}_{c}(\rz; (\cc)^{m}),$ and set $F:=|\T^{*}\overrightarrow{f}|^{q_{0}'}.$ Since $2 < q_0 < q_+$, it follows from \eqref{eq: cd1024} (applied with exponent $q_0$) and a duality argument that $F \in L^1(w)$. On the other hand, fix a ball $B$ of radius $r$ and let $\A_{r}$ be defined as before. Then $$F\leq 2^{q_{0}'-1}|(I-\A_{r})^{*}(\T^{*}\overrightarrow{f})|^{q_{0}'}+2^{q_{0}'-1}|\A_{r}^{*} (\T^{*}\overrightarrow{f})|^{q_{0}'}:=G_{B}+H_{B},$$ where the adjoint is taken in $L^{2}(w).$ This verifies condition \eqref{eq: cd5000} in Theorem \ref{theorem: cd1035}. 

Setting $G:=M_{w}(|\overrightarrow{f}|^{q_{0}'})$ and $q:=\frac{p_{0}'}{q_{0}'},$ we would like to prove \eqref{eq: cd1036} for the choice of $G$ and $q.$ To achieve this, we first note that  $\A_{r} \in \mho(L^{p_{0}}(w)\to L^{q_{0}}(w))$ by Proposition \ref{proposition: cd100}. Then, using this and duality, we can find a function $g\in L_{c}^{p_{0}}(B, \frac{dw}{w(B)})$ with norm 1 such that for all $x\in B,$ 
\begin{equation*}\begin{aligned} \l(\fint_{B}H_{B}^{q}dw\r)^{\frac{1}{qq_{0}'}}&\lesssim \frac{1}{w(B)}\int_{\rz}|\T^{*}\overrightarrow{f}| |\A_{r} g|dw\\[4pt] 
&\lesssim \sum_{j\geq 1} 2^{jD}\l(\fint_{C_{j}(B)}|\T^{*}\overrightarrow{f}|^{q_{0}'}dw\r)^{\frac{1}{q_{0}'}} \l(\fint_{C_{j}(B)}|\A_{r} g|^{q_{0}}dw\r)^{\frac{1}{q_{0}}}\\[4pt]
&\lesssim  M_{w}(F)^{\frac{1}{q_{0}'} }(x)(x)\sum_{j\geq 1} 2^{j(D+\theta_{1}+\theta_{2})} e^{-c(2^{j\frac{2m}{2m-1}})} \l(\fint_{B}| g|^{p_{0}}dw\r)^{\frac{1}{p_{0}}}\lesssim  M_{w}(F)^{\frac{1}{q_{0}'} }(x).\end{aligned}\end{equation*} 
Similarly, there exists $g\in L_{c}^{q_{0}}(B, \frac{dw}{w(B)})$ with norm 1 such that for all $x\in B,$ 
\begin{equation}\label{eq: cd5001}\begin{aligned} \l(\fint_{B}G_{B}^{q}dw\r)^{\frac{1}{q_{0}'}}\lesssim  M_{w}(|\overrightarrow{f}|^{q_{0}'})(x)^{1/q_{0}'}\sum_{j\geq 1} 2^{jD} \l(\fint_{C_{j}(B)}|\T(I-\A_{r}) g|^{q_{0}}dw\r)^{\frac{1}{q_{0}}}.\end{aligned}\end{equation} We proceed to analyze each term in the preceding sum. For $j=1,$ the $L^{q_{0}}(w)$-boundedness of $\T$ (by \eqref{eq: cd1024}) and of $e^{-r^{2m}\L_{w}}$ (as $q_{0} \in \widetilde{\J}(\L_{w})$) implies that \begin{equation}\label{eq: cd1038}\l(\fint_{4B}|\nabla^{m} \L_{w}^{-1/2}(I-e^{-r^{2m} \L_{w}})^{N}g|^{q_{0}} dw\r)^{1/q_{0}}\lesssim \l(\fint_{B}|g|^{q_{0}} dw\r)^{1/^{q_{0}} } =1.\end{equation} For $j \geq 2$,  we again employ the integral representation \eqref{eq: cd1028} and, by estimating as in \eqref{eq: cd1029} but with the roles of $B$ and $C_j(B)$ interchanged, conclude that \begin{equation*}\begin{aligned} &\l(\fint_{C_{j}(B)}\bigg |\int_{\Gamma_{\pm}}t^{1/2}\nabla^{m}e^{-z\L_{w}}g\eta_{\pm}(z, t)dz\bigg |^{q_{0}}dw\r)^{\frac{1}{q_{0}}}\\[4pt]
&\quad\quad\quad\quad\lesssim \int_{\Gamma_{\pm}}\l(\fint_{C_{j}(B)}|  z^{1/2} \nabla^{m} e^{-z\L_{w}}g |^{q_{0}}dw\r)^{\frac{1}{q_{0}}}\frac{t^{1/2}}{|z|^{1/2}}|\eta_{\pm}(z, t)| |dz|\\[4pt]
&\quad\quad\quad\quad\lesssim 2^{j\theta_{1}}\l(\fint_{B}|g|^{q_{0}}dw\r)^{1/q_{0}}\int_{\Gamma_{\pm}}\Upsilon\l(\frac{2^{j} r}{|z|^{1/2m}}\r)^{\theta_{2}}e^{-c\l(\frac{2^{j}r}{|z|^{\frac{1}{2m}}}\r)^{\frac{2m}{2m-1}}}\frac{t^{1/2}}{|z|^{1/2}}|\eta_{\pm}(z, t)| |dz|\\[4pt]
&\quad\quad\quad\quad\lesssim \l(\fint_{B}|g|^{q_{0}}dw\r)^{1/q_{0}} 2^{j\theta_{1}}\int_{0}^{\infty}\frac{r^{2mN}}{(s+t)^{N+1}}\Upsilon\l(\frac{2^{j} r}{s^{1/2m}}\r)^{\theta_{2}}e^{-c\l(\frac{2^{j}r}{s^{\frac{1}{2m}}}\r)^{\frac{2m}{2m-1}}}\frac{t^{1/2}}{s^{1/2}} ds.\end{aligned}
\end{equation*} 
With this inequality in hand, the argument leading to \eqref{eq: cd1030} yields  
\begin{equation}\label{eq: cd1039}\l(\fint_{C_{j}(B)}|\T(I-\A_{r}) g|^{q_{0}}dw\r)^{\frac{1}{q_{0}}}\lesssim 2^{j(\theta_{1}-2mN)}\end{equation} provided $2mN>\theta_{2}+1.$ Gathering \eqref{eq: cd5001}, \eqref{eq: cd1038} and \eqref{eq: cd1039} we arrive at $$\l(\fint_{B}G_{B}^{q}dw\r)^{\frac{1}{q_{0}'}}\lesssim  M_{w}(|\overrightarrow{f}|^{q_{0}'})(x)^{1/q_{0}'}\sum_{j\geq 1} 2^{j(D+\theta_{1}-2mN)} \lesssim M_{w}(|\overrightarrow{f}|^{q_{0}'})(x)^{1/q_{0}'}=G(x)^{1/q_{0}'},$$ provided $2mN>D+\theta_{1}+\theta_{2}+1.$ This gives \eqref{eq: cd1036} with $q=\frac{p_{0}'}{q_{0}'}$ and $G=M_{w}(|\overrightarrow{f}|^{q_{0}'}).$

Since $u\in \mbox{RH}_{(\frac{p_{0}'}{p'})'}(w),$ by Proposition \ref{proposition: cd3000}, there exists a $s<\frac{p_{0}'}{p'}$ such that $u\in \mbox{RH}_{s'}(w).$ If we set $t=\frac{p' }{q_{0}'}<q/s,$ then $u\in A_{t}(w),$ so $M_{w}$ is bounded on $L^{t}(udw).$ From this and \eqref{eq: cd1037}, it follows that
\begin{equation*}\begin{aligned} \|\T^{*}\overrightarrow{f}\|_{L^{p'}(udw)}^{q_{0}'}&\leq \|M_{w}(F)\|_{L^{t}(udw)} \\[4pt]
& \lesssim \|G\|_{L^{t}(udw)}\approx \|M_{w}(|\overrightarrow{f}|^{q_{0}'})\|_{L^{t}(udw)} \lesssim \|\overrightarrow{f}\|_{L^{p'}(udw)}^{q_{0}'}.\end{aligned}\end{equation*} This proves \eqref{eq: cd10380}, thus completing our proof. 

\hfill$\Box$


\section{Square function estimates for $e^{-t\L_{w}}$ and $t^{1/2}\nabla^{m} e^{-t\L_{w}}$}
As an application of the main results established above, we can study the weighted $L^{p}$ norm inequalities for two vertical square functions associated with the semigroups $e^{-t\L_{w}}$ and $t^{1/2}\nabla^{m} e^{-t\L_{w}}.$ These are defined, respectively, as $$g_{\L_{w}}f(x):=\l(\int_{0}^{\infty}|(t\L_{w})^{1/2}e^{-t\L_{w}}f(x)|^{2}\frac{dt }{t}\r)^{1/2}$$ and $$G_{\L_{w}}f(x):=\l(\int_{0}^{\infty}|t^{1/2} \nabla^{m} e^{-t\L_{w}}f(x)|^{2}\frac{dt}{t}\r)^{1/2}.$$ More precisely, the goal of this section is to prove the following two propositions. 

\begin{proposition}\label{proposition: cd50}\;  
Assume $p_{-}<p<p_{+}$. Then \begin{equation}\label{eq: cd51} \|g_{\L_{w}}f\|_{L^{p}(w)} \approx \|f\|_{L^{p}(w)}.\end{equation}
Conversely, if \eqref{eq: cd51} holds for some $p$, then $p \in \widetilde{\J}(\L_{w})$. In other words, the interior of the interval on which \eqref{eq: cd51} holds is exactly $(p_{-}, p_{+})$. Moreover, \begin{equation}\label{eq: cd52} \|g_{\L_{w}}f\|_{L^{p}(vdw)} \approx \|f\|_{L^{p}(vdw)}\end{equation} holds for any $v\in A_{\frac{p}{p_{-}}}(w)\cap \mbox{RH}_{(\frac{p_{+}}{p})'}(w).$
\end{proposition}

\begin{proposition}\label{proposition: cd1040}\;  Assume $q_{-}(\L_{w})<p<q_{+}(\L_{w}).$ Then 
 \begin{equation}\label{eq: cd1041} \|G_{\L_{w}} f\|_{L^{p}(w)}\lesssim \| f\|_{L^{p}(w)}\end{equation} and for any $v\in A_{\frac{p}{q_{-}(\L_{w})}}(w)\cap \mbox{RH}_{(\frac{q_{+}(\L_{w})}{p})'}(w),$ \begin{equation}\label{eq: cd1042}\|G_{\L_{w}} f\|_{L^{p}(vdw)}\lesssim \|f\|_{L^{p}(vdw)}.\end{equation} 
\end{proposition}

Central to the proofs of Propositions \ref{proposition: cd50} and \ref{proposition: cd1040} is the following Lemma \ref{lemma: cd53} on Hilbert-valued extensions. This requires some notation: let $\HH$ denote the Hilbert space $L^{2}((0, \infty), \frac{dt}{t}),$ endowed with the norm $$|||f|||=\l(\int_{0}^{\infty}|f(t)|^{2}\frac{d t}{t}\r)^{1/2}.$$ Furthermore, given a Borel measure $\mu$ on $\rz,$ we define $L^{p}_{\HH}(\mu)$ as the space of $\HH$-valued functions with the norm 
$$\|f\|_{L^{p}_{\HH}(w)}:=\l(\int_{\rz}|||  f(x, \cdot)|||^{p} d\mu\r)^{1/p}.$$ 


\begin{lemma}\label{lemma: cd53}\;(\cite[Lemma 7.4]{AM1})\;  Let $\D$ be a subspace of $\G,$ the space of measurable functions in $\rz,$ and let $S, T$ be two linear operators from $\D$ into $\G.$ Fix $1\leq p\leq q<\infty$ and suppose there exists $C_{0}>0$ such that for all $f\in \D,$ $$\|T f\|_{L^{q}(\mu)}\leq C_{0}\sum_{j\geq 1}\alpha_{j}\|S f\|_{L^{p}(F_{j}, \mu)},$$ where the $F_{j}$ are measurable subsets of $\rz$ and $\alpha_{j}\geq 0.$ Then there is an $\HH$-valued inequality with the same constant: for all $f: \rz\times (0, \infty)\to \cc$ such that for almost all $t>0, f(\cdot, t)\in \D,$ $$\|T f\|_{L^{q}_{\HH}(\mu)}\leq C_{0} \sum_{j\geq 1}\alpha_{j}\|S f\|_{L^{p}_{\HH}(F_{j}, \mu)}.$$ 
\end{lemma}
The extension of a linear operator $T$ on $\cc$-valued functions to $\HH$-valued functions is defined by $$(Th)(x, t):=T(h(\cdot, t))(x)\quad \mbox{for}\; x\in \rz \;\mbox{and}\; t>0.$$ This means $t$ is treated as a parameter, and $T$ acts only on the spatial variable. We begin with the proof of Proposition \ref{proposition: cd50}. As it is very similar to \cite[Proposition 5.1]{AM}, we outline the key differences that arise in the higher-order setting.

\noindent{\textbf{Proof of Proposition \ref{proposition: cd50}}:} \; Choose $0<\mu<\frac{\pi}{2},$ and let $\phi(z):=z^{1/2}e^{-z}.$ Then $\phi \in \H_{0}^{\infty}(\Sigma_{\mu}),$ so it follows from \eqref{eq: cd3007} that \begin{equation}\label{eq: cd54} \|g_{\L_{w}}f\|_{L^{2}(w)}=\l(\int_{0}^{\infty}\|\phi(t\L_{w})f\|_{L^{2}(w)}^{2}\frac{dt }{t}\r)^{1/2} \lesssim \|f\|_{L^{2}(w)}.\end{equation} Our first goal is to apply Theorem \ref{theorem: cd21}, in view of \eqref{eq: cd54}, to establish the inequality
\begin{equation}\label{eq: cd55}
 \|g_{\L_{w}}f\|_{L^{p}(w)}\lesssim \|f\|_{L^{p}(w)}, \quad p_{-} < p < 2.
\end{equation}


Let $q_{0}:=2$ and fix $p,$ $p_{-}<p<q_{0}.$ We use the operator $\A_{r}$, defined as before. Then, by Proposition \ref{proposition: cd15}, $\A_{r}$ is bounded on $L^{q_{0}}(w)$ for each $N.$ With these preparations, we now show that for any $f\in L^{\infty}_{c}$ with $\supp\; f \subset B$ and $j\geq 2,$ \eqref{eq: cd22} holds with $\T=g_{\L_{w}}.$ To do so, we set $\phi(z, t):=(t z)^{1/2}e^{-tz}(1-e^{-r^{2m} z})^{N}.$ Clearly, $\phi(\cdot, t)\in \H_{0}^{\infty}(\Sigma_{\mu})$ if $0<\mu<\frac{\pi}{2},$ and $$(t\L_{w})^{1/2}e^{-t\L_{w}}(I-\A_{r})f=\phi(\L_{w}, t)f.$$ Moreover, we can rewrite $\phi(\L_{w}, t)f$ in the form given by \eqref{eq: cd3004}-\eqref{eq: cd3005},
with functions $\eta_{\pm}(z, t)$ satisfying \begin{equation}\label{eq: cd1089}|\eta_{\pm}(z, t)|\lesssim \frac{t^{1/2} r^{2mN}}{(|z|+t)^{N+3/2}},\quad \mbox{for any}\; z\in \Gamma_{\pm},\end{equation} where $0<\V<\theta<\nu<\mu<\frac{\pi}{2}.$ A direct consequence of \eqref{eq: cd1089} is that \begin{equation}\label{eq: cd570} ||| \eta(z, \cdot)||| \lesssim \frac{r^{2mN}}{|z|^{N+1}}.\end{equation} Using \eqref{eq: cd570}, together with the fact that $e^{-z\L_{w}} \in \mho(L^{p}(w) \to L^{p}(w))$ for $z \in \Gamma_{\pm}$ (which is guaranteed by Corollary \ref{corollary: cd14}), we deduce that \begin{equation}\label{eq: cd60}\begin{aligned}&\l(\fint_{C_{j}(B)}|g_{\L_{w}}(I-\A_{r})f|^{p}dw\r)^{\frac{1}{p}}\\[4pt]
&\quad\quad \lesssim \l(\fint_{C_{j}(B)}\bigg |\l(\int_{0}^{\infty}\bigg|\int_{\Gamma_{\pm}}e^{-z\L_{w}}f\eta_{\pm}(z, t)dz\bigg|^{2}\frac{dt}{t}\r)^{1/2}\bigg |^{p}dw\r)^{\frac{1}{p}}\\[4pt]
&\quad\quad\lesssim  \l(\fint_{C_{j}(B)}\bigg |\int_{\Gamma_{\pm}}|e^{-z\L_{w}}f| |||\eta_{\pm}(z, \cdot)||||dz| \bigg |^{p}dw\r)^{\frac{1}{p}}\\[4pt]
&\quad\quad\lesssim  \int_{\Gamma_{\pm}} \l(\fint_{C_{j}(B)}|e^{-z\L_{w}}f|^{p} dw\r)^{\frac{1}{p}}\frac{r^{2mN}}{|z|^{N+1}}|dz| \\[4pt]
&\quad\quad\lesssim 2^{j\theta_{1}}\l(\fint_{B}|f|^{p}dw\r)^{1/p}\int_{\Gamma_{\pm}}\frac{r^{2mN}}{|z|^{N+1}}\Upsilon\l(\frac{2^{j} r}{|z|^{1/2m}}\r)^{\theta_{2}}e^{-c\l(\frac{2^{j}r}{|z|^{\frac{1}{2m}}} \r)^{\frac{2m}{2m-1}}}d|z|\\[4pt]
&\quad\quad\approx \l(\fint_{B}|f|^{p}dw\r)^{1/p} 2^{j(\theta_{1}-2mN)}\int_{0}^{\infty}\Upsilon(\tau)^{\theta_{2}} \tau^{2mN} e^{-c\tau^{\frac{2m}{2m-1}} }\frac{d\tau}{\tau} \\[4pt]
& \quad\quad\lesssim 2^{j(\theta_{1}-2mN)}\l(\fint_{B}|f|^{p}dw\r)^{1/p} ,\end{aligned}
\end{equation} where in the last step we also used the assumption $2mN>\theta_{2}+1.$ Furthermore, if we choose $N$ large enough so that $2mN>\theta_{1}+\theta_{2}+D+1,$ then $g_{\L_{w}}$ satisfies \eqref{eq: cd22} with $g(j):=C2^{j(\theta_{1}-2mN)}.$ On the other hand, \eqref{eq: cd23} has already been established in \eqref{eq: cd30} with $g(j) = C 2^{j(\theta_1 + \theta_2)} e^{-c (2^j)^{\frac{2m}{2m-1}}}$. Hence, applying Theorem \ref{theorem: cd21} yields \eqref{eq: cd55} for any $p$ with $p_{-}<p<2.$


In accordance with the strategy of Proposition \ref{proposition: cd15}, it remains to prove \eqref{eq: cd52} on the interval $(p_{-}, p_{+})$ by exploiting Theorem \ref{theorem: cd18}. 
We first prove condition \eqref{eq: cd19}. For this aim, we recall \eqref{eq: cd5005} and repeat the argument in  \eqref{eq: cd60} to conclude that for $j\geq 1,$
\begin{equation}\label{eq: cd61}\begin{aligned}
&\l(\fint_{B}|g_{\L_{w}}(I-\A_{r})f_{j}|^{p_{0}}dw\r)^{\frac{1}{p_{0}}}\quad (f_{j}:=f1_{C_{j}(B)})\\[4pt]
&\quad\quad\quad\lesssim  \int_{\Gamma_{\pm}} \l(\fint_{B}|e^{-z\L_{w}}f_{j}|^{p_{0}} dw\r)^{\frac{1}{p_{0}}}\frac{r^{2mN}}{|z|^{N+1}}|dz|\\[4pt]
&\quad\quad\quad\lesssim 2^{j\theta_{1}}\l(\fint_{C_{j}(B)}|f|^{p_{0}}dw\r)^{1/p_{0}}\int_{\Gamma_{\pm}}\frac{r^{2mN}}{|z|^{N+1}}\Upsilon\l(\frac{2^{j} r}{|z|^{1/2m}}\r)^{\theta_{2}}e^{-c\l(\frac{2^{j}r}{|z|^{\frac{1}{2m}}}\r)^{\frac{2m}{2m-1}}}d|z|\\[4pt]
&\quad\quad\quad\approx \l(\fint_{C_{j}(B)}|f|^{p_{0}}dw\r)^{1/p_{0}} 2^{j(\theta_{1}-2mN)}\int_{0}^{\infty}\Upsilon(\tau)^{\theta_{2}} \tau^{2mN} e^{-c\tau^{\frac{2m}{2m-1}} }\frac{d\tau}{\tau} \\[4pt]
&\quad\quad\quad\lesssim 2^{j(\theta_{1}-2mN)}\l(\fint_{2^{j+1}B}|f|^{p_{0}}dw\r)^{1/p_{0}}.\end{aligned}
\end{equation} Summing \eqref{eq: cd61} over all $j \geq 1$ and taking $g(j) := C 2^{j(\theta_1 - 2mN)}$ for sufficiently large $N$, we obtain the estimate \eqref{eq: cd19} for $\T = g_{\L_w}$ and $\S = I$.


We begin the proof of \eqref{eq: cd20} by invoking Proposition \ref{proposition: cd12}. This gives, for $1\leq k\leq N,$ $j\geq 1$ and $\supp\; g\subset C_{j}(B),$
\begin{equation}\label{eq: cd62}\l(\fint_{B}|e^{-kr^{2m} \L_{w}} g|^{q_{0}}dw\r)^{\frac{1}{q_{0}}}\leq C_{0} 2^{j(\theta_{1}+\theta_{2})} e^{-c(2^{j\frac{2m}{2m-1}})}\l(\fint_{C_{j}(B)}|f|^{p_{0}}dw\r)^{\frac{1}{p_{0}}},\end{equation} with $c, C_{0}$ independent of $k.$ Setting $T: L^{p_{0}}(w)\to L^{q_{0}}(w)$ as $$Tg =(C_{0}2^{j(\theta_{1}+\theta_{2})} e^{-c(2^{j\frac{2m}{2m-1}})})^{-1}\frac{w(2^{j+1}B)^{1/p_{0}}}{w(B)^{1/q_{0}}}1_{B}e^{-kr^{2m}\L_{w}}(g1_{C_{j}(B)}),$$ we then have, by \eqref{eq: cd62}, $$\|T g\|_{L^{q_{0}}(w)}\leq \l(\int_{C_{j}(B)}|f|^{p_{0}}dw\r)^{1/p_{0}}=\l(\int_{C_{j}(B)}|Sf|^{p_{0}}dw\r)^{1/p_{0}},$$ where $S=I.$ This allows us to apply Lemma \ref{lemma: cd53} to obtain that for all $g\in L^{p_{0}}_{\HH}(w)$ with $\supp\; g(\cdot, t)\subset C_{j}(B)$ ($t>0$), 
\begin{equation}\label{eq: cd63}\l(\fint_{B}||| e^{-kr^{2m} \L_{w}} g(x, \cdot)   |||^{q_{0}}dw\r)^{\frac{1}{q_{0}}}\leq C_{0} 2^{j(\theta_{1}+\theta_{2})} e^{-c(2^{j\frac{2m}{2m-1}})}\l(\fint_{C_{j}(B)}||| g(x, \cdot)   |||^{p_{0}}dw\r)^{\frac{1}{p_{0}}}.\end{equation} From \eqref{eq: cd63},  it follows that for any $g\in L^{p_{0}}_{\HH}(w),$ \begin{equation}\label{eq: cd64}\begin{aligned}
\l(\fint_{B}||| e^{-kr^{2m} \L_{w}} g(x, \cdot)   |||^{q_{0}}dw\r)^{\frac{1}{q_{0}} }&\lesssim \sum_{j\geq 1}\l(\fint_{B}||| e^{-kr^{2m} \L_{w}} g_{j}(x, \cdot)   |||^{q_{0}}dw\r)^{\frac{1}{q_{0}}}\\[4pt]
&\lesssim \sum_{j\geq 1} 2^{j(\theta_{1}+\theta_{2})} e^{-c(2^{j\frac{2m}{2m-1}})}\l(\fint_{C_{j}(B)}||| g(x, \cdot)   |||^{p_{0}}dw\r)^{\frac{1}{p_{0}}},
 \end{aligned} 
\end{equation} where $$g(x, t)=\sum_{j\geq 1}g(x, t)1_{C_{j}(B)}(x):=\sum_{j\geq 1}g_{j}(x, t).$$ 
In particular, we choose $g(x, t):=(t \L_{w})^{1/2}e^{-t\L_{w}}f(x),$ so $g_{\L_{w}}f(x)=||| g(x, \cdot) |||.$ We note that $p_{-} < p_0 < 2$ and, by \eqref{eq: cd55}, $g \in L^{p_0}_{\HH}(w)$. Moreover, since $(t \L_w)^{1/2}e^{-t\L_w}$ and $e^{-k r^{2m} \L_w}$ commute, we can write $$g_{\L_{w}}(e^{-kr^{2m}\L_{w}})f(x)=||| e^{-kr^{2m}\L_{w}}g(x, \cdot) |||.$$ Consequently, an application of \eqref{eq: cd64} leads to \begin{equation}\label{eq: cd65}\begin{aligned}
\l(\fint_{B}|g_{\L_{w}} \A_{r}f|^{q_{0}}dw\r)^{1/q_{0}}\lesssim 
\sum_{j\geq 1} 2^{j(\theta_{1}+\theta_{2})} e^{-c(2^{j})^{\frac{2m}{2m-1}}}\l(\fint_{2^{j+1}(B)}|g_{\L_{w}}f|^{p_{0}}dw\r)^{1/p_{0}},
 \end{aligned} 
\end{equation} which implies \eqref{eq: cd20} with $\T = g_{\L_w}.$ Therefore, Theorem \ref{theorem: cd18} applies, and \eqref{eq: cd52} is concluded. 

A careful examination of the preceding arguments reveals that they actually prove a more general result: the upper bounds in \eqref{eq: cd51}-\eqref{eq: cd52} remain valid when $g_{\L_w}$ is replaced by either $g^{\phi}_{\L_{w}}$ or $g^{\phi}_{\L_{w}, d}.$ Here, these generalized square functions are defined for any holomorphic function $\phi$ on the sector $\Sigma_{\pi/2}$ by $$g^{\phi}_{\L_{w}}f(x):=\l(\int_{0}^{\infty}|\phi(t \L_{w})f(x)|^{2}\frac{dt }{t}\r)^{1/2} \quad \mbox{and}\quad g^{\phi}_{\L_{w}, d}f(x):=\l(\sum_{k\in \Z}|\phi(2^{2km}\L_{w})|^{2}\r)^{1/2}$$ provided that $\phi$ satisfies the growth condition $$|\phi(z)| \lesssim |z|^{1/2} e^{-c|z|}\quad \mbox{uniformly on}\; \Sigma_{\mu}\; \mbox{for any}\; 0\leq \mu<\frac{\pi}{2}.$$ From the upper bound in \eqref{eq: cd51} for $g^{\phi}_{\L_{w}, d},$ it follows that for any sequence of functions $\{\beta_{k}\}_{k\in \Z}$ and $p\in (p_{-}, p_{+}),$ \begin{equation}\label{eq: cd6000} \|\sum_{k\in \Z} \psi(2^{2mk} \L_{w})\beta_{k}\|_{L^{p}(w)}\lesssim \|\l(\sum_{k\in \Z} |\beta_{k}|^{2}\r)^{1/2}\|_{L^{p}(w)},\end{equation} where $$\psi(z):=\frac{1}{\pi^{1/2}}\int_{ 1}^{\infty}ze^{-tz}\frac{dt}{t^{1/2}}.$$ A detailed proof of \eqref{eq: cd6000} can be found in \cite[Proposition 5.14]{DMR}.


We now prove the converse of \eqref{eq: cd51}-\eqref{eq: cd52}. Since the lower bound in \eqref{eq: cd51} is the special case of \eqref{eq: cd52} with $v \equiv 1$, we focus on proving the lower bound in \eqref{eq: cd52}. Note that Lemma \ref{lemma: h5} gives the duality relation: \begin{equation}\label{eq: cd80} 
p_{\pm}(\L_{w})'=p_{\mp}(\L_{w}^{*}).
\end{equation} Combining \eqref{eq: cd80} and \cite[Lemma 4.4]{AM2}, we see that for all $p\in (p_{-}, p_{+})$ and $v\in A_{\frac{p}{p_{-}}}(w)\cap \mbox{RH}_{(\frac{p_{+}}{p})'}(w),$ $$v^{1-p'}\in A_{\frac{p'}{ p_{-}(\L^{*}_{w}) }}(w)\cap \mbox{RH}_{(\frac{p_{+}(\L^{*}_{w}) }{p'})'}(w).$$ The remainder of the proof follows verbatim from the arguments on \cite[pp. 632-633]{DMR}, thereby completing the proof of Proposition \ref{proposition: cd50}.

\hfill$\Box$

\noindent{\textbf{Proof of Proposition \ref{proposition: cd1040}}:} \; Using Proposition \ref{proposition: cd1023} and Lemma \ref{lemma: cd53}, \eqref{eq: cd1042} can be reduced to \eqref{eq: cd51}; see \cite[Proposition 10.1]{DMR} for a proof. Once \eqref{eq: cd1042} is proved, \eqref{eq: cd1041} follows readily by taking $v\equiv 1.$  

\hfill$\Box$




We conclude this section by stating a reverse inequality for $G_{\L_{w}}$, although it will not be used in the subsequent proofs, even in our higher-order extension of \cite{ACMP}.

\begin{proposition}\label{proposition: cd1090}\; Let $q_{+}(\Delta_{w, m})'<p<\infty,$  where $\Delta_{w, m}:=(-1)^{m}w^{-1}\mbox{div}_{m}(w\nabla^{m}).$ Then 
 \begin{equation}\label{eq: cd1091} \| f\|_{L^{p}(w)}\lesssim \| G_{\L_{w}}f\|_{L^{p}(w)}.\end{equation} Furthermore, if $v\in A_{\frac{p}{q_{+}(\Delta_{w, m})'}}(w),$ \begin{equation}\label{eq: cd1092}\| f\|_{L^{p}(vdw)}\lesssim \|G_{\L_{w}} f\|_{L^{p}(vdw)}.\end{equation} 
\end{proposition}
The proof follows \cite[Proposition 10.4]{DMR} almost identically, relying on the property that $e^{-t\Delta_{w, m}}\in \mho(L^{1}(w)\to L^{\infty}(w))$. This property is equivalent to the Gaussian estimate for the kernel of $e^{-t\Delta_{w, m}}$, as shown in \cite[Proposition 2.2]{AM}. A forthcoming work will be devoted to a more general result, which can be viewed as either a higher-order generalization of \cite[Theorem 1]{CR2} or a weighted analogue of \cite[Definition 9]{AQ}. This result implies Proposition \ref{proposition: cd1090} and is summarized below:


\begin{theorem}\label{theorem: cd1093}\; If $\{a_{\alpha, \beta}\}_{|\alpha|=|\beta|=m}\in \Ez(w, c_{1}, c_{2}),$ then there exists a heat kernel $K_{t}(x, y)$ associated to $e^{-t\L_{w}}$ such that, for some $\mu=\nu+l$ with $l\in \{0, 1,...m-1\}$ and $\nu\in (0, 1),$ and for any $f\in C_{0}^{\infty}(\rz),$ all $t>0,$ all $x, y\in \rz$ and all multi-index $\gamma$ 
 \begin{equation}\label{eq: cd1094} |D_{x}^{\gamma}K_{t}(x, y)|+|D_{y}^{\gamma}K_{t}(x, y)|\leq \frac{C}{ w(B_{t^{1/2m}}(x))^{\frac{1}{2}+\frac{|\gamma|}{2n}} w(B_{t^{1/2m}}(y))^{\frac{1}{2}+\frac{|\gamma|}{2n}}}g_{m, c}\l(\frac{|x-y|}{t^{\frac{1}{2m}}}\r),\end{equation} when $|\gamma|\leq l,$
\begin{equation}\label{eq: cd1095}\begin{aligned} &|D_{x}^{\gamma}K_{t}(x+h, y)-D_{x}^{\gamma}K_{t}(x, y)|+|D_{y}^{\gamma}K_{t}(x, y+h)-D_{y}^{\gamma}K_{t}(x, y)| \\[4pt]
&\quad\quad \leq  \frac{C}{ w(B_{t^{1/2m}}(x))^{\frac{1}{2}+\frac{|\gamma|}{2n}} w(B_{t^{1/2m}}(y))^{\frac{1}{2}+\frac{|\gamma|}{2n}}}\l(\frac{|h|}{t^{1/2m}+|x-y|}\r)^{\nu} g_{m, c}\l(\frac{|x-y|}{t^{\frac{1}{2m}}}\r),
\end{aligned}\end{equation}
 when $|\gamma|=l$ and $2|h|\leq t^{1/2m}+|x-y|,$ where $g_{m, c}(s):=e^{-cs^{\frac{2m}{2m-1}}}$ for $s>0.$\end{theorem}







\section{Unweighted $L^{p}$ Kato estimates and their applications}
This section constitutes the culmination of our argument. We will derive unweighted $L^p$ estimates for operators associated to $\L_w$-such as the semigroup, its gradients, Riesz transforms, functional calculus, and square functions-when $p$ is near 2. We achieve this by imposing additional requirements on $w \in A_2$, which permit us to set $v \equiv w^{-1}$ in the $L^p(v,dw)$-estimates established in the previous sections. These unweighted $L^{p}$ estimates are then employed to solve the corresponding $L^{p}(\rz)$-Dirichlet, regularity and Neumann boundary value problems.


\begin{theorem}\label{theorem: cd1048}\;  Let $w\in A_{2},$ and $\eta\geq 1$ with $|p-2|<\ez\;(0<\ez<\frac{2m}{n+2m}),$ and assume $1\leq  r_{w}<1+\frac{p m}{n}$ and $s_{w}>\frac{nr_{w}}{p m}+1.$ Then $e^{-t\L_{w}}: L^{p}(\rz)\to L^{p}(\rz)$ is uniformly bounded for all $t>0.$ Likewise, both $\phi(\L_{w})$ (with $\phi$ bounded and holomorphic on $\Sigma_{\mu}, \mu\in (\V, \pi)$) and $g_{\L_{w}}$ are bounded operators on $L^{p}(\rz).$ More generally, these $L^{p}$ bounds remian valid under either of the following conditions: $(i)$\; $w\in A_{r}\cap \mbox{RH}_{\frac{n r}{p m}+1}$ with $1<r< 1+\frac{p m}{n};$ $(ii)$\; $w$ is a power weight $w_{\alpha}(x):=|x|^{\alpha}$ with $-\frac{p mn}{n+p m}<\alpha<p m.$

 \end{theorem}
{\it Proof.}\quad Let $p_{0}=(p^{*, m}_{w})', q_{0}=p^{*, m}_{w},$ and set $v=w^{-1}.$ Then, from Proposition \ref{proposition: cd12} and $0<\ez<\frac{2m}{n+2m},$ it holds that $$p_{-}\leq p_{0}<p<q_{0}\leq p_{+}.$$ Hence,
by Corollary \ref{corollary: cd13}, we have $e^{-t\L_{w}}\in \mho(L^{p}(\rz)\to L^{p}(\rz))$ whenever $w^{-1}\in A_{\frac{p}{p_{0}}}(w)\cap \mbox{RH}_{(\frac{q_{0}}{p})'}(w).$ Note that property $(x)$ of Proposition \ref{proposition: cd3000} implies $$w^{-1}\in A_{\frac{p}{p_{0}}}(w)\cap \mbox{RH}_{(\frac{q_{0}}{p})'}(w)\Longleftrightarrow w\in A_{\frac{q_{0}}{p}}\cap \mbox{RH}_{(\frac{p}{p_{0}})'}.$$ Moreover, by recalling the definition of $p^{*, m}_{w},$ we see 
$$\frac{q_{0}}{p}=\frac{nr_{w}}{nr_{w}-p m} \quad\mbox{and}\quad (\frac{p}{p_{0}})'=\frac{nr_{w}}{p m}+1.$$ Clearly, it follows from $r_{w}<1+\frac{p m}{n}$ that $w\in A_{\frac{q_{0}}{p}},$ and from $s_{w}>\frac{nr_{w}}{\eta m}+1$ that $w\in \mbox{RH}_{(\frac{\eta}{p_{0}})'}.$ 

If $w\in A_{r}\cap \mbox{RH}_{\frac{n r}{p m}+1}$ and $1<r< 1+\frac{p m}{n},$ it is easy to see that $r_{w}\leq r<1+\frac{p m}{n}$ and $s_{w}>\frac{n r}{p m}+1\geq \frac{nr_{w}}{p m}+1.$ Consequently, applying Lemma \ref{lemma: h6} yields that $e^{-t\L_{w}}$ is uniformly bounded on $L^{p}(\rz).$ The case of power weights is immediate from \eqref{eq: cd1043}, as $-\frac{p mn}{n+p m}<\alpha<p m.$

We can extend these arguments to $\phi(\L_{w})$ using Proposition \ref{proposition: cd15}, and to $g_{\L_{w}}$ using Proposition \ref{proposition: cd50}.

\hfill$\Box$

\begin{remark}\label{remark: cd1049}\;  We can easily construct weights satisfying the conditions on $r_{w}$ and $s_{w}$ in Theorem \ref{theorem: cd1048} that are not power weights. Indeed, define $w=u_{1}^{\frac{p m}{2p m+n}}u_{2}^{-1-\frac{2p m}{n}},$ where $u_{1}, u_{2}\in A_{1}.$ It then follows from properties $(ix)$ and $(viii)$ of Proposition \ref{proposition: cd3000} that $w\in A_{1+\frac{p m}{n}}\cap \mbox{RH}_{ \frac{n}{p m}+2}.$

 \end{remark}

As a direct consequence of Theorem \ref{theorem: cd1048}, we obtain the solvability of the Dirichlet problem on $\rdd:=\rz\times [0, \infty):$
\begin{equation}\label{eq: cd1063} \l\{
\begin{aligned}
&\; \partial_{t}^{2}u-\L_{w}u=0\quad\;\mbox{on}\;\rz, \quad\quad\quad\quad\quad\quad\\
&\;u|_{\partial \rdd}=f\;\quad\;\;\; \mbox{on}\;\partial \rdd=\rz.
\end{aligned}\r.
\end{equation}

\begin{theorem}\label{theorem: cd1064}\;  Assume that $w\in A_{2}, \;p\geq 1,$ and $p, r_{w}, s_{w}$ satisfy the conditions in Theorem \ref{theorem: cd1048}. Then, for any $f\in L^{p}(\rz),$ the Dirichlet problem \eqref{eq: cd1063} admits a solution given by $u(x, t):=e^{-t\L_{w}^{1/2}}f(x),$ and $u(\cdot, t)$ converges strongly to $f$ in $L^{p}(\rz)$ as $t\to 0^{+}.$ Morever, the solution $u$ satisfies the uniform bound
 \begin{equation}\label{eq: cd1065} \sup_{t>0}\|t^{k}\partial_{t}^{k}u(\cdot, t)\|_{L^{p}}\lesssim \|f\|_{L^{p}}, \quad \forall\; k\geq 0.
 \end{equation}
 \end{theorem}

{\it Proof.}\quad The function $u(x, t)$ defined above constitutes a formal solution to \eqref{eq: cd1063}, as can be verified through the theory of sectorial operators (see \cite{H, M, TK}). Moreover, for any admissible $p$ as in Theorem \ref{theorem: cd1048}, we can prove that $e^{-t\L_{w}}f\to f$ strongly in $L^{p}(\rz)$ as $t\to 0^{+}.$ This follows from the argument in \cite[Proposition 4.4; Corollary 4.5]{AM}, relying on Proposition \ref{proposition: cd12} and Lemma \ref{lemma: h1}.
Note that the functional calculus for $\L_w$ provides the integral representation \begin{equation}\label{eq: cd1066}e^{-t\L_{w}^{1/2}}=C\int_{0}^{\infty}e^{-\lambda}\lambda^{1/2}e^{-\frac{t^{2}}{4\lambda} \L_{w}}\frac{d\lambda}{\lambda}.\end{equation} Utilizing \eqref{eq: cd1066} and the uniform bound $\sup_{\lambda>0}\|e^{-\lambda \L_{w}}f\|_{L^{p}}\lesssim \|f\|_{L^{p}}$-proved by Corollary \ref{corollary: cd13} and Lemma \ref{lemma: h6}-we deduce that $e^{-t\L_{w}^{1/2}}f\to f$ strongly in $L^{p}(\rz)$ as $t\to 0^{+}.$ 

For such $p,$ both $\partial_{t}^{k}u(\cdot, t)$ and $\L_{w}^{k/2}u(\cdot, t)$ belong to $L^{p}(\rz)$ for each $k\geq 1$ and $t>0$ by \eqref{eq: cd17} and \eqref{eq: cd1066}, and they coincide in $L^{p}(\rz).$ In particular, the case $k=2$ gives $\partial_{t}^{2}u-\L_{w}u=0$ on $\rz.$ Consider the function $\phi_{t}(z):=(tz)^{k}e^{-tz^{1/2}},$ which is bounded and holomorphic on $\Sigma_{\mu}$ for $\mu\in (\V, \pi).$ The estimate \eqref{eq: cd1065} then follows by Theorem \ref{theorem: cd1048}.

\hfill$\Box$

We now turn to the $L^{p}(\rz)$-boundedness of the operators $t^{1/2}\nabla^{m}e^{-t\L_{w}},$ $G_{\L_{w}}$ and $\nabla^{m} \L_{w}^{-1/2}.$


\begin{theorem}\label{theorem: cd1054}\; Given $w\in A_{2}$ and $p\geq 1.$ Then $t^{1/2}\nabla^{m}e^{-t\L_{w}}: L^{p}(\rz)\to L^{p}(\rz)$ is uniformly bounded for all $t>0,$ under the following conditions on $p:$ 
\begin{equation}\label{eq: cd1056} |p-2|<\ez \;\; \mbox{with}\;\; 0<\ez<\min\{\frac{4m}{nr_{w}+2m}, q_{+}-2\}, \end{equation} 
and \begin{equation}\label{eq: cd1055}1\leq r_{w}<\frac{q_{+}}{p},\quad  s_{w}>\frac{p}{p -\frac{2nr_{w}}{nr_{w}+2m} }.\end{equation} (Note that $q_{+}=q_{+}(\L_{w})>2$ for any $w\in A_2$, as established in Section 6.3.) Moreover, in the same range of $p$, the operators $\nabla^m \L_w^{-1/2}$ and $G_{\L_w}$ are also bounded on $L^{p}(\rz).$

The conditions \eqref{eq: cd1056}-\eqref{eq: cd1055} are satisfied in any of the following scenarios: $(i)$ $w\in A_{1}\cap \mbox{RH}_{\frac{p}{p -\frac{2n}{n+2m} } }$ and $|p-2|<\ez$ with $0<\ez<\min\{ \frac{4m}{n+2m}, q_{+}-2\}$; $(ii)$ 
Given $\Theta \geq 1$, there exists $\ez_{0}=\ez_{0}(\Theta, c_{1}, c_{2}, n, m)$ such that $0<\ez_{0}\leq \frac{1}{2n},$ $[w]_{A_{2}}\leq \Theta$ and $w\in A_{1+\ez_{1}}\cap \mbox{RH}_{\frac{p}{p -\frac{2n(1+\ez_{1})}{n(1+\ez_{1})+2m} } }$ for some $0<\ez_{1}<\frac{\ez_{0}}{2},$ and the exponent $p$ satisfies $|p-2|<\ez$ with $0<\ez <\ez_{2},$ where 
\begin{equation*}\ez_{2}:=\l\{
\begin{aligned}
&\; \min\{\frac{4m}{n+2m}, q_{+}-2,\; \frac{1}{4},\; \frac{2m}{n-m}, \;\frac{\ez_{0}-2\ez_{1}}{1+\ez_{1}}\}\quad\;\;\mbox{if}\; m<n, \quad\quad\quad\quad\quad\quad\\
&\;\min\{\frac{4m}{n(1+\ez_{1})+2m},\; q_{+}-2,\; \frac{1}{4},\; \frac{\ez_{0}-2\ez_{1}}{1+\ez_{1}}\}\quad\;\;\; \mbox{if}\;m\geq n.
\end{aligned}\r.
\end{equation*}

In particular, for the power weight $w_{\alpha}:=|x|^{\alpha},$ there exists a $\ez_{3},$ depending only on $n,m, c_{1}, c_{2},$ $0<\ez_{3}<\frac{1}{2n},$ such that if $|p-2|<\ez_{4}$ and $-\frac{n(p-\frac{2n}{n+2m})}{p}<\alpha<\ez_{3},$ 
with $\ez_{4}>0$ given by  
\begin{equation*}\ez_{4}:=\l\{
\begin{aligned}
&\; \min\{\frac{4m}{n+2m},\; q_{+}-2, \;\frac{1}{4}, \;\frac{2m}{n-m}, \;\frac{2\ez_{3}}{1+\ez_{3}}\}\quad\;\;\mbox{if}\; m<n, \quad\quad\quad\quad\quad\quad\\
&\;\min\{\frac{4m}{n(1+\ez_{3})+2m}, \;q_{+}-2, \;\frac{1}{4}, \;\frac{2\ez_{3}}{1+\ez_{3}}\}\quad\;\;\; \mbox{if}\;m\geq n,
\end{aligned}\r.\end{equation*} then the $L^p(\rz)$-boundedness of the aforementioned operators holds for $w_\alpha$,
\end{theorem} 
{\it Proof.}\quad We prove the theorem for $t^{1/2}\nabla^m e^{-t\L_w}$ by using Proposition \ref{proposition: cd100} (Corollary \ref{corollary: cd1010}). The proofs for $\nabla^m \L_w^{-1/2}$ and $G_{\L_w}$ follow similarly, by replacing Proposition \ref{proposition: cd100} with Proposition \ref{proposition: cd1023} and Proposition \ref{proposition: cd1040}, respectively.


Corollary \ref{corollary: cd1010} shows that to prove $t^{1/2}\nabla^{m}e^{-t\L_{w}}: L^{p}(\rz)\to L^{p}(\rz),$ it suffices to verify that $$w\in A_{\frac{q_{+}}{p}}\cap \mbox{RH}_{(\frac{p}{q_{-}})'}.$$ So we must ensure that $r_{w}<\frac{q_{+}}{p}$ holds. Furthermore, if we can show $s_{w}>\frac{p}{p -\frac{2nr_{w}}{nr_{w}+2m} },$ then $w\in \mbox{RH}_{(\frac{p}{q_{-}})'}$ follows from Proposition \ref{proposition: cd12}, since $q_{-}=p_{-}\leq (2_{w}^{*, m})'.$ 

In case $(i),$ we clearly have $r_{w}=1$ and $s_{w}>\frac{p}{p -\frac{2n}{n+2m} }.$ Hence \eqref{eq: cd1055} holds because $|p-2|<q_{+}-2.$  The proof for case (ii) is more involved. Since $w\in A_{1+\ez_{1}}\cap\mbox{RH}_{\frac{p}{p -\frac{2n(1+\ez_{1})}{n(1+\ez_{1})+2m} } }$ and $|p-2|<\ez_{2},$ we readily obtain $s_{w}>\frac{p}{p -\frac{2n(1+\ez_{1})}{n(1+\ez_{1})+2m} }>\frac{p}{p -\frac{2nr_{w}}{nr_{w}+2m} }.$ Thus, the proof reduces to showing $p r_{w} <q_{+}.$ In view of the inequality $q_{+}>q_{w}$ from Section 6.3, it is enough to show $$p r_{w}<q_{w}:=\min\{p_{+}(\L_{w}), r_{w}', p_{0}\}.$$ Here $p_{0}$ is given by \eqref{eq: cd1022}. To the end, we need to determine a suitable threshold $\ez_{0}.$  

Observe that $r_{w}<1+\ez_{1}<1+\ez_{0}<1+\frac{1}{p n}< 1+\frac{1}{p },$ which implies that $p r_{w}<r_{w}'.$ On the other hand, by Proposition \ref{proposition: cd12} and $|\eta-2|<\ez_{2},$ we also have $p r_{w}<2_{w}^{*, m}\leq p_{+}.$ It therefore remains to prove that $p r_{w}<p_{0}.$ Before proceeding, we point out that the definition of $p_0$ in \eqref{eq: cd1022} is inadequate, primarily because of the logarithmic term in Theorem \ref{theorem: cd8}. To address this, we refine the definition of $p_0$. This is possible under the assumption that $[w]_{A_{2}}\leq \Theta.$ 

By the bound $[w]_{A_2} \leq \Theta$ and property $(iii)$ of Proposition \ref{proposition: cd3000}, there exist positive constants $C_{0}=C_{0}(n, \Theta)$ and $\delta=\delta(n, \Theta)$ (small) such that $$[w]_{A_{2-\delta}}\leq C_{0},$$ see \cite{G1}.
With $q_{0}:=2-\frac{1}{Nn}$ and $N$ (depending only on $m, n, \Theta,$) sufficiently large, we have $$[w]_{A_{q_{0}}}\leq [w]_{A_{2-\delta}}\leq C_{0},$$ provided that $N>\frac{1}{\delta n}.$ This yields $$\frac{2n(q_{0}+\log[w]_{A_{q_{0}}}) }{n(q_{0}+\log[w]_{A_{q_{0}}})+2m}<q_{0}$$ if $N>\frac{n+2m+n\log C_{0}}{2mn}.$ Hence, for any $N>\max\{\frac{1}{\delta n}, \frac{n+2m+n\log C_{0}}{2mn}\},$ we get $$\max\{r_{w}, \frac{2n(q_{0}+\log[w]_{A_{q_{0}}}) }{n(q_{0}+\log[w]_{A_{q_{0}}})+2m}\}<q_{0}<2\leq n.$$


Recall that $\frac{1}{(q_{0})^{*, m}_{w}}:=\frac{1}{q_{0}}-\frac{m}{n(q_{0}+\log[w]_{A_{q_{0}}}) }$ if $q_{0}<\frac{n(q_{0}+\log[w]_{A_{q_{0}}})}{m},$ and $(q_{0})^{*, m}_{w}=\infty$ otherwise. Obviously, $2<(q_{0})^{*, m}_{w}.$ Invoking Theorem \ref{theorem: cd8} and repeating the argument leading to \eqref{eq: cd1022}, we obtain  
\begin{equation*}p_{0}:=2+\frac{2-q_{0}}{2^{\frac{4}{q_{0}} +1} C_{1}^{2} C_{2}^{2}(2^{D}[w]_{A_{2}})^{\frac{6}{q_{0}} +17} }\approx 2+\frac{1 } {C N n [w]^{\frac{2}{q_{0}}}_{A_{q_{0}}} [w]_{A_{2}}^{m+\frac{6}{q_{0}} +17} }.\end{equation*} Here, $D=D(n),$ $C=C(n,m, c_{1}, c_{2}, N)$ and $C_{1}, C_{2}$ are as defined in \eqref{eq: cd1018} and \eqref{eq: cd1020}, respectively. Since $ [w]_{A_{2}}\leq  \Theta$ and $[w]_{A_{q_{0}}}\leq C_{0},$ it holds that 
$$p_{0}\geq 2+\frac{1} {N n CC_{0}^{\frac{2}{q_{0}}} \Theta^{m+\frac{6}{q_{0}} +17 }}=2+2\ez_{0},$$ where $\ez_{0}:=(2NnC)^{-1}$ depends only on $n,m, c_{1}, c_{2}, \Theta.$ Clearly, $0<\ez_{0}<\frac{1}{2  n }$ and $p r_{w}<p_{0},$ as $|p-2|<\frac{\ez_{0}-2\ez_{1}}{1+\ez_{1}}$ and $\ez_{1}<\frac{\ez_{0}}{2}.$ This proves case $(ii).$

We now consider the power weight $w_{\alpha}(x):=|x|^{\alpha}.$ If $|p-2|<\ez$ with $\ez<\min\{ \frac{4m}{n+2m}, q_{+}-2\}$ and $-\frac{n(p-\frac{2n}{n+2m})}{p}<\alpha\leq 0,$ then $r_{w_{\alpha}}=1$ and $s_{w_{\alpha}}=-\frac{n}{\alpha},$ so condition $(i)$ is satisfied. This yields the desired estimates.

If $0<\alpha<\frac{1}{2},$ then $r_{w_{\alpha}}=1+\frac{\alpha}{n}<1+\frac{1}{2n}$ and $s_{w_{\alpha}}=\infty,$ so $w_{\alpha}\in A_{2}.$ It is well-known that $$[|x|^{\alpha}]_{A_{p}}=\frac{n}{n+\alpha}\l(\frac{n(p-1)}{n(p-1)-\alpha}\r)^{p-1}, \quad 1<p<\infty, \quad -n < \alpha < n(p-1).$$ Consequently, $$[w_{\alpha}]_{A_{2}}\leq 2:=\Theta, \quad \forall \;0<\alpha<\frac{1}{2}.$$ Applying the preceding argument, we can find a constant $\ez_{0},$ depending only on $n,m, c_{1}, c_{2},$ such that $0<\ez_{0}<\frac{1}{2n}$ and $p_{0}\geq 2+2\ez_{0}.$ Define $\ez_{3}:=\frac{\ez_{0}}{4}.$ Then, for $0<\alpha<\ez_{3}<\frac{1}{2n},$ we have $w_{\alpha}\in A_{1+\ez_{3}}.$ Moreover, for such  $\alpha$ and for any $p$ with $|p-2| <\ez_{4}$ (where $\ez_{4}$ is defined as above), 
we see $w_{\alpha}\in \mbox{RH}_{\frac{p}{p -\frac{2n(1+\ez_{3})}{n(1+\ez_{3})+2m} } }.$ Thus, condition $(ii)$ is satisfied, which leads to the desired estimates as well.


\hfill$\Box$

\begin{remark}\label{remark: cd1057}\;  If $u\in A_{1}$ and $p>\frac{2n}{n+2m},$ then $w:=u^{\frac{p-\frac{2n}{n+2m}}{p}}\in A_{1}\cap \mbox{RH}_{\frac{p}{p -\frac{2n}{n+2m}} }.$ Moreover, if $|p-2|<\ez$ with $\ez<\min\{ \frac{4m}{n+2m}, q_{+}-2\},$ the weight $w$ satisfies condition $(i)$ in Theorem \ref{theorem: cd1054}. Clearly, $w$ is not a power weight. Besides, given $u\in A_{2}$ and $0<\theta<1,$ let $w:=u^{\theta}.$ Then, property $(vii)$ in Proposition \ref{proposition: cd3000} implies $w\in A_{1+\theta}.$ Furthermore, there exists a $\gamma,$ depending only on $n, [u]_{A_{2}},$ such that $u\in \mbox{RH}_{1+\gamma },$ or equivalently,  $u^{-1}\in A_{(1+\gamma)'}(udx);$ see property $(x)$ in Proposition \ref{proposition: cd3000}. From this, applying property $(vii)$ in Proposition \ref{proposition: cd3000} again yields $u^{-\theta}\in A_{\theta (1+\gamma)'+1-\theta}(udx),$ and hence $w\in \mbox{RH}_{(\theta (1+\gamma)'+1-\theta)' }.$ Note that $(\theta (1+\gamma)'+1-\theta)'\to \infty$ as $\theta\to 0^{+}.$ Thus, by repeating the argument leading to $p_{0}$ in Theorem \ref{theorem: cd1054} and choosing $\theta$ sufficiently small (depending on $n, m, c_{1}, c_{2}, [u]_{A_{2}}$), the weight $w$ satisfies condition $(ii).$

\end{remark} 

Combining Theorem \ref{theorem: cd1048} and Theorem \ref{theorem: cd1054}, we obtain the solvability of the Neumann problem 
\begin{equation}\label{eq: cd1070} \l\{
\begin{aligned}
&\; \partial_{t}^{2}u-\L_{w}u=0\quad\;\mbox{on}\;\rz, \quad\quad\quad\quad\quad\quad\\
&\;\partial_{t}u|_{\partial \rdd}=f\;\quad\;\;\; \mbox{on}\;\partial \rdd=\rz, .
\end{aligned}\r.
\end{equation}

\begin{theorem}\label{theorem: cd1071}\; Given $w\in A_{2}$ and $p \geq 1.$ Suppose that $p, r_{w}, s_{w}$ satisfy the conditions in Theorem \ref{theorem: cd1054}. Then for any $f\in L^{p}(\rz),$ $u(x, t):=-\L_{w}^{-1/2}e^{-t\L_{w}^{1/2}}f(x)$ solves the Neumann problem \eqref{eq: cd1070} with $\partial_{t}u(\cdot, t)\to f$ strongly in $L^{p}(\rz)$ as $t\to 0^{+},$ and satisfies for all $k\geq 1:$ 
\begin{equation}\label{eq: cd1072}\sup_{t>0}\l(\|t^{k-1}\partial_{t}^{k}u(\cdot, t)\|_{L^{\eta}}+\|\nabla^{m}u(\cdot, t)\|_{L^{\eta}}\r) \lesssim \|f\|_{L^{\eta}}.\end{equation}
\end{theorem} 
{\it Proof.}\quad Following the argument in Theorem \ref{theorem: cd1064}, the function $u(x, t):=-\L_{w}^{-1/2}e^{-t\L_{w}^{1/2}}f(x)$ defines a formal solution to \eqref{eq: cd1070}, with $\partial_{t}u(\cdot, t)\to f$ strongly in $L^{p}(\rz)$ as $t\to 0^{+}.$ Then, it follows from Theorem \ref{theorem: cd1054} and Theorem \ref{theorem: cd1048} that \begin{equation*}\begin{aligned}\|\nabla^{m}u(\cdot, t)\|_{L^{p}}\lesssim \|\nabla^{m} \L_{w}^{-1/2}e^{-t\L_{w}^{1/2}}f\|_{L^{p}}\lesssim\|e^{-t\L_{w}^{1/2}}f\|_{L^{p}}\lesssim \|f\|_{L^{p}}
\end{aligned}\end{equation*} while  $$\|t^{k-1}\partial_{t}^{k}u(\cdot, t)\|_{L^{p}}\lesssim \|f\|_{L^{p}},\;\forall\; k\geq 1,$$ follows from \eqref{eq: cd1065}. The proof is complete. 

\hfill$\Box$ 

The following theorem establishes unweighted $L^{p}$ reverse inequalities for the square root of $\L_{w}.$

\begin{theorem}\label{theorem: cd1050}\;  Given $w\in A_{2},$ let $p \geq 1$ such that $|p-2|<\ez$ with $0<\ez<\min\{\frac{2m}{n+2m}, 2-r_{w}\}.$ Assume that $1\leq  r_{w}< 1+\frac{p m}{n}$ and $s_{w}>\max\{(\frac{p}{r_{w}})', \frac{nr_{w}}{p m}+1\}.$ Then \begin{equation}\label{eq: cd1051}\|\L_{w}^{1/2} f\|_{L^{p}(\rz)}\leq C \| \nabla^{m}f\|_{L^{p}(\rz)},\quad \forall \;f\in \ss(\rz).\end{equation} In particular, \eqref{eq: cd1051} holds for any $p$ satisfying $|p-2|<\ez$ with $0<\ez<\min\{\frac{2m}{n+2m}, 2-r_{w}\},$ provided one of the following conditions is met: $(i)$ $w\in A_{1}\cap \mbox{RH}_{\max\{p', \frac{n }{p m}+1 \}};$ $(ii)$ $w\in A_{r}\cap \mbox{RH}_{\max\{ (\frac{p}{r})', \frac{n r}{p m}+1 \}}$ for some $1<r<\min\{ p, 1+\frac{p m}{n}\};$ $(iii)$ $w=w_{\alpha}(x):=|x|^{\alpha}$ with $\max\{-\frac{n}{p}, -\frac{p mn}{n+p m}\}<\alpha<p m.$
 \end{theorem}
{\it Proof.}\quad By \eqref{eq: cd83} in Proposition \ref{proposition: cd81}, if $r_{w}\leq p_{-},$ the proof is identical to that of Theorem \ref{theorem: cd1048}; otherwise we proceed as in Theorem \ref{theorem: cd1048} with the choices $p_{0}=r_{w}$ and $q_{0}=p^{*, m}_{w}.$

\hfill$\Box$ 

\begin{remark}\label{remark: cd1052}\; It is clear that $\max\{ (\frac{p}{r})', \frac{n r}{p m}+1 \}=\frac{n r}{p m}+1$ holds if $r\leq p (1-\frac{m}{n}).$ Moreover, this condition is guaranteed when $n\geq \frac{2p m}{p-1},$ since then $1+\frac{p m}{n} \leq p (1-\frac{m}{n}).$ In this case, the conditions in the second part of Theorem \ref{theorem: cd1050} simplify to those of Theorem \ref{theorem: cd1048}.\end{remark} 

\begin{remark}\label{remark: cd1053}\;  Tracking carefully the proofs of Theorem \ref{theorem: cd1048} and Theorem \ref{theorem: cd1050}, we find that the condition $1\leq  r_{w}<1+\frac{p m}{n}$ may be relaxed to the potentially weaker condition $1\leq  r_{w}<\frac{p_{+}}{p}$ by taking $q_{0}=p_{+}$ in the argument. 
\end{remark}


Synthesizing the results of Theorem \ref{theorem: cd1054}, Theorem \ref{theorem: cd1050} and Remark \ref{remark: cd1053} we conclude with the following unweighted Kato estimate for higher-order degenerate elliptic operators:  
\begin{theorem}\label{theorem: cd1058}\; Let $\L_{w}$ be as in \eqref{eq: z00}-\eqref{eq: z2} with $w\in A_{2}.$ If there exists a $\ez>0$ small enough such that \begin{equation}\label{eq: cd1059} |p-2|<\ez, \; 1\leq r_{w}<\frac{q_{+}(\L_{w})}{p} \; \mbox{and}\; s_{w}>\max\{(\frac{p}{r_{w}})', \frac{nr_{w}+p m}{p m}, \frac{p}{p -\frac{2nr_{w}}{nr_{w}+2m} }\},\end{equation} then, for every $f\in H^{m}(\rz),$ we have the Kato estimate \begin{equation}\label{eq: cd1060} \|L_{w}^{1/2} f\|_{L^{p}(\rz)}\approx \| \nabla^{m}f\|_{L^{p}(\rz)},\end{equation} where the implicit constants depend only on $n, m, c_{1}, c_{2}$ and $ [w]_{A_{2}}$. 

In particular, \eqref{eq: cd1060} holds for any $|p-2|<\ez$ with $\ez$ sufficiently small, in each of the following scenarios: $(i)$ $w\in A_{1}\cap \mbox{RH}_{\max\{p', \frac{n }{p m}+1, \frac{p}{p -\frac{2n}{n+2m} }  \}};$ $(ii)$ 
Given $\Theta \geq 1$, there exist $\ez_{0}=\ez_{0}(\Theta, c_{1}, c_{2}, n, m),$ $0<\ez_{0}\leq \frac{1}{2n},$ such that $w\in A_{1+\ez_{1}}\cap \mbox{RH}_{\max\{ (\frac{p}{(1+\ez_{1})})', \frac{n (1+\ez_{1})}{p m}+1,  \frac{p}{p -\frac{2n(1+\ez_{1})}{n(1+\ez_{1})+2m} }  \}},$ $0<\ez_{1}<\frac{\ez_{0}}{2}$ and $[w]_{A_{2}}\leq \Theta.$

Finally, there exists a $\ez_{2}=\ez_{2}(n, m, c_{1}, c_{2})\in (0,\; \frac{1}{2n})$ such that for $p$ near 2, \eqref{eq: cd1060} holds for $w_{\alpha}=|x|^{\alpha}$ whenever the exponent $\alpha$ satisfies $$\max\{-\frac{n}{p}, -\frac{p mn}{n+p m}, -\frac{n(p-\frac{2n}{n+2m})}{p}\}<\alpha<\ez_{2}.$$
\end{theorem}


\begin{remark}\label{remark: cd1062}\; In particular, for the power weight $w_{-\gamma}:=|x|^{-\gamma}$ with $-\ez_{2}<\gamma< \frac{2 mn}{n+2 m},$ Theorem \ref{theorem: cd1058} gives $$\|\L_{w_{-\gamma}}^{1/2} f\|_{L^{2}(\rz)}\approx \| \nabla^{m}f\|_{L^{2}(\rz)},$$ where $\L_{w_{-\gamma}}$ is defined by \eqref{eq: z00}-\eqref{eq: z2}.  When $\gamma=0,$ we recover the classical Kato square root problem for higher-order elliptic operators, which was settled in \cite{AHMT1}.
\end{remark} 
 
To conclude this section, we address the solvability of the regularity problem \eqref{eq: cd2004} on $\rdd,$ using Theorem \ref{theorem: cd1058}. 
\begin{theorem}\label{theorem: cd1067}\; Let $w\in A_{2},$ with $p, r_{w}, s_{w}$ satisfying the requirements of Theorem \ref{theorem: cd1058}. Then, for any $f\in H^{m, p}(\rz),$ $u(x, t)$ (defined as in Theorem \ref{theorem: cd1064}) is a solution to the regularity problem \eqref{eq: cd2004} with $\nabla^{l}u(\cdot, t)\to \nabla^{l}f$ strongly in $L^{p}(\rz)$ as $t\to 0^{+}$ for all $0\leq l\leq m-1$. Furthermore, for all $k\geq 1$ and $0\leq l\leq m,$
\begin{equation}\label{eq: cd1068}\sup_{t>0}\l(\|t^{k-1}\partial_{t}^{k}u(\cdot, t)\|_{L^{p}}+\|\nabla^{l}u(\cdot, t)\|_{L^{p}}\r) \lesssim \|f\|_{H^{m, p}}.\end{equation}
\end{theorem} 
{\it Proof.}\quad As established in the first part of Theorem \ref{theorem: cd1064}, we see that $e^{-t\L_{w}^{1/2}}f\to f$ strongly in $L^{p}(\rz)$ as $t\to 0^{+}$. 
From this, along with Theorem \ref{theorem: cd1058}, it follows that 
\begin{equation*}\begin{aligned}\|\nabla^{m}u(\cdot, t)\|_{L^{p}} \lesssim \|\nabla^{m} \L_{w}^{-1/2}(\L_{w}^{1/2}e^{-t\L_{w}^{1/2}}f)\|_{L^{p}}&\lesssim\|\L_{w}^{1/2}e^{-t\L_{w}^{1/2}}f\|_{L^{p}}\\[4pt] 
&\lesssim \|e^{-t\L_{w}^{1/2}}\L_{w}^{1/2} f\|_{L^{p}}\lesssim \|\L_{w}^{1/2} f\|_{L^{p}}\lesssim \|\nabla^{m}f\|_{L^{p}}.
\end{aligned}\end{equation*} 
Similarly, for all $k\geq 1,$ \begin{equation*}\begin{aligned}\|t^{k-1}\partial_{t}^{k}u(\cdot, t)\|_{L^{p}} \lesssim \| (t\L_{w}^{1/2})^{k-1}\L_{w}^{1/2}e^{-t\L_{w}^{1/2}}f\|_{L^{p}}&\lesssim\|(t\L_{w}^{1/2})^{k-1}e^{-t\L_{w}^{1/2}}(\L_{w}^{1/2}f)\|_{L^{p}}\\[4pt] 
&\lesssim \|\L_{w}^{1/2} f\|_{L^{p}}\lesssim \|\nabla^{m}f\|_{L^{p}}.
\end{aligned}\end{equation*} 

Recall from Remark \ref{remark: cd1053} that \eqref{eq: cd1065} also holds under the hypotheses of Theorem \ref{theorem: cd1067}.
Thus, \begin{equation}\label{eq: cd1069}\sup_{t>0}\|t^{k}\nabla^{m} \partial^{k-1}_{t} u(\cdot, t)\|_{L^{p}}\lesssim \|f\|_{L^{p}},\quad \forall\; k\geq 1.\end{equation} Furthermore, interpolation gives for all $0\leq l\leq m:$ $$\sup_{t>0}\|\nabla^{l} u(\cdot, t)\|_{L^{p}}\lesssim \sup_{t>0}\|\nabla^{m}u(\cdot, t)\|_{L^{p}}^{\frac{l}{m}}\cdot \sup_{t>0}\|u(\cdot, t)\|_{L^{p}}^{1-\frac{l}{m}}\lesssim  \|f\|_{L^{p}}+\|\nabla^{m}f\|_{L^{p}}.$$ We therefore obtain that for all $0\leq l\leq m-1,$ $\nabla^{l}u(\cdot, t)$ converges to $\nabla^{l}f$ strongly in $L^{p}(\rz)$ as $t\to 0^{+}$ when $f\in H^{m, p}(\rz).$

\hfill$\Box$

\begin{remark}\label{remark: cd8000}\; When $m=1,$ Theorem \ref{theorem: cd1054}, Theorem \ref{theorem: cd1071} and Theorem \ref{theorem: cd1067} reduce to \cite[Theorem 12.2]{DMR}, \cite[Theorem 12.10]{DMR} and  \cite[Theorem 12.6]{DMR}, respectively. 
\end{remark} 

\section{Appendix}
The first two lemmas, while auxiliary, are crucial to the core argument. The first of these generalizes \cite[Lemma 3.3]{CR2}. 


\begin{lemma}\label{lemma: h2}\; Assume that $\{a_{\alpha, \beta}\}_{|\alpha|=|\beta|=m}\in \Ez(w, c_{1}, c_{2}).$ Then, $\{z a_{\alpha, \beta}\}_{|\alpha|=|\beta|=m}\in \Ez(w, \lambda_{z}, \Lambda_{z})$ for any $z\in \Sigma_{\frac{\pi}{2}-\V},$ where $\V$ is given by \eqref{eq: cd9001}.
 \end{lemma}
{\it Proof.}\quad Fix $f\in \D(\L_{w}),$ and define $$\S:=\sum_{|\alpha|=|\beta|=m} \int_{\rz}a_{\alpha, \beta}( x)\partial^{\alpha}f(x)  \overline{\partial^{\beta}f(x)}dx= <\L_{w} f, f>.$$ Its imaginary and real parts are denoted by $\R:=\mbox{Im}\;\S$ and $\T:=\mbox{Re}\;\S,$ respectively. Using the definition of $\V$ and \eqref{eq: z1}, we have \begin{equation*}\begin{aligned}\mbox{Re}\; (z \S) &=|z|(\cos (\mbox{arg} z)\T-\sin (\mbox{arg} z)\R)\\[4pt] 
&=|z|\T (\cos (\mbox{arg} z)-\sin (\mbox{arg} z)\frac{\R}{\T})\\[4pt] 
&\geq c_{1}|z| \|\nabla^{m} u\|_{L^{2}(w)}^{2} (\cos (\mbox{arg} z)-|\sin (\mbox{arg} z)|\tan \V).
\end{aligned}\end{equation*} 
Consequently, since $|\sin (\mbox{arg} z)|<\frac{\cos (\mbox{arg} z)}{\tan \V}$ and $\D(\L_{w})$ is dense in $H^{m}(w),$ the identity $\lambda_{z}=c_{1}|z| (\cos (\mbox{arg} z)-|\sin (\mbox{arg} z) |\tan \V)$ is valid. On the other hand, \eqref{eq: z2} implies $$\bigg|\sum_{|\alpha|=|\beta|=m}z a_{\alpha, \beta}(x)\xi_{\alpha} \overline{\zeta_{\beta}}\bigg|\leq c_{2}|z| |\xi||\zeta|w(x),$$ which immediately gives $\Lambda_{z}=c_{2}|z|.$

\hfill$\Box$

\begin{lemma}\label{lemma: z1}\; Let $s>0$, $\alpha\geq 0,$ and $\beta>0$ with $\alpha\neq \beta.$ Then, for any $0<c'<c,$ $$\sum_{k=0}^{\infty}2^{k\alpha}\Upsilon(2^{k}s)^{\beta} e^{-c s^{\frac{2m}{2m-1}}} \lesssim \Upsilon(s)^{\max\{\alpha, \beta\}} e^{-c's^{\frac{2m}{2m-1}}}.$$ 
 \end{lemma}
{\it Proof.}\quad In light of \cite[Lemma 6.3]{AM}, the proof is routine, and we skip it.

\hfill$\Box$

We now present a detailed proof of the reverse H\"{o}lder inequality with sharp constants for solutions to $\L_{w},$ a result referenced in Section 6.3.

\begin{lemma}\label{lemma: cd1011}\; Fix $B_{0}:=B(x_{0}, R),$ and suppose $w\in A_{2}.$ Consider any solution $u\in H^{m}(B_{0}, w)$ to $\L_{w} u=0$ in $B_{0}.$ Then, for any $0<r<R,$ we have \begin{equation}\label{eq: cd1012}
\int_{B(x_{0}, r)}|\nabla^{m} u|^{2}dw \leq \sum_{k=0}^{m-1} \frac{C}{(R-r)^{2m-2k}}\int_{B_{0}\setminus B(x_{0}, r)}|\nabla^{k}(u-P_{B_{0}}(u))|^{2}dw,
\end{equation} where the constant $C$ depends only on $c_{1}, c_{2}, m, n.$
  \end{lemma}
{\it Proof.}\quad Let $\phi$ be a smooth, nonnegative, real-valued test function supported in $B_{0},$ identically 1 on $B(x_{0}, r),$ and satisfying $|\nabla^{k}\phi|\leq C_{k} (R-r)^{-k}$ for any $0\leq k\leq m.$ Testing the equation $\L_{w} u=0$ in $B_{0}$ against the function $\psi:=\phi^{4m}\widetilde{u},$ where $\widetilde{u}= (u-P_{B_{0}}(u)),$ yields $$\sum_{|\alpha|=|\beta|=m}\int_{B_{0}}a_{\alpha, \beta}(x)\partial^{\alpha}\widetilde{u}(x) \cdot \overline{\partial^{\beta}\psi(x)}dx=0.$$ From this, along with the product rule, it holds that $$-\int_{B_{0}}a_{\alpha, \beta}(x)\partial^{\alpha}\widetilde{u}(x) \sum_{\gamma<\beta}C_{\beta}^{\gamma}\partial^{\beta-\gamma}\phi^{2m}\partial^{\gamma} (\phi^{2m}\overline{\widetilde{u}(x)} )dx=\int_{B_{0}}a_{\alpha, \beta}(x)\phi^{2m} \partial^{\alpha}\widetilde{u}(x)\partial^{\beta} (\phi^{2m}\overline{\widetilde{u}(x)})dx.$$ Note from \cite[Lemma 3.8]{B1} that there exists functions $\Phi_{\beta, \xi}$ supported in $B_{0}\setminus B(x_{0}, r)$ with $|\Phi_{\beta, \xi}|\leq C(R-r)^{|\xi|-|\beta|}$ such that 
we may write $$\sum_{\gamma<\beta}C_{\beta}^{\gamma}\partial^{\beta-\gamma}\phi^{2m}\partial^{\gamma} (\phi^{2m}\overline{\widetilde{u}(x)})=\sum_{\xi <\beta} \phi^{2m} \Phi_{\beta, \xi}\partial^{\xi}\overline{\widetilde{u}(x)}.$$ Thus
\begin{equation*}\begin{aligned} \int_{B_{0}}a_{\alpha, \beta}(x) \partial^{\alpha}(\phi^{2m}\widetilde{u}(x))\partial^{\beta} (\phi^{2m}\overline{\tilde{u}(x)})dx&=\int_{B_{0}}a_{\alpha, \beta}(x)\sum_{\gamma<\alpha} C_{\alpha}^{\gamma}\partial^{\alpha-\gamma}\phi^{2m}\partial^{\gamma} \widetilde{u}(x) \partial^{\beta} (\phi^{2m}\overline{\tilde{u}(x)}) dx\\[4pt] 
&\quad\quad -\int_{B_{0}}a_{\alpha, \beta}(x)\partial^{\alpha}\widetilde{u}(x) \sum_{\xi <\beta} \phi^{2m} \Phi_{\beta, \xi}\partial^{\xi}\overline{\widetilde{u}(x)}dx\\[4pt]
&=\int_{B_{0}}a_{\alpha, \beta}(x)\sum_{\gamma<\alpha} C_{\alpha}^{\gamma}\partial^{\alpha-\gamma}\phi^{2m}\partial^{\gamma} \widetilde{u}(x) \partial^{\beta} (\phi^{2m}\overline{\widetilde{u}(x)}) dx\\[4pt]
&\quad\quad-\int_{B_{0}}a_{\alpha, \beta}(x)\partial^{\alpha}(\widetilde{u}(x)\phi^{2m})\sum_{\xi <\beta} \Phi_{\beta, \xi}\partial^{\xi}\overline{\widetilde{u}(x)} dx\\[4pt]
&\quad\quad- \sum_{\gamma<\alpha} C_{\alpha}^{\gamma}(\partial^{\alpha-\gamma}\phi^{2m}\partial^{\gamma} \widetilde{u}(x) )\sum_{\xi <\beta} \Phi_{\beta, \xi}\partial^{\xi}\overline{\widetilde{u}(x)} dx
\end{aligned}\end{equation*} Invoking \eqref{eq: z1}-\eqref{eq: z2} and applying Young's inequality, we get \begin{equation*}\begin{aligned} c_{1}\int_{B_{0}} |\nabla^{m} (\phi^{2m}\widetilde{u}(x) )|^{2}dw&\leq \frac{c_{1}}{2} \int_{B_{0}} |\nabla^{m} (\phi^{2m}\widetilde{u}(x) )|^{2}dw\\[4pt]  &+\sum_{k=0}^{m-1} \frac{C}{(R-r)^{2m-2k}}\int_{B_{0}\setminus B(x_{0}, r)}|\nabla^{k}(\widetilde{u}(x))|^{2}dw,\end{aligned}\end{equation*} where $C$ depends only on $c_{1}, m , n,  c_{2}.$ This yields \eqref{eq: cd1012}. 

\hfill$\Box$

The bound on the right-hand side of \eqref{eq: cd1012} can be improved to depend solely on $\|u\|_{L^{2}(w)}.$To achieve this, we adapt the approach from \cite[Theorem 3.10]{B1}.

\begin{corollary}\label{corollary: cd1013}\;  Let $B_{0}:=B(x_{0}, R)$ with $x_{0}\in \rz$ and $R>0.$ Given $w\in A_{2},$ assume that $\tilde{u}\in H^{m}(B_{0}, w)$ satisfies for any $0<\rho<r<R,$ \begin{equation}\label{eq: cd1014}
\int_{B(x_{0}, \rho)}|\nabla^{m} \widetilde{u}|^{2}dw \leq \sum_{k=0}^{m-1} \frac{C}{(r-\rho)^{2m-2k}}\int_{B(x_{0}, r)\setminus B(x_{0}, \rho) } |\nabla^{k}\widetilde{u}|^{2}dw.
\end{equation} Then $\widetilde{u}$ satisfies the following improved estimates: \begin{equation}\label{eq: cd1016}
\int_{B(x_{0}, r)}|\nabla^{m} \widetilde{u}|^{2}dw \leq  \frac{C[w]^{m}_{A_{2}}}{(R-r)^{2m}}\int_{B(x_{0}, R) \setminus B(x_{0}, r) } |\widetilde{u}|^{2}dw
\end{equation} and, for any $0\leq j\leq m-1,$ 
\begin{equation}\label{eq: cd1015}
\int_{B(x_{0}, r)}|\nabla^{j} \widetilde{u}|^{2}dw \leq \frac{C[w]^{j(m+1-j)}_{A_{2}}}{(R-r)^{2j}}\int_{B(x_{0}, R) } |\widetilde{u}|^{2}dw.
\end{equation}  Here, the constant $C$ depends only on $c_{1}, m , n,  c_{2}.$\end{corollary}
{\it Proof.}\quad Let $A(r, \xi)$ (with $\xi>0$) denote the annulus $B(x_{0}, r+\xi)\setminus B(x_{0}, r-\xi)$ for the proof of \eqref{eq: cd1016}, and the ball $B(x_{0}, r+\xi)$ for that of \eqref{eq: cd1015}, respectively. 

To prove \eqref{eq: cd1016}-\eqref{eq: cd1015}, it suffices to establish the estimate \begin{equation}\label{eq: cd1017}\int_{A(r, \xi) }|\nabla^{k} \widetilde{u}|^{2}dw \leq \sum_{j=0}^{k-1} \frac{C_{k}}{(\eta-\xi)^{2k-2j}}\int_{A(r, \eta) } |\nabla^{j}\widetilde{u}|^{2}dw\end{equation} for all $1\leq k\leq m,$ $R/2<r<R$ and $0<\xi<\min\{R-r, r\}.$ Indeed, from \eqref{eq: cd1014} and \eqref{eq: cd1017}, \eqref{eq: cd1016} follows immediately. For $k=m,$ the inequality \eqref{eq: cd1017} is precisely \eqref{eq: cd1014}. Hence, we only need to prove that if \eqref{eq: cd1017} holds for some $k+1<m,$ then it also holds for $k.$

Consider a sequence $\{\rho_{j}\}$ satisfying $\xi:=\rho_{0}<\rho_{1}<...<\eta,$ which will be fixed momentarily. For this sequence, we set
$A_{j}=A(r, \rho_{j}),$ $\delta_{j}:=\rho_{j+1}-\rho_{j},$ and $\widetilde{A}_{j}:=A(r, \rho_{j}+\frac{\delta_{j}}{2}).$ Thus $A_{j}\subset \widetilde{A}_{j} \subset A_{j+1}.$ We also choose a nonnegative, smooth function $\phi_{j},$ supported in $\widetilde{A}_{j}$ and identically 1 on $A_{j},$ satisfying $\|\nabla \phi_{j}\|_{\infty}\leq \frac{C}{\delta_{j}}$ and $\|\nabla^{2} \phi_{j}\|_{\infty}\leq \frac{C}{\delta_{j}^{2}}$ for some absolute constant $C.$ Clearly, for all $j\geq 0,$ $$\int_{A_{j}}|\nabla ^{k}\widetilde{u} |^{2}dw\leq  \int_{\tilde{A}_{j}}|\nabla (\phi_{j} \nabla^{k-1}\widetilde{u}) |^{2}dw.$$ The following key interpolation inequality was proved in \cite{Hy}: for all $f\in H^{2}(w),$ \begin{equation}\label{eq: cd1019} \|\nabla f\|_{L^{2}(w)}^{2}\leq C[w]_{A_{2}} \|\nabla^{2} f\|_{L^{2}(w)}\| f\|_{L^{2}(w)}. \end{equation} By \eqref{eq: cd1019}, we have \begin{equation*}\begin{aligned} \int_{A_{j}}|\nabla ^{k}\widetilde{u} |^{2}dw&\leq C[w]_{A_{2}}^{1/2} \l(\int_{\widetilde{A}_{j}}|\nabla^{2} (\phi_{j} \nabla^{k-1}\widetilde{u}) |^{2}dw\r)^{1/2}\l(\int_{\widetilde{A}_{j}}|\phi_{j} \nabla^{k-1}\widetilde{u} |^{2}dw\r)^{1/2}\\[4pt]  
&\leq C[w]_{A_{2}}^{1/2} \l(\int_{\widetilde{A}_{j}}|\nabla^{k+1}\widetilde{u} |^{2}+\frac{1}{\delta_{j}^{2}}|\nabla^{k}\widetilde{u} |^{2}+\frac{1}{\delta_{j}^{4}}|\nabla^{k-1}\widetilde{u} |^{2}dw\r)^{\frac{1}{2}}\l(\int_{\widetilde{A}_{j}}| \nabla^{k-1}\widetilde{u} |^{2}dw\r)^{\frac{1}{2}}.\end{aligned}\end{equation*} An application of \eqref{eq: cd1017} to control $|\nabla^{k+1}\widetilde{u} |^{2}$ leads to $$\int_{A_{j}}|\nabla ^{k}\widetilde{u} |^{2}dw\leq C[w]_{A_{2}}^{1/2} \l(\sum_{i=0}^{k}\frac{C_{k}}{\delta_{j}^{2k+2-2i}}\int_{A_{j+1}}|\nabla^{i}\widetilde{u} |^{2}dw\r)^{\frac{1}{2}}\l(\int_{\widetilde{A}_{j}}| \nabla^{k-1}\widetilde{u} |^{2}dw\r)^{\frac{1}{2}}.$$ This, by Young's inequality, further implies $$\int_{A_{j}}|\nabla ^{k}\widetilde{u} |^{2}dw\leq \frac{1}{2} \sum_{i=0}^{k}\frac{1}{\delta_{j}^{2k-2i}}\int_{A_{j+1}}|\nabla^{i}\widetilde{u} |^{2}dw+\frac{C_{k}[w]_{A_{2}}}{\delta_{j}^{2}}\int_{\widetilde{A}_{j}}| \nabla^{k-1}\widetilde{u} |^{2}dw.$$ We separate the term for $i=k$ from the sum. This, together with $[w]_{A_{2}}\geq 1,$ yields that $$\int_{A_{j}}|\nabla ^{k}\widetilde{u} |^{2}dw\leq C_{k}[w]_{A_{2}} \sum_{i=0}^{k-1}\frac{1}{\delta_{j}^{2k-2i}}\int_{A_{j+1}}|\nabla^{i}\widetilde{u} |^{2}dw+\frac{1}{2}\int_{A_{j+1}}| \nabla^{k}\widetilde{u} |^{2}dw.$$ Then, using an iteration argument, we arrive at  \begin{equation*}\begin{aligned} \int_{A_{0}}|\nabla ^{k}\widetilde{u} |^{2}dw&\leq \sum_{j=0}^{\infty}2^{-(j-1)}\l(C_{k}[w]_{A_{2}} \sum_{i=0}^{k-1}\frac{1}{\delta_{j}^{2k-2i}}\int_{A_{j+1}}|\nabla^{i}\widetilde{u} |^{2}dw\r)\\[4pt]  
&\leq C_{k}[w]_{A_{2}} \sum_{i=0}^{k-1}\l(\sum_{j=0}^{\infty}2^{-(j-1)}\frac{1}{\delta_{j}^{2k-2i}}\r) \int_{A_{\infty}}|\nabla^{i}\widetilde{u}|^{2}dw.\end{aligned}\end{equation*} Let $0<\tau<1,$ and set $\rho_{0}=\xi$ with $$\rho_{j}:=\xi+(\eta-\xi)(1-\tau)\sum_{i=1}^{j}\tau^{i}\quad \mbox{for}\; j\geq 1.$$ Then $\lim_{j\to \infty}\rho_{j}=\eta.$ We therefore obtain $$\int_{A_{0}}|\nabla ^{k}\widetilde{u} |^{2}dw\leq C_{k, \tau}[w]_{A_{2}} \sum_{i=0}^{k-1}\l(\sum_{j=0}^{\infty}\frac{1}{(2\tau^{2k-2i})^{j}}\frac{1}{(\eta-\xi)^{2k-2i}}\r) \int_{A_{\infty}}|\nabla^{i}\widetilde{u}|^{2}dw.$$ Choosing $\tau$ such that $2\tau^{2k}>1$ and $\tau<1$ proves \eqref{eq: cd1017}. In particular, \eqref{eq: cd1015} is a direct consequence of \eqref{eq: cd1016} and \eqref{eq: cd1019}.

\hfill$\Box$

\section*{Availability of data and material}
 Not applicable.
 
 \section*{Competing interests}
 The author declares that they have no competing interests.


\begin{thebibliography}{10}























\bibitem{A} P. Auscher, \emph{ On necessary and sufficient conditions for $L^{p}$-estimates of Riesz transforms associated to elliptic operators on $\rz$ and related estimates, } Mem. Amer. Math. Soc. 871, American Mathematical Society, Providence, RI, 2007.

\bibitem{A1} P. Auscher, \emph{On $L^{p}$ estimates of square roots of second roder elliptic operators on $\rz$, }  Publ. Mat. 48 (2004), 159-186


\bibitem{ACMP}  P. Auscher, L. Chen, J. M. Martell, C. Prisuelos-Arribas, \emph{The regularity problem for degenerate elliptic operators in weighted spaces,} Rev. Mat. Iberoam. 39 (2023), no. 2, 563-610.




\bibitem{AHLT} P. Auscher, S. Hofmann, M. Lacey, A. McIntosh, P. Tchamitchian, \emph{The solution of the
Kato square root problem for second order elliptic operators on $\rz,$} Ann. of Math. (2) 156 (2002), no. 2,
633-654.



\bibitem{AHMT1} P. Auscher, S. Hofmann, A. McIntosh, P. Tchamitchian, \emph{The Kato square root problem for higher order elliptic operators and systems on $\rz$ } J. Evol. Equ. 1(4), 361-385 (2001).

\bibitem{AM}  P. Auscher, J. M. Martell, \emph{Weighted norm inequalities, off-diagonal estimates and elliptic
operators, II: Off-diagonal estimates on spaces of homogeneous type,} J. Evol. Equ. 7:2 (2007), 265-316.

\bibitem{AM1} P. Auscher, J. M. Martell, \emph{Weighted norm inequalities, off-diagonal estimates and elliptic
operators, III: Harmonic analysis of elliptic operators,} J. Funct. Anal. 241:2 (2006), 703-746.

\bibitem{AM3}  P. Auscher and J. M. Martell, \emph{ Weighted norm inequalities, off-diagonal estimates and elliptic
operators, IV: Riesz transforms on manifolds and weights }, Math. Z. 260:3 (2008), 527-539.





\bibitem{AM2} P. Auscher, J. M. Martell, \emph{Weighted norm inequalities, off-diagonal estimates and elliptic
operators, I: General operator theory and weights,} Adv. Math. 212:1 (2007), 225-276. 

\bibitem{AQ} P. Auscher, M.  Qafsaoui, \emph{Equivalence between regularity theorems and heat kernel estimates
for higher order elliptic operators and systems under divergence form,} J. Funct. Anal. 177 (2000), no. 2, 310-364.


\bibitem{ARR} P. Auscher, A. Ros\'{e}n, D. Rule, \emph{ Boundary value problems for degenerate elliptic equations and systems} Ann. Sci. \'{E}c. Norm. Sup\'{e}r. (4), 48(4): 951-1000, 2015.












\bibitem{AT} P. Auscher, P. Tchamitchian,  \emph{Square root problem for divergence operators and related topics,} Ast\'{e}risque 249 (1998), Societe Mathematique de France. 

\bibitem{B1} A. Barton, \emph{Gradient estimates and the fundamental solution for higher-order elliptic systems with rough coefficients,} Manuscripta Math. 151 (2016), no. 3-4, 375-418. 

\bibitem{B2} A. Barton, S. Hofmann, S. Mayboroda, \emph{ Dirichlet and Neumann boundary values of solutions to higher order elliptic equations, }Ann. Inst. Fourier (Grenoble) 69 (2019), no. 4, 1627-1678.








\bibitem{BB} A. Bj\"{o}rn, J. Bj\"{o}rn, \emph{Nonlinear potential theory on metric spaces, } EMS Tracts in Mathematics 17, European Mathematical Society, Z\"{u}rich.







\bibitem{CMPR} M. E. Cejas, C. Mosquera, C. P\'{e}rez, E. Rela, \emph{Self-improving Poincar\'{e}-Sobolev type functionals in product spaces,} J. Anal. Math. 149 (2023), no. 1, 1-48. 


\bibitem{C11} L. Chen, J. M. Martell, C. Prisuelos-Arribas, \emph{ Conical square functions for degenerate elliptic operators,} Adv. Cal. Var. 13(1):75-113, 2020.

\bibitem{C12} L. Chen, J. M. Martell, C. Prisuelos-Arribas, \emph{The regularity problem for uniformly elliptic operators in weighted spaces,} Potential Anal. 58 (2023), no. 3, 409-439. 


\bibitem{Chua} S. Chua, \emph{Extension theorems on weighted Sobolev spaces,} Indiana Univ. Math. J. 41 (1992), no. 4, 1027-1076.


\bibitem{CR} D. Cruz-Uribe, C. Rios, \emph{The solution of the Kato problem for degenerate elliptic operators with
Gaussian bounds,} Trans. Amer. Math. Soc. 364 (2012), no. 7, 3449-3478.

\bibitem{CR1} D. Cruz-Uribe, C. Rios, \emph{The Kato problem for operators with weighted ellipticity,} Trans. Amer.
Math. Soc. 367 (2015), no. 7, 4727-4756.

\bibitem{CR2} D. Cruz-Uribe, C. Rios, \emph{Gaussian bounds for degenerate parabolic equations,} Journal of Functional Analysis 255 (2008) 283-312.


\bibitem{DMR} D. Cruz-Uribe, J. Martell, C. Rios. \emph{On the Kato problem and extensions for degenerate elliptic
operators,} Anal. PDE 11 (2018), no. 3, 609-660.


\bibitem{D0} E. B. Davies, \emph{Uniformly elliptic operators with measurable coefficients,} J. Funct. Anal. 132 (1995), no. 1, 141-169. 

\bibitem{D1} E. B. Davies, \emph{Limits on $L^{p}$ regularity of selfadjoint elliptic operators,} J. Differential Equations 135, no. 1, (1997), 83-102. 





\bibitem{FJK} E. B. Fabes, D. Jerison, C. Kenig,  \emph{The Wiener test for degenerate elliptic equations,} Ann. Inst. Fourier (Grenoble), 32(3):151-182, 1982

\bibitem{FJK1}  E. B. Fabes, D. Jerison, C. Kenig, \emph{Boundary behavior of solutions to degenerate elliptic equations,} In Conference on harmonic analysis in honor of Antoni Zygmund, Vol. I, II (Chicago, Ill., 1981), Wadsworth Math. Ser. Wadsworth, Belmont, CA, 1983, pp. 577-589.

\bibitem{FJK2} Fabes, E. B., Kenig, C. E., and Serapioni, R. P. \emph{The local regularity of solutions of degenerate elliptic equations,} Comm. Partial Differential Equations, 7(1): 77-116, 1982.











\bibitem{G1} L. Grafakos, \emph{Modern Fourier Analysis,} Grad. Texts in Math., vol. 250, Springer, 2009.




\bibitem{GW} C. E. Gutierrez, R. L. Wheeden, \emph{Sobolev interpolation inequalities with weights,} Trans Amer Math Soc, 1991, 323: 263-281.






\bibitem{H} M. Haase,  \emph{The functional calculus for sectorial operators,} Oper. Theory Adv. Appl., vol. 169, Birkh\"{a}user, Basel, 2006.




\bibitem{Hy} T. P. Hyt\"{o}nen,  \emph{The sharp weighted bound for general Calder\'{o}n-Zygmund operators, } Ann. of Math. (2) 175 (2012), no. 3, 1473-1506.






\bibitem{G2} R. Johnson, C. J. Neugebauer, \emph{Change of variable results for $A_{p}$ and reverse H\"{o}lder $\mbox{RH}_{r}-$classes,} Trans. Amer. Math. Soc. 328 (2) (1991) 639-666.






\bibitem{TK} T. Kato, \emph{ Perturbation theory for linear operators,} Die Grundlehren der mathematischen Wissenschaften, Band 132, Springer-Verlag New York, Inc., New York, 1966.











































\bibitem{Ma} Vladimir G. Maz'ja, \emph{Sobolev spaces,} Springer-Verlag, New York, 1985.









\bibitem{M} A. McIntosh, \emph{Operators which have an $H^{\infty}$ functional calculus. In Miniconference on operator theory and partial differential equations,} Proc. Centre Math. Anal. Austral. Nat. Univ., vol. 14, Austral. Nat. Univ., Canberra, (1986), 210-231.


\bibitem{NM}  N. Miller, \emph{Weighted Sobolev spaces and pseudodifferential operators with smooth symbols,} Trans. Amer. Math.
Soc. 269:1 (1982), 91-109.





\bibitem{CP1} C. P\'{e}rez and E. Rela, \emph{Degenerate Poincar\'{e}-Sobolev inequalities,} Trans. Amer. Math. Soc. 372, no. 9, 6087-6133 (2019)





















	



	







\bibitem{YZ} D. Yang,  J. Zhang \emph{ Weighted $L^{p}$ estimates of Kato square roots associated to degenerate elliptic operators,} Publ. Mat. 61: 2 (2017), 395-444.





\bibitem{GZ} G. Zhang, \emph{The Kato problem for weighted elliptic and parabolic systems of higher order,} arXiv: 2503. 01636.















	
\end{thebibliography}
\end{document}